%
%
%
\newif\ifsect\newif\iffinal
\secttrue\finalfalse
\def\thm #1: #2{\medbreak\noindent{\bf #1:}\if(#2\thmp\else\thmn#2\fi}
\def\thmp #1) { (#1)\thmn{}}
\def\thmn#1#2\par{\enspace{\sl #1#2}\par
        \ifdim\lastskip<\medskipamount \removelastskip\penalty 55\medskip\fi}
\def\square{{\msam\char"03}}
\def\qedn{\thinspace\null\nobreak\hfill\square\par\medbreak}
\def\pf{\ifdim\lastskip<\smallskipamount \removelastskip\smallskip\fi
        \noindent{\sl Proof\/}:\enspace}
\def\itm#1{\item{\rm #1}\ignorespaces}

\def\bar#1{\overline{#1}}
\input boxedeps
\input boxedeps.cfg
\HideDisplacementBoxes
\def\Figuraeps #1 (#2){\message{Figura #1}
	\midinsert      
	\centerline{\BoxedEPSF{#2.eps}}
	\bigskip
	\centerline{\bf Figure~#1}
	\endinsert}
\def\Figurascaledeps #1 (#2)(#3){\message{Figura #1}
	\midinsert      
	\centerline{\BoxedEPSF{#2.eps scaled #3}}
	\bigskip
	\centerline{\bf Figure~#1}
	\endinsert}
\newcount\parano
\newcount\eqnumbo
\newcount\thmno
\newcount\versiono
\newcount\remno
\newbox\notaautore
\def\neweqt#1$${\xdef #1{(\number\parano.\number\eqnumbo)}
    \eqno #1$$
    \global \advance \eqnumbo by 1}
\def\newrem #1\par{\global \advance \remno by 1
    \medbreak
{\bf Remark \the\parano.\the\remno:}\enspace #1\par
\ifdim\lastskip<\medskipamount \removelastskip\penalty 55\medskip\fi}
\def\newthmt#1 #2: #3{\xdef #2{\number\parano.\number\thmno}
    \global \advance \thmno by 1
    \medbreak\noindent
    {\bf #1 #2:}\if(#3\thmp\else\thmn#3\fi}
\def\neweqf#1$${\xdef #1{(\number\eqnumbo)}
    \eqno #1$$
    \global \advance \eqnumbo by 1}
\def\newthmf#1 #2: #3{\xdef #2{\number\thmno}
    \global \advance \thmno by 1
    \medbreak\noindent
    {\bf #1 #2:}\if(#3\thmp\else\thmn#3\fi}
\def\forclose#1{\hfil\llap{$#1$}\hfilneg}
\def\newforclose#1{
	\ifsect\xdef #1{(\number\parano.\number\eqnumbo)}\else
	\xdef #1{(\number\eqnumbo)}\fi
	\hfil\llap{$#1$}\hfilneg
	\global \advance \eqnumbo by 1
	\iffinal\else\rsimb#1\fi}
\def\forevery#1#2$${\displaylines{\let\eqno=\forclose
        \let\neweq=\newforclose\hfilneg\rlap{$\qquad\quad\forall#1$}\hfil#2\cr}$$}
\def\noNota #1\par{}
\def\today{\ifcase\month\or
   January\or February\or March\or April\or May\or June\or July\or August\or
   September\or October\or November\or December\fi
   \space\number\year}
\def\inizia{\ifsect\let\neweq=\neweqt\else\let\neweq=\neweqf\fi
\ifsect\let\newthm=\newthmt\else\let\newthm=\newthmf\fi}
\def\bititolo{\empty}
\gdef\begin #1 #2\par{\xdef\titolo{#2}
\ifsect\let\neweq=\neweqt\else\let\neweq=\neweqf\fi
\ifsect\let\newthm=\newthmt\else\let\newthm=\newthmf\fi
\iffinal\let\Nota=\noNota\fi
\centerline{\titlefont\titolo}
\if\bititolo\empty\else\medskip\centerline{\titlefont\bititolo}
\xdef\titolo{\titolo\ \bititolo}\fi
\bigskip
\centerline{\bigfont
\autore \ifvoid\notaautore\else\footnote{${}^1$}{\unhbox\notaautore}\fi}
\bigskip\if\istituto!\centerline{\it\today}\else
\centerline{\istituto}
\centerline{\indirizzo}
\centerline{E-mail: \email}
\medskip
\centerline{#1\ \anno}\fi
\bigskip\bigskip
\ifsect\else\global\thmno=1\global\eqnumbo=1\fi}
\def\anno{2007}
\def\raggedleft{\leftskip2cm plus1fill \spaceskip.3333em \xspaceskip.5em
\parindent=0pt\relax}

\font\titlefont=cmssbx10 scaled \magstep1
\font\bigfont=cmr12
\font\eightrm=cmr8

\font\sc=cmcsc10
\font\bbr=msbm10
\font\sbbr=msbm7
\font\ssbbr=msbm5
\font\msam=msam10

\font\bfm=cmmib10

\def\ca #1{{\cal #1}}

\nopagenumbers
\binoppenalty=10000
\relpenalty=10000
\newfam\amsfam
\textfont\amsfam=\bbr \scriptfont\amsfam=\sbbr \scriptscriptfont\amsfam=\ssbbr
\newfam\boldifam
\textfont\boldifam=\bfm
\let\de=\partial
\let\eps=\varepsilon
\let\phe=\varphi

\def\Aut{\mathop{\rm Aut}\nolimits}
\def\Hom{\mathop{\rm Hom}\nolimits}

\def\Sing{\mathop{\rm Sing}\nolimits}
\def\End{\mathop{\rm End}\nolimits}

\def\Res{\mathop{\rm Res}\nolimits}
\def\Re{\mathop{\rm Re}\nolimits}
\def\Im{\mathop{\rm Im}\nolimits}

\def\id{\mathop{\rm id}\nolimits}

\def\ord{\mathop{\rm ord}\nolimits}

\mathchardef\void="083F
\mathchardef\ellb="0960
\mathchardef\taub="091C
\def\C{{\mathchoice{\hbox{\bbr C}}{\hbox{\bbr C}}{\hbox{\sbbr C}}
{\hbox{\sbbr C}}}}
\def\R{{\mathchoice{\hbox{\bbr R}}{\hbox{\bbr R}}{\hbox{\sbbr R}}
{\hbox{\sbbr R}}}}
\def\N{{\mathchoice{\hbox{\bbr N}}{\hbox{\bbr N}}{\hbox{\sbbr N}}
{\hbox{\sbbr N}}}}
\def\P{{\mathchoice{\hbox{\bbr P}}{\hbox{\bbr P}}{\hbox{\sbbr P}}
{\hbox{\sbbr P}}}}

\def\Z{{\mathchoice{\hbox{\bbr Z}}{\hbox{\bbr Z}}{\hbox{\sbbr Z}}
{\hbox{\sbbr Z}}}}

\newcount\notitle
\notitle=1
\headline={\ifodd\pageno\rhead\else\lhead\fi}
\def\rhead{\ifnum\pageno=\notitle\iffinal\hfill\else\hfill\tt Version
\the\versiono; \the\day/\the\month/\the\year\fi\else\hfill\eightrm\titolo\hfill
\folio\fi}
\def\lhead{\ifnum\pageno=\notitle\hfill\else\eightrm\folio\hfill\autore\hfill
\fi}
\newbox\bibliobox
\def\setref #1{\setbox\bibliobox=\hbox{[#1]\enspace}
    \parindent=\wd\bibliobox}
\def\biblap#1{\noindent\hang\rlap{[#1]\enspace}\indent\ignorespaces}
\def\art#1 #2: #3! #4! #5 #6 #7-#8 \par{\biblap{#1}#2: {\sl #3\/}.
    #4 {\bf #5} (#6)\if.#7\else, \hbox{#7--#8}\fi.\par\smallskip}
\def\book#1 #2: #3! #4 \par{\biblap{#1}#2: {\bf #3.} #4.\par\smallskip}
\def\coll#1 #2: #3! #4! #5 \par{\biblap{#1}#2: {\sl #3\/}. In {\bf #4,}
#5.\par\smallskip}
\def\pre#1 #2: #3! #4! #5 \par{\biblap{#1}#2: {\sl #3\/}. #4, #5.\par\smallskip}
\def\End{\mathop{\rm End}\nolimits}
\def\Bittersweet{} 
\def\raggedleft{\leftskip2cm plus1fill \spaceskip.3333em \xspaceskip.5em
\parindent=0pt\relax}
\def\Nota #1\par{\medbreak\begingroup\Bittersweet\raggedleft
#1\par\endgroup
\ifdim\lastskip<\medskipamount \removelastskip\penalty 55\medskip\fi}
\newcount\defno
\def\smallsect #1. #2\par{\bigbreak\noindent{\bf #1.}\enspace{\bf #2}\par
    \global\parano=#1\global\eqnumbo=1\global\thmno=1\global\defno=0\global\remno=0
    \nobreak\smallskip\nobreak\noindent\message{#2}}
\def\newdef #1\par{\global \advance \defno by 1
    \medbreak
{\bf Definition \the\parano.\the\defno:}\enspace #1\par
\ifdim\lastskip<\medskipamount \removelastskip\penalty 55\medskip\fi}
\def\autore{Marco Abate\footnote{${}^1$}{\eightrm Support by the Centro di Ricerca Matematica ``Ennio de Giorgi" (Pisa, Italy) and by the Institut Mittag-Leffler (Djursholm, Sweden) is gratefully
acknowledged.}, Francesca Tovena\footnote{}{\eightit Mathematics Subject Classification 2010: \eightrm Primary 32H50; 32S65, 37F75, 34M35. \eightit Keywords: \eightrm local holomorphic dynamics; homogeneous vector fields; meromorphic connections; geodesics; normal forms.}}
\def\istituto{!}
\def\indirizzo{via Buonarroti 2, 56127 Pisa, Italy}
\def\email{abate@dm.unipi.it}
\def\bititolo{and holomorphic homogeneous vector fields}
\finaltrue
\versiono=11
\month=12
\begin {November} Poincar\'e-Bendixson theorems for meromorphic connections

{\narrower
{\sc Abstract.} We first study the dynamics of the geodesic flow of
a meromorphic connection on a Riemann surface, and prove a 
Poincar\'e-Bendixson theorem describing
recurrence properties and $\omega$-limit sets of geodesics for a
meromorphic connection on $\P^1(\C)$. We then show how to associate
to a homogeneous vector field $Q$ in $\C^n$ a rank 1 singular holomorphic 
foliation~$\ca F$ of~$\P^{n-1}(\C)$ and a (partial) meromorphic connection~$\nabla^o$ 
along~$\ca F$ so that integral curves of~$Q$ are described by the geodesic flow of $\nabla^o$ along the leaves of~$\ca F$, which are Riemann surfaces. The combination of these results yields powerful tools for a detailed study of the dynamics of homogeneous vector fields. For instance, in dimension two we obtain 
a description of recurrence properties of integral curves
of $Q$, and of the behavior of the geodesic flow
in a neighbourhood of a singularity, classifying the possible
singularities both from a formal 
point of view and (for generic singularities) from a holomorphic point of view. We also get examples of unexpected new phenomena,
we put in a coherent context scattered results previously known, and we obtain (as far as we know for the first time) a
complete description of the dynamics in a full neighbourhood of the origin for a substantial class of 2-dimensional holomorphic maps tangent to the identity. Finally, as an example of application of our methods we study in detail the dynamics of quadratic homogeneous vector fields in $\C^2$. 

}

\smallsect 0. Introduction

In this paper we shall study the dynamics of two apparently unrelated objects:
geodesics for meromorphic connections on a Riemann surfaces, and integral curves
of homogeneous vector fields in~$\C^n$. 

Meromorphic connections on Riemann surfaces have been well-studied, particularly from a 
algebraic geometrical point of view (see, e.g., [M]);  however, to our knowledge the
dynamical properties of the {\sl real} geodesic curves associated to a meromorphic connection
have not been investigated before.

One of the main results of this paper is (as the title suggests) a type of Poincar\'e-Bendixson 
theorem describing the recurrence properties and $\omega$-limit sets of geodesics for a 
meromorphic connection on $\P^1(\C)$. The classical Poincar\'e-Bendixson theorem (see, e.g., [HK, Theorem~14.1.1]) deals with integral curves of vector fields defined on open subsets
of the sphere~$S^2$ (notice that, as a differentiable manifold, $\P^1(\C)$ is diffeomorphic
to~$S^2$): recurrent integral curves are necessarily periodic, and 
the $\omega$-limit set of an integral curve either contains singular points or is a periodic integral curve. In our Poincar\'e-Bendixson theorem (see Theorem~4.6) the poles of the
connection replace the singular points of the vector field:\def\autore{Marco Abate, Francesca Tovena}

\newthm Theorem \ztre: Let $\sigma\colon[0,\eps_0)\to S$ be a maximal
geodesic for a meromorphic connection~$\nabla^o$ on $\P^1(\C)$, where
$S=\P^1(\C)\setminus\{p_0,\ldots,p_r\}$ and $p_0,\ldots,p_r$ are the
poles of~$\nabla^o$. Then either
\smallskip
\itm{(i)} $\sigma(t)$ tends to a pole of~$\nabla^o$ as $t\to\eps_0$; or
\itm{(ii)} $\sigma$ is closed, and then surrounds poles $p_1,\ldots,p_g$
 with $\sum\limits_{j=1}^g\Re\Res_{p_j}(\nabla^o)=-1$; or
\item{\rm (iii)} the $\omega$-limit set of $\sigma$ in $\P^1(\C)$ is given by the support of a closed geodesic surrounding poles
$p_1,\ldots,p_g$ with $\sum\limits_{j=1}^g\Re\Res_{p_j}(\nabla^o)=-1$; or
\item{\rm (iv)} the $\omega$-limit set of $\sigma$ in $\P^1(\C)$ is a simple cycle of saddle connections (see below) surrounding
poles $p_1,\ldots,p_g$ with $\sum\limits_{j=1}^g\Re\Res_{p_j}(\nabla^o)=-1$; or
\itm{(v)} $\sigma$ intersects itself infinitely many times, and in this
case every simple loop of $\sigma$ surrounds a set of poles whose
sum of residues has real part belonging to $(-3/2,-1)\cup(-1,-1/2)$.
\smallskip
\noindent In particular, a recurrent geodesic either intersects itself
infinitely many times or is closed.

Here, a {\sl saddle connection} is a geodesic connecting two (not
necessarily distinct) poles
of $\nabla^o$; a {\sl simple cycle of saddle connections} is a Jordan
curve composed of saddle connections. Furthermore, the residue of the connection 
at a pole~$p$ is defined as the residue at~$p$ of the meromorphic 1-form representing the 
connection with respect to any holomorphic local coordinate in~$p$. 

We have examples (see Examples~6.1, 8.1 and~8.2) of cases (i), (ii), (iii) and (v); we do not know 
yet whether case (iv) can actually be realized (but we are able to exclude it
in several situations; see Remark~8.3).
Notice furthermore that (see, e.g., [IY, Theorem~I
\negthinspace I\negthinspace I.17.33]) the only
limitation on the residues of a meromorphic connection on~$\P^1(\C)$
is that their sum should be~$-2$; more precisely,  given any finite set of
pairs $\{(p_1,r_1),\ldots,(p_g,r_g)\}\subset\P^1(\C)\times\C$ with $r_1+\cdots+r_g=-2$
there exists a meromorphic connection~$\nabla^o$ with poles $\{p_1,\ldots,p_g\}$ 
and $\Res_{p_j}(\nabla^o)=r_j$ for $j=1,\ldots,g$. Since cases (ii)--(v) of Theorem~\ztre\
impose additional conditions on the residues (in particular, in the last case
the condition should be satisfied by each of the infinitely many simple loops of the
geodesic; see also Proposition~4.7), it follows that maximal geodesics of a meromorphic connection often display no 
recurrence phenomena at all, being simply saddle connections.  

It is important to notice that in general a meromorphic connection~$\nabla^o$ on a Riemann surface 
is not the Chern connection of a Hermitian metric (unless all residues are real:
see Proposition~1.2 and Corollary~3.7). Furthermore, even when it is, the associated Hermitian 
metric is never complete (except in trivial cases: see Corollary~2.5); so the behavior of our 
geodesics is subtly different from the behavior of the usual geodesics in Riemannian geometry
(for instance, we can have closed geodesics which are not periodic: geodesics for meromorphic
connections are not necessarily of constant speed). Nevertheless,
we shall be able to associate to $\nabla^o$ a conformal family of local Hermitian flat
metrics that shall be very useful.

The proof of Theorem~\ztre\ depends on three main ingredients, developed in the
first four sections of this paper. The first ingredient is a detailed study
of geodesics for holomorphic connections on a simply connected Riemann surface (the 
universal cover of the complement of the poles); since a
holomorphic connection is necessarily flat, it turns out that its geodesics behave locally as 
Euclidean segments. The second ingredient is a Gauss-Bonnet theorem, 
used to control the monodromy of the geodesics. Indeed, it turns out that, even
though we do not have a global metric, we can still use a Gauss-Bonnet formula, relating
the external angles of geodesics polygons to the residues of the poles inside the 
polygon (see Theorem~4.1), which in turn control the monodromy of the connection
(Proposition~3.6) and the speed of the geodesics (Lemma~4.4). Finally, the proof is 
completed by a delicate argument introducing a sort of Poincar\'e return map on a transversal
defined at a point in the $\omega$-limit set of the geodesic.

The second half of the paper is devoted to the dynamics of homogeneous vector fields. A homogeneous vector field in $\C^n$ is a vector field of the form
$$
Q=Q^1{\de\over\de z^1}+\cdots+Q^n{\de\over\de z^n}\;,
\neweq\eqzuno
$$
where $Q^1,\ldots,Q^n$ are homogeneous polynomials of degree $\nu+1\ge 2$. The complex 
foliation generated by $Q$ is not that difficult to study (see, e.g., Theorem~6.2); but here we are 
mostly interested in the dynamics of the {\it real} integral curves of~$Q$. The reason is
that we arrived to this problem because we wanted to study the dynamics of holomorphic 
maps tangent to the identity (that is, of germs of holomorphic self-maps of $\C^n$
fixing the origin and with differential at the origin equal to the identity), and homogeneous vector fields provide good examples of those. Indeed, it is
easy to see that the time-1 map of a vector field of the form~\eqzuno\ has a homogeneous 
expansion of the form
$$
f(z)=z+Q_{\nu+1}(z)+\cdots\;,
$$
where $Q_{\nu+1}=(Q^1,\ldots,Q^n)$, and thus $f$ is tangent to the identity. In dimension one, 
the classical Leau-Fatou theorem (see, e.g., [Mi] and [A4]) gives a complete description of the 
dynamics of a holomorphic function tangent to the identity in a full neighbourhood 
of the origin. Using this, in 1978 Camacho [C] (see also [Sh]) proved that, from a topological
point of view, time-1 maps of homogeneous vector fields provide a complete list of models 
for the dynamics:

\newthm Theorem \zCamacho: (Camacho, 1978 [C]) Let
$
f(z)=z+a_{\nu+1}z^{\nu+1}+\cdots
$
with $a_{\nu+1}\ne 0$, be a germ of holomorphic function tangent to the identity. Then $f$ is
locally topologically conjugated to the time-1 map of the homogeneous vector field 
$$
Q=z^{\nu+1}{\de\over\de z}\;.
$$

In dimension greater than one nothing of the sort is (as yet) known. More precisely, understanding
the (topological) dynamics of holomorphic germs tangent to the identity in a full
neighbourhood of the origin is one of the main open problems in local dynamics of several complex variables. There are versions of the Leau-Fatou flower theorem in several  variables,
obtained by \'Ecalle (see [\'E1--4]) and Hakim (see [H1, 2]) in any dimension but for generic
germs, and for all germs in dimension 2 by the first author (see [A2]).
But these theorems mostly give the 
existence of 1-dimensional invariant sets only, and are quite far from providing a
description of the dynamics in a full neighbourhood of the origin. In fact, as far as we know, 
before the present paper such a description was available for a handful of examples only.

On the other hand, along the lines of Camacho's Theorem~\zCamacho, it is conjectured that in any dimension a {\it generic} (e.g., with only
non-degenerate characteristic directions; see below) germ tangent to
the identity is locally topologically conjugated to the time-1 map of a homogeneous
vector field. To build such a conjugation, one usually needs a precise description of
the dynamics of the model; so we decided to study in detail the dynamics of 
time-1 maps of homogeneous vector fields. The orbit of a point~$p$ under the action of such
a map is contained in the {\it real} integral curve issuing from~$p$, and so we are led to
the study of the dynamics of real integral curves.

To state the (somewhat unexpected) link between integral curves of a homogeneous vector
field and geodesics of a meromorphic connection we need a few definitions. 
Let $Q$ be a homogeneous vector field of the form~\eqzuno, and denote by $[\cdot]\colon\C^n\setminus\{O\}\to\P^{n-1}(\C)$ the canonical projection.
Following \'Ecalle (see [\'E1--4]) and Hakim (see [H1, 2]) we say that a direction
$[v]\in\P^1(\C)$ is a {\sl characteristic direction} of $Q$ if the line $L_v=\C v$
is $Q$-invariant (and in that case we say that $L_v$ is a {\sl characteristic leaf\/}).
If $Q$ is identically zero along~$L_v$ we say that $[v]$ is {\sl degenerate;}
otherwise, it is {\sl non-degenerate.} It turns out (see, e.g., [AT1]) that either
$Q$ has ${1\over\nu}[(\nu+1)^n-1]$ characteristic directions, counting multiplicities, or infinitely many
directions are characteristic. In particular, if all directions
are characteristic we shall say that $Q$ is {\sl dicritical}. When $n=2$,
it turns out that either $Q$ is dicritical or it has $\nu+2$ characteristic
directions, counting multiplicities.
 
The dynamical meaning of characteristic directions is expressed by the following fact
(see [H1]): if $\gamma$ is an integral curve of $Q$ converging to the origin tangentially
to some direction~$[v]\in\P^{n-1}(\C)$ then $[v]$ is a characteristic
direction. Notice however that (as noted by Rivi [R]) 
there might exist orbits converging to the origin without being tangent to any direction;
see Example~6.1 (and Corollary 8.5, giving an explanation of this phenomenon).

The dynamics inside a characteristic leaf is 1-dimensional, and easy to study (see Lemma~5.4).
So we are interested in the dynamics outside characteristic leaves of (necessarily)
non-dicritic vector fields. Our second main result says that 
integral curves outside characteristic leaves are given by geodesics
for a suitable meromorphic connection on suitable Riemann surfaces, foliating a projective space.
Let $\ca F$ be a
rank 1 singular holomorphic foliation of a complex manifold~$M$. Let~$\Sing(\ca F)$ be the singular 
set of~$\ca F$, and set $M^o=M\setminus\Sing(\ca F)$. A {\sl (partial) meromorphic connection 
along~$\ca F$} is a $\C$-linear map $\nabla^o\colon\ca F|_{M^o}\to\ca F|_{M^o}^*\otimes\ca F|_{M^o}$ satisfying the usual Leibniz condition (see Section~5 for details); roughly speaking, $\nabla^o$ allows to
differentiate sections of~$\ca F$ along directions tangent to the 
foliation. In particular, $\nabla^o$ induces a (classical) meromorphic
connection on each (1-dimensional) leaf of the foliation. 

If $\sigma\colon I\to M^o$ is a curve contained in a leaf of the foliation
(that is, $\sigma'(t)\in\ca F_{\sigma(t)}$ for all $t\in I$), and $\nabla^o$
is a meromorphic connection along~$\ca F$, then we can consider
$\nabla^o_{\sigma'}\sigma'$; we shall say that $\sigma$ is a
{\sl $\nabla^o$-geodesic} if $\nabla^o_{\sigma'}\sigma'\equiv O$. In 
other words, a curve $\sigma$ contained in a leaf~$L$ is a 
$\nabla^o$-geodesic if and only if it is a geodesic for the meromorphic
connection induced on~$L$. 

The link between integral curves and geodesics is then provided by the following 
result (see Theorem~5.3):

\newthm Theorem \zdmez: Let $Q$ be a non-dicritical homogeneous vector field in $\C^n$, and 
let $\hat S_Q$ be the complement in~$\C^n$ of the characteristic leaves of~$Q$.
Then there exists a rank-$1$ singular holomorphic foliation~$\ca F$ 
on~$\P^{n-1}(\C)$, whose singular points are characteristic directions of~$Q$, and a meromorphic connection $\nabla^o$ along~$\ca F$ such that a curve
$\gamma\colon I\to\hat S_Q$ is an integral curve of~$Q$ if and only if its
direction $[\gamma]\colon I\to\P^{n-1}(\C)$ is a $\nabla^o$-geodesic.

This result then reduces the study of integral curves of a homogeneous vector field to the study
of the foliation $\ca F$ and to the study of geodesics for a meromorphic connection on 
Riemann surfaces.

The proof of Theorem~\zdmez\ depends on another ingredient. Let $\pi\colon M\to\C^n$ be the blow up of the origin (see [A1] for a 
description of the blow-up construction adapted to dynamical purposes); the 
exceptional divisor (that is, the preimage of the origin under $\pi$) is canonically
identified with $\P^{n-1}(\C)$. Let $p\colon N^{\otimes\nu}_{\P^{n-1}(\C)}\to\P^{n-1}(\C)$
be the $\nu$-th tensor power of the normal bundle~$N_{\P^{n-1}(\C)}$ of the exceptional
divisor in $M$. Then it is possible to define (see Proposition~5.2) a
$\nu$-to-one holomorphic covering map $\chi_\nu\colon\C^n\setminus\{O\}\to 
N^{\otimes\nu}_{\P^{n-1}(\C)}\setminus\P^{n-1}(\C)$ (where we are identifying $\P^{n-1}(\C)$
with the zero section of $N^{\otimes\nu}_{\P^{n-1}(\C)}$) so that $p\circ\chi_\nu(z)=[z]$
for all $z\in\C^n\setminus\{O\}$.

Usually, the push-forward of a vector field is not a vector field. However, the homogeneity of~$Q$
implies (see Theorem~5.3) that $d\chi_\nu(Q)$ is a holomorphic vector field~$G$ globally
defined on the total space of $N^{\otimes\nu}_{\P^{n-1}(\C)}$; so using~$\chi_\nu$ we
transform the study of integral curves of~$Q$ in the study of integral curves of~$G$.

In [ABT1] we introduced (in a more general setting) a canonical morphism $X_Q\colon
N^{\otimes\nu}_{\P^{n-1}(\C)}\to T\P^{n-1}(\C)$ associated to the (non-dicritical) 
homogeneous vector field~$Q$. The zeroes of~$X_Q$ are exactly the characteristic directions of~$Q$;
so $X_Q$ is an isomorphism outside the characteristic directions,
and thus it defines a rank 1 singular holomorphic foliation~$\ca F$
of~$\P^{n-1}(\C)$. Furthermore, again in [ABT1]
we showed how to use~$Q$ to define (in an essentially unique way; see [ABT2]) a partial
holomorphic connection $\nabla$ along~$\ca F$ on $N^{\otimes\nu}_{\P^{n-1}(\C)}$.  

In Section~5 we shall describe this construction in our context, adding a few new ideas. In particular, we shall remark that using~$X_Q$ we can push~$\nabla$ to~$T\P^{n-1}(\C)$ 
obtaining a meromorphic connection~$\nabla^o$ along~$\ca F$ in the
sense mentioned above, and
we shall show (Proposition~5.1) that a curve $\sigma$ in $\P^{n-1}(\C)$ is a 
$\nabla^o$-geodesic if and only if the image of $\sigma$ is contained 
in a leaf of $\ca F$ and the curve $X_Q^{-1}(\sigma')$ in the total space of 
$N^{\otimes\nu}_{\P^{n-1}(\C)}$ is an integral curve of the global vector field~$G=d\chi_\nu(Q)$. So everything fits together, and
we can use meromorphic connections on Riemann surfaces to study the dynamics of homogeneous vector fields.

To exemplify the strength of this method, in the rest of the paper we
specialize to the case of $n=2$, where the foliation~$\ca F$ has only one leaf, the complement of the characteristic directions.
Then, as a corollary of Theorem~\ztre\ we immediately get (see Theorem~6.3) a description of recurrent 
integral curves:

\newthm Theorem \ztremez: Let $Q$ be a homogeneous holomorphic vector
field on $\C^2$, and let $\gamma$
be a recurrent maximal integral curve of $Q$. Then either $\gamma$ is periodic or 
$[\gamma]$ intersect itself infinitely many times.

More can be said along this line (see, e.g., Section~9); but to fully
understand the dynamics we need to know what happens to integral curves nearby the
characteristic lines, that is to the geodesics nearby the poles. So in
Sections~7 and~8 we turn to a detailed study of the geodesic field $G$ and its singularities in dimension~2, showing that we must distinguish between three types of singularities: apparent, Fuchsian (which is the generic case)
and irregular. We shall be able to give a complete formal classification of 
all cases (see Proposition~7.1 and Theorem~8.1), and a complete holomorphic classification of the first
two cases (Proposition~7.1 and Theorem~8.3); in particular, it is worthwhile to remark that in the Fuchsian case 
resonances appear. More precisely, we shall prove the following (see Theorem~8.3):

\newthm Theorem \zqua: Let $p_0\in\P^1(\C)$ be a Fuchsian pole of~$G$, that is
assume that in local coordinates
$(U_\alpha,z_\alpha)$ centered at~$p_0$, denoting by $v_\alpha$ the induced coordinate along the fibers
of~$N^{\otimes\nu}_{\P^1(\C)}$, we can write
$$
G=z_\alpha^{\mu}(a_0+a_1 z_\alpha+\cdots)v_\alpha{\de\over\de z_\alpha}-
z_\alpha^{\mu-1}(b_0+b_1 z_\alpha+\cdots)v_\alpha^2{\de\over\de v_\alpha}\;,
$$
with $\mu\ge 1$ and $a_0$, $b_0\ne 0$. Put $\rho=b_0/a_0=\Res_{p_0}(\nabla)$. Then we can
find a chart $(U,z)$ centered in~$p_0$ in which $G$ is given by
$$
z^{\mu-1}\left(z v{\de\over\de z}-\rho v^2{\de\over\de v}\right)
$$
if $\mu-1-\rho\notin\N^*$, or by
$$
z^{\mu-1}\left(zv{\de\over\de z}-\rho(1+a z^n)  v^2{\de\over\de v}\right)
$$
for a suitable $a\in\C$ if $n=\mu-1-\rho\in\N^*$.

It is easy to check that non-degenerate characteristic directions with non-zero residue
are Fuchsian (with $\mu=1$). Using Theorem~\zqua\ we then get a complete description of the 
dynamics in a neighbourhood of Fuchsian characteristic directions (see Proposition~8.4 and Corollary~8.5); we also have a complete description of the dynamics
in a neighbourhood of apparent singularities (see Corollaries~7.2 and~7.3). Putting this
together with our Poincar\'e-Bendixson theorems 
we get a
complete description of the dynamics for a substantial class of 2-dimensional homogenous vector fields
(see Corollary~8.6):

\newthm Theorem \zsei: Let $Q$ be a non-dicritical 
homogeneous vector field on $\C^2$. Assume
that all characteristic directions of $Q$ are non-degenerate with non-zero residue. Assume moreover that for no set of characteristic
directions the real part of the sum of the residues of $\nabla^o$ is equal to~$-1$.\hfil\break\indent Let $\gamma\colon[0,\eps_0)\to\C^2$ be a maximal integral curve of $Q$. Then:
\smallskip
\itm{(a)} If $\gamma(0)$ belongs to a characteristic leaf $L_{v_0}$, then the whole image of $\gamma$ is contained 
in~$L_{v_0}$. Moreover, either $\gamma(t)\to O$ (and this happens for
an open Zariski dense set of initial conditions in~$L_{v_0}$), or $\|\gamma(t)\|\to+\infty$.
\itm{(b)} If $\gamma(0)$ does not belong to a characteristic leaf, then either
\itemitem{\rm (i)} $\gamma$ converges to the origin tangentially to a characteristic 
direction $[v_0]$ whose residue has negative real part; or
\itemitem{\rm (ii)} $\|\gamma(t)\|\to+\infty$ tangentially to a characteristic
direction $[v_0]$ whose residue has positive real part; or
\itemitem{\rm(iii)} $[\gamma]\colon[0,\eps_0)\to\P^1(\C)$ intersects itself
infinitely many times. 
\smallskip
Furthermore, if {\rm (iii)} never occurs then {\rm (i)} holds for an open Zariski dense set
of initial conditions.

In particular, if the residues of $Q$ does not satisfy the condition
corresponding to the one described in Theorem~\ztre.(v) --- and we already remarked that
this is often the case --- then 
case (b.iii) of Theorem~\zsei\ cannot occur, and 
we get a description of the dynamics of the time-1 map of $Q$ in a full neighbourhood of the origin; as mentioned before,
as far as we know this is the first such description for a substantial class
of maps tangent to the identity. We also remark that we do have a description of the dynamics
even when the real part of the sum of the residues might be~$-1$, or for some
classes of degenerate characteristic directions, and so the scope of our results is larger
than Theorem~\zsei; see Section~8 for details and Section~9 for examples.

In our opinion, this approach not only offers effective tools for studying the dynamics
of homogeneous vector fields (and thus hopefully of maps tangent to the identity), 
but it also gives a better understanding of what is going on. For instance, Hakim's theorem [H2]
on the existence of parabolic basins in this context is explained by the fact that
non-degenerate characteristic directions whose residue has negative real part are attractors (see 
Corollary~8.5, that actually extends Hakim's theorem to some degenerate characteristic
directions in dimension two). Or, Rivi's [R] example of orbits going to the origin without being tangent to
any direction turns out to be related to the existence of characteristic directions with
purely imaginary residue (see again Corollary~8.5). We are also able to give examples of unexpected phenomena.
In dimension one, the Leau-Fatou flower theorem implies that a map tangent to the
identity has no small cycles: there is a neighbourhood
of the origin containing no periodic points beside the origin itself. So it was expected
that even in several complex variables maps tangent to the identity could not have small 
cycles. Surprisingly, this turns out to be false: in Corollaries~7.3 and~8.5 we shall give
examples having periodic points of arbitrarily high period accumulating at the origin (see also
Corollary~6.4).

This paper is organized as follows. In Section~1 we shall introduce the local metrics
and the global (metric) foliation associated to a holomorphic connection on a line
bundle over a Riemann surface~$S$. Along the way,
we shall characterize the holomorphic connections which are
the Chern connection of a Hermitian metric on a line bundle over $S$.
In Section~2 we shall study in depth the geodesic flow of a holomorphic connection
over a simply connected Riemann surface. In Section~3 we shall use the monodromy representation of a holomorphic connection over a multiply connected Riemann
surface to study the geodesic flow there. In Section~4 we shall consider 
meromorphic connections on $\P^1(\C)$, and use the previous
material to prove Theorem~\ztre. In Section 5 we shall clarify the relations between
maps tangent to the identity, homogeneous vector fields and meromorphic connections
in any dimension,
proving in particular Theorem~\zdmez. In Section~6 we shall specialize to dimension~2,
proving Theorem~\ztremez. In Section 7 we shall begin the study of the geodesic flow nearby the singularities in dimension two;
the formal and holomorphic classifications, as well as their dynamical consequences and the proofs of Theorem~\zqua\ and~\zsei,
are contained in Section~8. Finally, in Section~9 we shall discuss 
in detail 2-dimensional quadratic homogenous vector fields.

Developing the ideas leading to this paper has been a long process, carried
out not only in our home institutions but (mostly) in several other places.
We would like to thank the Department of Mathematics of Niigata, Kyoto and Barcelona Universities,
the IMPA (Rio de Janeiro, Brazil) and, in particular, the Institut
Mittag-Leffler (Djursholm, Sweden) for their warm hospitality and 
productive environment. 

\smallsect 1. The metric, horizontal and geodesic foliations

Let us begin recalling a few standard facts about holomorphic connections
on line bundles over Riemann surfaces (see, e.g., [K]).

\newdef Let $E$ be a complex line bundle on a Riemann surface $S$. A
{\sl holomorphic connection} on $E$ is a $\C$-linear map 
$\nabla\colon\ca E\to\Omega^1_S\otimes\ca E$ satisfying the Leibniz rule
$$
\nabla(se)=ds\otimes e+s\nabla e
$$
for all $s\in\ca O_S$ and $e\in\ca E$, where $\ca E$ denotes the sheaf of 
germs of holomorphic sections of~$E$, while $\ca O_S$ is the structure sheaf of~$S$ and $\Omega^1_S$ is the sheaf of holomorphic
1-forms on~$S$. A {\sl horizontal section} of~$\nabla$ is a section $e\in\ca E$ such that $\nabla e\equiv O$.

Let $\{(U_\alpha,z_\alpha,e_\alpha)\}$
be an atlas of~$S$ trivializing $E$, where $(U_\alpha,z_\alpha)$ are local charts
of $S$ and $e_\alpha$ is a holomorphic generator of $E|_{U_\alpha}$. Over
$U_\alpha$, a
holomorphic connection~$\nabla$ is represented by a holomorphic 1-form
$\eta_\alpha\in\Omega^1_S(U_\alpha)$ such that
$$
\nabla e_\alpha=\eta_\alpha\otimes e_\alpha\;.
$$
If $\{\xi_{\alpha\beta}\}$ is the cocycle representing the
cohomology class $\xi\in H^1(S,\ca O^*)$ of~$E$, 
over~$U_\alpha\cap U_\beta$ we have
$$
e_\beta=e_\alpha \xi_{\alpha\beta}
$$
and
$$
\eta_\beta=\eta_\alpha+{1\over\xi_{\alpha\beta}}\,\de\xi_{\alpha\beta}\;.
\neweq\equbeta
$$
Recalling the short exact sequence of
sheaves
$$
O\to\C^\ast\to {\cal O}^\ast\mathop{\longrightarrow}\limits^{\de \log}  \Omega_S^1\to O\;,
$$
we see that equality \equbeta\ shows that the existence of a holomorphic
connection
$\nabla$ is equivalent to the vanishing of the image of $\xi$ under the map 
$\de\log\colon H^1(S, \ca O^*)\to H^1(S, \Omega_S^1)$ induced on cohomology. So, the class
$\xi$ is the image of a class $\hat \xi\in H^1(S, \C^\ast)$: we now recall how to find a
representative
$\hat{\xi}_{\alpha\beta}$ of $\hat{\xi}$.
Up to shrinking the $U_\alpha$'s, we can find
holomorphic functions~$K_\alpha\in\ca O(K_\alpha)$ such that $\eta_\alpha= \de K_\alpha$ on~$U_\alpha$. Then \equbeta\ implies that on
$U_\alpha\cap U_\beta$
$$
\hat{\xi}_{\alpha\beta}= {\exp(K_\alpha) \over \exp(K_\beta)}\, \xi_{\alpha\beta}
\neweq\eqhatxi
$$
is a complex non zero constant defining a cocycle representing $\xi$. We furthermore remark that 
$$
\nabla \bigl(\exp(-K_\alpha) e_\alpha\bigr)\equiv O\;,
\neweq\equhor
$$
that is $\exp(-K_\alpha) e_\alpha$ is a horizontal section on~$U_\alpha$.

\newdef The homomorphism $\rho\colon \pi_1(S)\to \C^\ast$ corresponding
to the class $\hat\xi$ under the canonical isomorphism $H^1(S, \C^*)\cong\Hom\bigl(H_1(S,\Z),\C^*\bigr)=\Hom\bigl(\pi_1(S),\C^*\bigr)$ is
the {\sl monodromy representation} of the holomorphic connection $\nabla$.
We shall say that $\nabla$ has {\sl real periods} if the image of $\rho$
is contained in~$S^1$, that is if $\hat\xi$ is the image of a class
in $H^1(S,S^1)$ under the natural inclusion $S^1\hookrightarrow\C^*$.

In Proposition~3.5 we shall explicitely compute the monodromy representation
when $S\subseteq\C$, explaining the rationale behind the terminology. 

Now, it is well-known
that to a Hermitian metric $g$ on a complex vector
bundle over a complex manifold~$M$ can be associated a unique 
(1,0)-connection $\nabla$ (not necessarily holomorphic) such that $\nabla g\equiv O$, the {\sl Chern connection} of $g$. We would like 
to study the converse problem: given a holomorphic connection $\nabla$,
does there exist a Hermitian metric $g$ so that $\nabla g\equiv O$?

\newdef Let $E$ be a complex line bundle on a Riemann surface $S$, and 
$\nabla\colon\ca E\to\Omega^1_S\otimes\ca E$ a holomorphic
connection on~$E$. We say that a Hermitian metric $g$ on $E$ is {\sl adapted} 
to~$\nabla$ if $\nabla g\equiv O$, that is if
$$
X\bigl(g(R,T)\bigr)=g(\nabla_X R,T)+g(R,\nabla_{\bar X}T)
$$
and
$$
\bar{X}\bigl(g(R,T)\bigr)=g(\nabla_{\bar X} R,T)+g(R,\nabla_X T)
$$
for all (not necessarily holomorphic) sections $R$,~$T$ of~$E$, and all vector 
fields $X$ on~$S$. 


Let us see what this condition means in local coordinates. With respect to an atlas $\{(U_\alpha,z_\alpha,
e_\alpha)\}$ trivializing $E$, a hermitian metric $g$ on $E$ is locally represented by a positive
$C^\infty$ function $n_\alpha\in C^\infty(U_\alpha,\R^+)$ given by
$$
n_\alpha=g(e_\alpha,e_\alpha)\;.
$$
Then it is easy to see that $g$ is adapted to $\nabla$ over $U_\alpha$ if and only if
$$
\de n_\alpha=n_\alpha\eta_\alpha\;,
\neweq\eqHuno
$$
where $\eta_\alpha$ is the holomorphic 1-form representing $\nabla$.
Our first result shows how to solve this equation:

\newthm Proposition \uuno: 
Let $E$ be a complex line bundle on a Riemann surface $S$, and $\nabla\colon\ca E\to
\Omega^1_S\otimes\ca E$ a holomorphic connection on~$E$. Let
$(U_\alpha,z_\alpha,e_\alpha)$ be a local chart trivializing $E$, and 
define $\eta_\alpha\in\Omega^1_S
(U_\alpha)$ by setting $\nabla e_\alpha=\eta_\alpha\otimes e_\alpha$.  Assume that we
have a holomorphic primitive~$K_\alpha$ of~$\eta_\alpha$ on $U_\alpha$. Then
$$
n_\alpha=\exp(2\Re K_\alpha)=\exp(K_\alpha+\bar{K_\alpha})
\neweq\eququa
$$
is a positive solution of $\eqHuno$. Conversely, if $n_\alpha$ is a 
positive solution of $\eqHuno$ then for any $z_0\in U_\alpha$ and simply connected
neighbourhood~$U\subseteq U_\alpha$ of~$z_0$ there is a holomorphic primitive $K_\alpha\in\ca
O(U)$ of~$\eta_\alpha$ over $U$ such that
$n_\alpha=\exp(2\Re K_\alpha)$ in $U$. Furthermore, $K_\alpha$ is unique up to a purely immaginary
additive constant. Finally, two (local) solutions of $\eqHuno$ differ (locally) by a
positive multiplicative constant. 

\pf If $n_\alpha$ is given by \eququa\ then it is always real-valued and positive, and we
have
$$
\de n_\alpha=\exp(K_\alpha)\exp(\bar{K_\alpha})\de K_\alpha=n_\alpha\eta_\alpha\;;
$$
thus $n_\alpha$ satisfies \eqHuno.

Conversely, assume that $n_\alpha$ is a positive solution of \eqHuno\ on~$U_\alpha$. Then
$\log n_\alpha$ is well-defined and satisfies $\de\log n_\alpha=\eta_\alpha$. 
In particular, $\bar\de\de\log n_\alpha=\bar\de\eta_\alpha=0$, because $\eta_\alpha$
is holomorphic. So $\log n_\alpha$ is harmonic; hence for each
$z_0\in U_\alpha$ and simply connected neighbourhood $U\subseteq U_\alpha$ of~$z_0$
we can find $\tilde n_\alpha\colon U\to\R$
harmonic conjugate to ${1\over 2} \log n_\alpha$,
so that 
$$
K_\alpha={1\over 2}\log n_\alpha+i \tilde n_\alpha
$$
is holomorphic in~$U$. But then
$$
\de(i\tilde n_\alpha)=\bar{-\bar\de(i\tilde n_\alpha)}=\bar{\bar\de({\textstyle{1\over2}}
\log n_\alpha)}=\de({\textstyle{1\over2}}\log n_\alpha)=
{\textstyle{1\over2}}\eta_\alpha\;;
$$
therefore $K_\alpha$ is a holomorphic primitive of~$\eta_\alpha$ such that $n_\alpha=\exp(2\Re K_\alpha)$ over $U$, as required.

If $K_1$ is another holomorphic function defined in $U$ so
that $n_\alpha=\exp(2\Re K_1)$, we find 
$$
\exp\bigl(2\Re(K_\alpha-K_1)\bigr)\equiv 1\;,
$$
and hence
$\Re K_\alpha-\Re K_1\equiv 0$. The holomorphicity of $K_\alpha$ and $K_1$ then forces~$K_\alpha-K_1\equiv ia$ for some $a\in\R$.

Finally, if $n_1$ and $n_2$ are two positive solutions of \eqHuno,
we have $\de\log(n_1/n_2)\equiv O$. Since $\log(n_1/n_2)$ is real-valued, this forces
$d\log(n_1/n_2)\equiv 0$, and thus $n_1/n_2$ is (locally) constant.
\qedn

\newrem It is well known that a holomorphic 1-form $k\,dz$ defined in an open set~$U\subseteq S$ has a (necessarily holomorphic) primitive if and only if
$$
\int_\gamma k\,dz=0
$$
for all closed loops~$\gamma$ in~$U$. 
However, the obstructions to the existence of $n_\alpha$ are slightly weaker,
because we just need the exponential of the real part of a primitive. To be more explicit,
let assume that $U=\Delta^*$ is a pointed disk, and use the Laurent expansion to write
$$
k(z)=k^*(z)+{\rho\over z}\;,
$$
where $\rho$ is the residue of $k$ at the origin. Then $k^*$ has vanishing residue at the
origin, and thus it admits a primitive $K^*$ on~$U$. Locally, a primitive 
of~$\rho/z$ is
of the form $\rho\log z$; therefore setting
$$
K(z)=K^*(z)+\rho\log z
$$
we have a locally defined (multivalued) primitive of~$k$. We are interested in the
exponential of the real part of~$K$; since
$$
\Re K(z)=\Re K^*(z)+(\Re\rho)\log|z|-(\Im\rho)\arg(z)\;,
$$
we find
$$
\exp\bigl(2\Re K(z)\bigr)=\exp\bigl(2\Re
K^*(z)\bigr)|z|^{\Re\rho}\exp\bigl(-(\Im\rho)\arg(z)\bigr)\;.
$$
But $\arg(z)$ is defined up to an integer multiple of~$2\pi$; therefore, the
indeterminacy of $\exp(2\Re K)$ is a multiplicative factor of the form
$e^{-2\pi h\Im\rho}$, with $h\in\Z$. In particular, if the residue $\rho$ is real,
then we get a well-defined solution of~$\eqHuno$ in the whole $\Delta^*$. We 
shall prove (see Proposition~1.2 and Corollary~3.5) that, roughly speaking, this will be the only obstruction to the existence
of a global metric adapted to $\nabla$.

\newrem The Gaussian curvature of a local metric~$g$ adapted to~$\nabla$ is identically 
zero. Indeed, $g$ is of the form $\exp(2\Re K)g_0$, where $g_0$ is
the euclidean metric. The Gaussian curvature of a metric of the form~$h g_0$ 
is $-{1\over h}\Delta\log h$. In our case $h=\exp(2\Re K)=|\exp(K)|^2$
is the modulus squared of a holomorphic function; so $\log h$
is harmonic, and hence $\Delta h=0$.

Having solved the local problem, let us see when we get a global hermitian metric
adapted to~$\nabla$. Let $\{(U_\alpha,z_\alpha,e_\alpha)\}$
be an atlas of~$S$ trivializing $E$. Up to shrinking the $U_\alpha$'s, we can take
a holomorphic primitive~$K_\alpha$ of the holomorphic 1-form~$\eta_\alpha$ representing~$\nabla$ on~$U_\alpha$. Taking the logarithm of the modulus
of~\eqhatxi\ on $U_\alpha\cap
U_\beta$ we get
$$
\Re(K_\alpha-K_\beta)+\log|\xi_{\alpha\beta}|=\log|\hat\xi_{\alpha\beta}|\;.
\neweq\equsei
$$
%
%
Hence:

\newthm Proposition \udue: 
Let $E$ be a complex line bundle on a Riemann surface $S$, and $\nabla\colon\ca E\to
\Omega^1_S\otimes\ca E$ a holomorphic connection on~$E$. Then there exists a
Hermitian metric adapted to~$\nabla$ if and only if $\nabla$ has real periods.

\pf Let $g$ be a Hermitian metric on $E$. If
$\{(U_\alpha,z_\alpha,e_\alpha)\}$ is any atlas trivializing~$E$, setting
$n_\alpha=g(e_\alpha,e_\alpha)$ over~$U_\alpha$, we
must have
$$
n_\beta=|\xi_{\alpha\beta}|^2 n_\alpha
\neweq\equcin
$$
over $U_\alpha\cap U_\beta$. If $g$ is adapted to $\nabla$, up to shrinking the 
$U_\alpha$'s, we can assume that
$n_\alpha=\exp(2\Re K_\alpha)$ over~$U_\alpha$, where $K_\alpha$ is a holomorphic
primitive of the form $\eta_\alpha$ representing~$\nabla$ over $U_\alpha$. Then
\equcin\ says that
$$
\Re(K_\beta-K_\alpha)=\log|\xi_{\alpha\beta}|
$$
over $U_\alpha\cap U_\beta$, and hence $\nabla$ has real periods.

Conversely, assume that $\nabla$ has real periods. Then we can find an atlas $\{(U_\alpha,z_\alpha,
e_\alpha)\}$ trivializing $E$, holomorphic primitives $K_\alpha$, and  constants $c_\alpha\in\C^*$ such that $\hat\xi_{\alpha\beta}=(c_\beta/c_\alpha)\tilde\xi_{\alpha\beta}$ with $\tilde\xi_{\alpha\beta}\in S^1$, so that
$$
\Re(K_\beta-K_\alpha)-\log|\xi_{\alpha\beta}|=\log|c_\alpha|-\log|c_\beta|\;.
$$
Then $\tilde K_\alpha=K_\alpha+\log|c_\alpha|$ is a holomorphic primitive of
$\eta_\alpha$ such that \equcin\ is satisfied by $n_\alpha=\exp(2\Re\tilde K_\alpha)$,
and thus setting $g(e_\alpha,e_\alpha)=n_\alpha$ we get a global Hermitian metric adapted
to~$\nabla$. \qedn


Actually, we shall be more interested in the case when does not exist a metric
adapted to $\nabla$. Indeed, the first main result of this section is the following

\newthm Proposition \utre: Let $\nabla$ be a holomorphic connection on a complex line
bundle~$E$ over a Riemann surface~$S$. Then there exists a real rank~$3$ 
non-singular foliation of~$E\setminus S$ (where we are identifying $S$ with the zero section of~$E$)
whose leaves are the level sets of any local or global Hermitian metric on~$E$ adapted
to~$\nabla$. 

\pf Choose an atlas $\{(U_\alpha,z_\alpha,e_\alpha)\}$ of $S$ trivializing~$E$ (with connected
intersections) and such that on each $U_\alpha$ we can find a holomorphic
primitive~$K_\alpha$ of the
holomorphic form~$\eta_\alpha$ representing~$\nabla$. Set $n_\alpha=\exp(2\Re K_\alpha)$, and
define $g_\alpha\colon p^{-1}(U_\alpha)\to\R^+$ by setting
$g_\alpha(v)=n_\alpha\bigl(p(v)\bigr)|v_\alpha|^2$, where $p\colon E\to S$ is the
canonical projection, and $v_\alpha\in\C$ is so that~$v=v_\alpha e_\alpha$. Clearly,
$g_\alpha$ is a submersion out of the zero section, and thus its level sets define a real rank 3
non-singular foliation of~$p^{-1}(U_\alpha)\setminus U_\alpha$; we must show that
the foliation is independent of~$\alpha$. 
But indeed if $z_0\in U_\alpha\cap U_\beta$ and $v\in E_{z_0}$ formula \equsei\ yields
$c_{\alpha\beta}\in\R$ such that
$$
\eqalign{
g_\beta(v)&=n_\beta(z_0)|v_\beta|^2=n_\alpha(z_0)\exp\bigl(2\Re(K_\beta-K_\alpha)(z_0)\bigr)
|v_\beta|^2=|\hat\xi_{\alpha\beta}|^{-2}
n_\alpha(z_0)|\xi_{\alpha\beta}(z_0)|^2|v_\beta|^2\cr
&=|\hat\xi_{\alpha\beta}|^{-2}n_\alpha(z_0)|v_\alpha|^2
=|\hat\xi_{\alpha\beta}|^{-2}g_\alpha(v)\;;
\cr}
$$
therefore $g_\alpha$ and $g_\beta$ differ by a multiplicative constant, and thus they have
the same level sets in $p^{-1}(U_\alpha\cap U_\beta)$.

Finally, if $g$ is a Hermitian metric adapted to~$\nabla$ defined on an open set~$U$,
Proposition~\uuno\ says that on each $U\cap U_\alpha$ the function $n=g(e_\alpha,e_\alpha)$ locally
is a positive multiple of~$n_\alpha$, and hence the level sets of the norm induced by~$g$
coincide with the level sets of~$g_\alpha$.\qedn

\newdef The foliation just defined induced by a holomorphic connection $\nabla$ on a
complex line bundle $E$ over a Riemann surface~$S$ is the {\sl metric foliation}
of~$\nabla$ on~$E$.

\newrem When $g$ is a globally defined Hermitian metric adapted to~$\nabla$, the
leaves of the metric foliation are simply the sets  $\{v\in E\setminus S\mid g(v,v)=
{\rm const.}\}$, and thus they are diffeomorphic to $S^1\times S$
and closed in the total space of~$E$ (zero section included). But when $\nabla$
does not admit an adapted metric the leaves might have a more complicated behavior;
in particular, they can accumulate the zero section (and thus they are not closed in
the total space of $E$); see Theorems~3.2 and~3.3.

The choice of a local chart $(U_\alpha,z_\alpha,e_\alpha)$ trivializing
the line bundle $p\colon E\to S$ yields local coordinates $(z_\alpha,v_\alpha)$ on
the total space of~$E$, and thus a local frame for $TE$. Let us denote by
$\{\de_\alpha,\de/\de v_\alpha\}$ this local frame, where $\de_\alpha$ is the tangent
vector corresponding to the coordinate~$z_\alpha$, and by $\{p^*(d
z_\alpha),dv_\alpha\}$ the dual co-frame. From
$e_\beta=e_\alpha\xi_{\alpha\beta}$ in $U_\alpha\cap U_\beta$ we get
$$
v_\alpha=(\xi_{\alpha\beta}\circ p)v_\beta\;,
\neweq\equv
$$
and
$$
dv_\alpha=v_\beta\,p^*(\de\xi_{\alpha\beta})+(\xi_{\alpha\beta}\circ p)
\,dv_\beta=v_\beta\left({\de\xi_{\alpha\beta}\over\de
z_\beta}\circ p\right)p^*(d z_\beta)+(\xi_{\alpha\beta}\circ p)\,dv_\beta\;.
\neweq\eqdvalpha
$$
Furthermore, we have
$$
p^*(d z_\alpha)=(\psi_{\alpha\beta}\circ p)p^*(d z_\beta)\;,
$$
where $\psi_{\alpha\beta}=\de z_\alpha/\de z_\beta$.
It follows that
$$
{\de\over\de v_\beta}=(\xi_{\alpha\beta}\circ p){\de\over\de v_\alpha}\qquad\hbox{and}
\qquad \de_\beta=(\psi_{\alpha\beta}\circ p)\de_\alpha+v_\alpha
\left({1\over\xi_{\alpha\beta}}
{\de\xi_{\alpha\beta}\over\de z_\beta}\right)\circ p\,{\de\over\de v_\alpha}\;.
$$
In particular, the local sections~$v_\alpha \de/\de v_\alpha$ give a
globally defined section~$R$ of~$TE$, whose integral curves are the real lines
through the origin in each fiber of~$E$. Analogously, the integral curves of
$iR$ are circumferences around the origin in each fiber of~$E$, and gives a rank 1 real
foliation of the total space of~$E$, singular along the zero section. 

Using these notations, it is easy
to see that in a local chart $(U_\alpha,z_\alpha,e_\alpha)$ trivializing $E$ the metric foliation
is generated by the real 1-form
$$
\varpi_\alpha=\Re\bigl(|v_\alpha|^2\,p^*\eta_\alpha+\bar{v_\alpha}\,dv_\alpha\bigr)\;,
\neweq\eqbaromega
$$
where $\eta_\alpha$ is the holomorphic 1-form representing~$\nabla$. In particular, the tangent
space to the foliation is generated by $H_\alpha$, $iH_\alpha$ and $iR$, where
$H_\alpha$ is the local section of~$TE$ defined by
$$
H_\alpha=\de_\alpha-(k_\alpha\circ p) v_\alpha{\de\over\de v_\alpha}\;,
\neweq\equHalpha
$$
with $k_\alpha=\eta_\alpha(\de/\de z_\alpha)$. In particular, it
is clear that the metric foliation is transversal to the fibers of~$E$.

The local fields $H_\alpha$ define a complex rank 1 foliation of the total 
space of~$E$. Indeed, we have
$$
\eqalign{
H_\beta&=\de_\beta-(k_\beta\circ p) v_\beta{\de\over\de v_\beta}\cr
&=
(\psi_{\alpha\beta}\circ p)\de_\alpha+v_\alpha\left({1\over\xi_{\alpha\beta}}
{\de\xi_{\alpha\beta}\over\de z_\beta}\right)\circ p\,{\de\over\de v_\alpha}-\left(\psi_{\alpha\beta}
\left[k_\alpha+{1\over\xi_{\alpha\beta}}{\de\xi_{\alpha\beta}\over\de
z_\alpha}\right]\right)\circ p\,{1\over\xi_{\alpha\beta}\circ p}
v_\alpha(\xi_{\alpha\beta}\circ p){\de\over\de v_\alpha}\cr
&=(\psi_{\alpha\beta}\circ p)H_\alpha\;.
\cr}
\neweq\eqduno
$$
Furthermore, it is easy to check that a local section $s_\alpha$
of $E$ is an integral curve of $H_\alpha$ if and only if $\nabla s_\alpha\equiv O$,
that is if and only if $s_\alpha$ is a horizontal section.

\newdef The complex rank 1 non-singular foliation on $E\setminus S$ induced by the local
fields $H_\alpha$ is the {\sl horizontal foliation} of the holomorphic
connection~$\nabla$. Clearly, the leaves of the horizontal foliation
are transversal to the fibers of~$E$, and are contained in the leaves of the metric
foliation. 

It is also easy to describe this foliation using a global holomorphic 1-form on the
total space of $E$. Indeed, let
as always $\{(U_\alpha,z_\alpha,e_\alpha)\}$ be an atlas trivializing~$E$, and let
$(z_\alpha,v_\alpha)$ denote the corresponding local coordinates on~$E|_{U_\alpha}$. Denote by $p\colon E\to S$ the projection.
Then \eqdvalpha\ and \equbeta\ yield
$$
\eqalign{
p^*\eta_\beta+{1\over
v_\beta}\,dv_\beta&=p^*\eta_\alpha+{1\over\xi_{\alpha\beta}\circ p}\,
p^*\de\xi_{\alpha\beta}
+{\xi_{\alpha\beta}\circ p\over v_\alpha}\left[{1\over\xi_{\alpha\beta}\circ p}\,dv_\alpha-
{v_\alpha\over\xi_{\alpha\beta}^2\circ p}\,p^*\de\xi_{\alpha\beta}\right]\;,\cr
&=p^*\eta_\alpha+{1\over v_\alpha}\,dv_\alpha\;.\cr}
$$
Therefore setting
$$
\omega=p^*\eta_\alpha+{1\over v_\alpha}\,dv_\alpha
$$
on $E|_{U_\alpha}\setminus U_\alpha$, we get a global holomorphic
1-form on~$E\setminus S$. Furthermore, $\omega(H_\alpha)\equiv 0$,
and so $\omega$ induces the horizontal foliation, as claimed.

\newrem We clearly have 
$$
\varpi_\alpha=|v_\alpha|^2\Re\omega\;,
$$
where $\varpi$ is given by \eqbaromega.

\newrem We can extend the horizontal foliation to a non-singular foliation of the whole
total space of~$E$ just by adding the zero section as a new leaf. Indeed, we have
$$
v_\beta p^*\eta_\beta+dv_\beta=v_\beta\,\omega={1\over\xi_{\alpha\beta}
\circ p}v_\alpha\,\omega
={1\over\xi_{\alpha\beta}\circ p}\bigl[v_\alpha p^*\eta_\alpha+dv_\alpha\bigr]\;;
$$
therefore the local forms $v_\alpha p^*\eta_\alpha+dv_\alpha$ define a complex rank 1 non-singular foliation on~$E$ which coincides with the horizontal foliation off the zero
section.

Later on we shall need local parametrizations for the leaves of the horizontal
foliation. We need a holomorphic map~$\phe\colon V\to E$ defined on some
open set~$V\subseteq\C$ and such that 
$$
\phe'=H_\alpha\circ\phe\;.
$$
Writing in local coordinates $\phe(\zeta)=\bigl(z_\alpha(\zeta),v_\alpha(\zeta)\bigr)$ we
see that we need
$$
z_\alpha'\equiv 1\qquad\hbox{and}\qquad v'_\alpha=-k_\alpha v_\alpha\;.
$$
Hence
$$
z_\alpha(\zeta)=\zeta+c_0\;,\qquad
v_\alpha(\zeta)=c_1\exp\bigl(-K_\alpha(\zeta+c_0)\bigr)\;,
\neweq\eqparhf
$$
where $c_0$,~$c_1\in\C$ and $K_\alpha$ is a holomorphic primitive of~$\eta_\alpha$
on~$V+c_0$; compare with \equhor. 

Since the local fields $H_\alpha$ do not glue together in the intersections, they do
not define a {\it real} rank 1 foliation of~$E\setminus S$. As discussed in the
introduction, in the cases we shall be interested in we shall have another ingredient 
available, allowing us to introduce a real rank 1 foliation of $E\setminus S$.

Assume that we have an isomorphism $X\colon E\to TS$.
The idea is to use $X$ to define $\nabla$-geodesics, and then consider the geodesic flow.

\newdef Let $E$ be a line bundle on a Riemann surface $S$, and assume we have a
holomorphic connection $\nabla$ on $E$ and a isomorphism $X\colon E\to TS$.
We say that a smooth curve $\sigma\colon I\to S$, where $I\subseteq\R$ is an interval, is a {\sl geodesic} (with respect to $\nabla$ and
$X$) if $\nabla_{\sigma'}X^{-1}(\sigma')\equiv O$. If $\sigma$ is a geodesic, then
$X^{-1}(\sigma')$ is a curve in the total space of~$E$; we shall momentarily show that it is an 
integral curve of a vector field on $E$. 

\newrem The reason we are explicitely using the isomorphism $X$ instead of just
considering geodesics for a holomorphic connection on $TS$ is that in the
applications we have in mind the line bundle $E$ will be the restriction to $S$ of a
line bundle $\hat E$ defined on a larger Riemann surface $\hat S\supset S$. We 
shall have a morphism $X\colon\hat E\to T\hat S$, but this will be an isomorphism
only over $S$. Furthermore, we shall be interested in the behavior of geodesics
in $\hat S$, and of the geodesic flow in the total space of $\hat E$; and to study those
it will be important to work in $E$ using $X$ instead of working in $TS$. However,
in the next three sections we shall deal with $E=TS$ and $X=\id$ only.

If $(U_\alpha,z_\alpha,e_\alpha)$ is a local chart trivializing~$E$, then there is a holomorphic 
function $X_\alpha\in\ca O^*(U_\alpha)$ such that
$$
X(e_\alpha)=X_\alpha{\de\over\de z_\alpha}\;,
$$ 
and it is easy to check that changing coordinates $X_\alpha$ changes according to the rule
$$
X_\beta={\xi_{\alpha\beta}\over\psi_{\alpha\beta}}X_\alpha\;.
\neweq\eqtuno
$$
Then we have 

\newthm Proposition \uqua:
Let $\nabla$ be a holomorphic connection on a complex line bundle~$p\colon E\to S$ over a Riemann
surface~$S$, and $X\colon E\to TS$ an isomorphism. Let $\{(U_\alpha,z_\alpha,e_\alpha)\}$ 
be an atlas trivializing~$E$. Then:
\smallskip
\itm{(i)} setting 
$$
G|_{p^{-1}(U_\alpha)}=(X_\alpha\circ p) v_\alpha H_\alpha=(X_\alpha
\circ p) v_\alpha\de_\alpha-(X_\alpha
k_\alpha)\circ p\,(v_\alpha)^2{\de\over\de v_\alpha}
$$ 
we define a global
holomorphic section~$G$ of~$TE$, vanishing only on the zero section;
\itm{(ii)} a curve $\sigma\colon I\to S$ is a geodesic if and only if $X^{-1}(\sigma')$ is
an integral curve of~$G$.

\pf (i) follows immediately from \eqtuno, \eqduno\ and \equv. Denoting by $z_\alpha(t)$
the expression of the curve $\sigma\colon I\to S$ in the local chart~$(U_\alpha,z_\alpha)$,
it is easy to see that $\sigma$ is a geodesic if and only if
$$
\left({z'_\alpha\over X_\alpha}\right)'+k_\alpha X_\alpha \left({z'_\alpha\over
X_\alpha}
\right)^{\!2}\equiv 0\;.
\neweq\eqgde
$$
On the other hand, a curve $t\mapsto\bigl(z_\alpha(t),v_\alpha(t)\bigr)$ is an integral curve of~$G$
if and only if
$$
\cases{z'_\alpha=X_\alpha(z_\alpha) v_\alpha\;,\cr
\noalign{\smallskip}
v'_\alpha=-k_\alpha(z_\alpha) X_\alpha(z_\alpha) v_\alpha^2\;.
\cr}
\neweq\eqsist
$$
Since $X^{-1}(\sigma')$ is expressed in local coordinates by
$\bigl(z_\alpha,z'_\alpha/X_\alpha(z_\alpha)\bigr)$, assertion (ii) follows.\qedn

\newdef The global holomorphic field $G$ is the {\sl geodesic field} associated
to~$\nabla$ and~$X$. The rank~1 non-singular real foliation of $E\setminus S$ given by the
integral curves of $G$ is the {\sl geodesic foliation} associated to~$\nabla$
and~$X$. Clearly, the leaves of the geodesic foliation are contained in the leaves of the horizontal foliation.

\newrem Since $G$ is a global field, the leaves of the geodesic foliation, 
being integral curves of $G$, are equipped with
a canonical parametrization. In principle, we can get such a parametrization by quadratures and taking inverses. Indeed, let $t\mapsto\bigl(z_\alpha(t),v_\alpha(t)\bigr)$ be a local integral curve of~$G$. By \eqparhf\ we must have $v_\alpha(t)=c_1\exp\bigl(-K_\alpha\bigl(z_\alpha(t)\bigr)\bigr)$;
hence the first equation in \eqsist\ yields
$$
{\exp\bigl(K_\alpha(z_\alpha)\bigr)\over X_\alpha(z_\alpha)}\,z'_\alpha=c_1\;.
$$
If $F_\alpha$ is a primitive of $\exp(K_\alpha)/X_\alpha$ we then get
$F_\alpha\bigl(z_\alpha(t)\bigr)=c_1t+c_2$; since $F'_\alpha\ne 0$ always
we finally get
$$
z_\alpha(t)=F^{-1}_\alpha(c_1t+c_2)\;.
$$ 
We shall use this procedure in the last section of this paper.

\newrem The leaves of the 
geodesic foliation are contained in the leaves of the horizontal foliation, and thus in the
leaves of the metric foliation. Furthermore, they are
transversal to the fibers of~$E$. So we have cut the total space of $E$ off the zero
section in three real foliations, of real rank 3, 2 and 1 respectively, one inside the other,
and all transversal to the fibers of $E$.

\newrem Clearly, the field $iG$ defines another real rank 1 non-singular foliation of $E\setminus S$; but we shall not use it in this paper.

The main goal of this paper will be the study of the dynamics of the geodesic foliation,
and then the application of our results to the study of the dynamics of homogeneous
vector fields in $\C^n$. Along the way we shall also get a few (usually easier) results
on the dynamics of the metric and horizontal foliations. 

\smallsect 2. Simply connected Riemann surfaces

Let $\tilde S$ be a simply connected Riemann surface, and assume that we have a
holomorphic connection $\tilde\nabla$ on~$T\tilde S$. In particular
(see, e.g., [IY, Theorem~I\negthinspace I\negthinspace I.17.33]), $\tilde S$
cannot be $\P^1(\C)$, and so $\tilde S$ is biholomorphic either to $\C$ or to the
unit disk~$\Delta$. In both cases, $T\tilde S=\tilde S\times\C$, and we have a
global coordinate $z$  on~$\tilde S$. We would like to study the metric, horizontal
and geodesic foliations associated to $\tilde\nabla$ on $T\tilde S$.

We use $\de/\de z$ as global section of $T\tilde S$, giving an explicit isomorphism
between $T\tilde S$ and $\tilde S\times\C$.
Let $\tilde\eta=\tilde k\,dz$ be the (global) holomorphic 1-form associated to~$\tilde\nabla$,
and $\tilde K\colon\tilde S\to\C$ a (global) holomorphic primitive of $\tilde\eta$.
By Proposition~\uuno, the function $\tilde g\colon\tilde S\times\C\to\C$ given by
$$
\tilde g(z;v)=\exp\bigl(2\Re\tilde K(z)\bigr)\,|v|^2
\neweq\eqdzero
$$
is the norm squared of a Hermitian metric (that we shall also denote by $\tilde g$) adapted to 
$\tilde\nabla$. In particular, the leaves of the metric foliation are just given by the level sets 
of $\tilde g$:
$$
\exp\bigl(2\Re\tilde K(z)\bigr)\,|v|^2= {\rm const.}\in\R^+\;.
$$
Furthermore, we have a similar description for the leaves of the horizontal foliation;
indeed, \eqparhf\ says that they can be expressed by
$$
\exp\bigl(\tilde K(z)\bigr)\,v={\rm const.}\in\C^*\;,
\neweq\eqhfp
$$
that is as the level sets of the holomorphic function 
$(z;v)\mapsto\exp\bigl(\tilde K(z)\bigr)\,v$.

As a consequence, the dynamics of both the metric and the horizontal foliations over a
simply connected Riemann surface is pretty trivial. In particular, each leaf of the
horizontal foliation intersect each fiber in exactly one point, and we have an explicit
biholomorphism between $\tilde S$ and the leaf through the point $(z_0;
v_0)$ given by
$$
z\mapsto \bigl(z;\exp\bigl(-\tilde K(z)\bigr)\exp\bigl(\tilde K(z_0)\bigr)v_0\bigr)\;.
$$
Analogously, each leaf of the metric foliation intersect each fiber in exactly one 
circumference, and we have an explicit
diffeomorphism between $S^1\times\tilde S$ and the leaf through the point $(z_0;
v_0)$ given by
$$
(e^{2\pi i\theta},z)\mapsto \bigl(z;\exp\bigl(-\Re\tilde K(z)\bigr)\exp\bigl(\Re\tilde 
K(z_0)\bigr)|v_0|e^{2\pi i\theta}\bigr)\;.
$$

The rest of this section is devoted to the study of the geodesic foliation.
As a, somewhat unexpected (at least by us), consequence,  we shall see 
(Corollary~2.5) that the metric
$\tilde g$ adapted to $\tilde\nabla$ is {\it never} complete 
(unless $\tilde S=\C$ and $\tilde\nabla$ is trivial), preventing the use of
standard theorems like the Hopf-Rinow theorem, even though in this case
our geodesics are the usual Riemannian geodesics of the Riemannian metric $\Re\tilde g$.

Our first result is the following

\newthm Proposition \dzero: Let $\tilde\nabla$ be a holomorphic connection 
on $T\tilde S=\tilde S\times\C$, where $\tilde S\cong\C$ or $\Delta$ is a 
simply connected Riemann surface. Let $\tilde\eta$ be the holomorphic $1$-form associated
to~$\tilde\nabla$, and $\tilde K\colon\tilde S\to\C$ a holomorphic primitive of
$\tilde\eta$. Finally, let $J\colon\tilde S\to\C$ be a holomorphic primitive of
$\exp(\tilde K)$.
Then $J\colon\tilde S\to\C$ is a local
isometry, where $\tilde S$ is endowed with the metric~$\tilde g$ 
adapted to~$\tilde\nabla$ corresponding to~$\tilde K$, and $\C$ is endowed with the euclidean metric. 

\pf First of all, since $\tilde S$ is simply connected, $J$ exists. Now,
by \eqdzero, the $\tilde g$-length of~$v\in T_zS$ is 
$$
\exp\bigl(\Re\tilde K(z)\bigr)|v|=|J'(z)v|\;,
$$
and hence $J$ is a local isometry.\qedn

In particular, and this will be important in the sequel, $J$ sends the
geodesic segments contained in any open set $U\subseteq\tilde S$ where $J$ 
is injective onto the line segments contained in $J(U)$. Notice that $J$ is locally invertible because $J'=\exp(\tilde K)$ is never vanishing.

Using $J$, we can say a lot more on the geodesics of $\tilde\nabla$.
We begin with 

\newthm Proposition \duno: Let $\tilde\nabla$ be a holomorphic connection on 
$T\tilde S=\tilde S\times\C$, where $\tilde S\cong\C$ or $\Delta$ is a simply
connected Riemann surface. Let $\tilde\eta$ be the holomorphic $1$-form associated
to~$\tilde\nabla$, and $\tilde K\colon\tilde S\to\C$ a holomorphic primitive of
$\tilde\eta$. Finally, let $\tilde J\colon\tilde S\to\C$ be a holomorphic primitive of
$\exp(\tilde K)$. Then a smooth curve $\sigma\colon I\to\tilde S$ is a geodesic
if and only if there are $c_0$, $w_0\in\C$ such that
$$
J\bigl(\sigma(t)\bigr)= c_0t+w_0\;.
$$
In particular, the geodesic issuing from $z_0\in\tilde S$ along the direction~$
v_0\in\C^*$ is given by $\sigma(t)=J^{-1}\bigl(c_0t+J(z_0)\bigr)$, where
$c_0=\exp\bigl(\tilde K(z_0)\bigr)v_0$ and $J^{-1}$ is an analytic continuation
of the local inverse of $J$ nearby $J(z_0)$ chosen so that $J^{-1}\bigl(J(z_0)
\bigr)=z_0$.

\pf The first assertion follows directly from the previous
proposition; but let us describe another proof giving an useful formula.

We know that if $\sigma$ is a geodesic then the support of $\sigma'$
is contained in a leaf of the horizontal foliation. Recalling \eqhfp,
this means that $\sigma'$ must satisfy the differential equation
$$
\sigma'=c_0 \exp\bigl(-\tilde K(\sigma)\bigr)
\neweq\eqddue
$$
for some $c_0\in\C$.

Conversely, if $\sigma$ satisfies this equation it is easy to check that
it satisfies \eqgde\ too (remember that $X_\alpha\equiv 1$ here), and thus
it is a geodesic. But we have
$$
\sigma'=c_0 \exp\bigl(-\tilde K(\sigma)\bigr)\quad\Longleftrightarrow\quad
\exp\bigl(\tilde K(\sigma)\bigr)\sigma'\equiv c_0\quad\Longleftrightarrow\quad
(J'\circ\sigma)\sigma'\equiv c_0\quad\Longleftrightarrow\quad J\bigl(\sigma(t)
\bigr)=c_0t+w_0\;,
$$
and the first assertion follows. The second is an easy consequence of the
fact that $c_0=\exp\bigl(\tilde K(z_0)\bigr)v_0$.\qedn

\newrem In particular, the proof shows that a curve $\sigma\colon[0,\eps)\to
\tilde S$ is a geodesic if and only if
$$
\sigma'(t)=\exp\bigl(-\tilde K\bigl(\sigma(t)\bigr)\bigr)
\exp\bigl(\tilde K\bigl(\sigma(0)\bigr)\bigr)\sigma'(0)
$$ 
if and only if
$$
J\bigl(\sigma(t)\bigr)=\exp\bigl(\tilde K\bigl(\sigma(0)\bigr)\bigr)
\sigma'(0)t+J\bigl(\sigma(0)\bigr)\;.
$$

The first important fact we deduce from this result is that geodesics
cannot accumulate points in $\tilde S$. To better express this fact
let us recall two standard definitions.

\newdef We say that a curve $\gamma\colon[0,\eps)\to\tilde S$ (with $\eps\in(0,+\infty]$) {\sl tends to
the boundary} of~$\tilde S$ if $\gamma(t)$ eventually leaves every compact
subset of $\tilde S$. In other words, $\gamma$ does not tend to the 
boundary if and only if there is a sequence~$t_k\uparrow\eps$ such that
$\gamma(t_k)\to \tilde z_0\in\tilde S$.

\newdef An {\sl asymptotic value} of a holomorphic function $J\colon
\tilde S\to\C$ is a $w_0\in\C$ such that there exists a curve
$\gamma\colon[0,1)\to\tilde S$ tending to the boundary of $\tilde S$
with $J\bigl(\gamma(t)\bigr)\to w_0$ as $t\to 1$.

\newthm Proposition \ddue: Let $\tilde\nabla$ be a holomorphic connection on 
$T\tilde S=\tilde S\times\C$, where $\tilde S\cong\C$ or $\Delta$ is a simply
connected Riemann surface, and let $\sigma_{v_0}\colon[0,\eps_{v_0})\to
\tilde S$ be the maximal geodesic issuing from $z_0\in\tilde
S$ in the direction~$v_0\in\C^*$. Then $\sigma_{v_0}$ tends
to the boundary of~$\tilde S$. Furthermore, if $\eps_{v_0}<+\infty$
then $w_0=J(z_0)+J'(z_0)v_0\eps_{v_0}$ is an 
asymptotic value of $J$.

\pf Let us first consider the case $\eps_{v_0}=+\infty$. If 
$\sigma_{v_0}$ does not tend to the boundary we can find a
sequence $t_k\to+\infty$ so that $\sigma_{v_0}(t_k)\to \tilde z_0\in
\tilde S$. Hence $J\bigl(\sigma_{v_0}(t_k)\bigr)\to J(\tilde z_0)\in\C$;
but $J\bigl(\sigma_{v_0}(t_k)\bigr)=J(z_0)+J'(z_0)
v_0 t_k$ is unbounded, contradiction.

Assume then $\eps_{v_0}<+\infty$, and put $\gamma(t)=\sigma_{v_0}(\eps_{v_0}t)$ 
and $\gamma_1(t)=J\bigl(\gamma(t)\bigr)=
J(z_0)+J'(z_0)v_0\eps_{v_0} t$. 
Clearly $J\bigl(\gamma(t)\bigr)\to w_0$ as $t\to 1$; so to end the proof
it suffices to show that $\gamma$ tends to the boundary of~$\tilde S$.

If $\gamma$ does not tend to the boundary, then
there is a sequence $t_k\to 1$ with $\gamma(t_k)\to\tilde z_0\in\tilde S$;
hence $\gamma_1(t_k)=J\bigl(\gamma(t_k)\bigr)\to J(\tilde z_0)$, and thus $w_0=J(\tilde z_0)\in J(\tilde
S)$.
Let $F\colon D\to\tilde S$ be the local inverse of~$J$ with $F(w_0)=\tilde 
z_0$, where $D$ is a disk centered at~$z_0$. By Proposition~\duno, there is
an inverse $F_1$ of~$J$ defined in a neighbourhood $U$ of the support
of $\gamma_1$ 
with $F_1\circ\gamma_1\equiv\gamma$; up to shrinking~$U$ we
can also assume that $U\cap D$ is connected. We have $\gamma_1(t_k)
\in D\cap U$ eventually. Furthermore, $F(D)$ is an open
neighbourhood of~$\tilde z_0$; hence $\gamma(t_k)\in F(D)$ eventually. But
$$
J\bigl(F_1\bigl(\gamma_1(t_k)\bigr)\bigr)=
\gamma_1(t_k)=J\bigl(F\bigl(\gamma_1(t_k)\bigr)\bigr)\;;
$$
since $J$ is injective in $F(D)$, it follows that $F_1\bigl(\gamma_1(t_k)\bigr)=F\bigl(\gamma_1(t_k)\bigr)$ eventually. But then $F$ and $F_1$ are
two branches of the inverse of $J$ defined in the connected open set
$U\cap D$ and assuming the same value at $\gamma_1(t_k)$; it follows
that $F\equiv F_1$ on $U\cap D$. Therefore the curve 
$t\mapsto F\bigl(J(z_0)+J'(z_0)v_0 t\bigr)$ is a 
geodesic extending $\sigma_{v_0}$ beyond $\eps_{v_0}$,
against the massimality of $\eps_{v_0}$.\qedn

The next step consists in studying the set of points reached by geodesics
issuing from a given base point. 

\newdef For $(z_0,v_0)\in T\tilde S\setminus\tilde S$, let
$\sigma_{v_0}\colon[0,\eps_{v_0})\to\tilde S$ denote the
maximal geodesic issuing from $z_0$ in the direction~$v_0$,
with $\eps_{v_0}\in(0,+\infty]$. 
Put $\ca D_{z_0}=\{v\in\C\mid \eps_{v}>1\}$ and
define $\exp_{z_0}\colon\ca D_{z_0}\to\tilde S$ by
setting $\exp_{z_0}(v)=\sigma_{v}(1)$. 

Then:

\newthm Proposition \dtre: Let $\tilde\nabla$ be a holomorphic connection on 
$T\tilde S=\tilde S\times\C$, where $\tilde S\cong\C$ or $\Delta$ is a simply
connected Riemann surface, and fix $z_0\in\tilde S$. Then:
\smallskip
\itm{(i)} $J\circ\exp_{z_0}=J(z_0)+J'(z_0)\id$;
\itm{(ii)} $\exp_{z_0}$ is a biholomorphism with its image;
\itm{(iii)} $J$ is globally injective on the open simply connected set 
$\exp_{z_0}(\ca
D_{z_0})$ of points that can be joined to~$z_0$ by a
geodesic, and the inverse $J^{-1}\colon J\bigl(\exp_{z_0}(\ca
D_{z_0})\bigr)=J(z_0)+J'(z_0)\ca D_{z_0}\to
\tilde S$ is given by
$$
J^{-1}(w)=\exp_{z_0}\left({w-J(z_0)\over J'(z_0)}
\right)\;.
$$

\pf (i) Take $v\in\ca D_{z_0}$. Then Proposition~\duno\ yields
$$
J\bigl(\exp_{z_0}(v)\bigr)=J\bigl(\sigma_{v}(1)
\bigr)=J(z_0)+J'(z_0)v\;,
$$
as stated.

(ii) Part (i) implies that $\exp_{z_0}$ is injective. The 
holomorphicity follows from the fact that $\sigma_{v}$ solves
the Cauchy problem
$$
\cases{\sigma'=\exp\bigl(-\tilde K(\sigma)\bigr)J'(z_0)v\;,\cr
\sigma(0)=z_0\;,\ \sigma'(0)=v\;,}
$$
and thus $\sigma_{v}(1)$ depends holomorphically on~$v$.
Then part (i) yields
$$
(J'\circ\exp_{z_0})\exp'_{z_0}\equiv J'(z_0)\;;
\neweq\eqdtre
$$
and since $J'=\exp(\tilde K)$ it follows that $\exp_{z_0}'$ is
never vanishing. Being globally injective, $\exp_{z_0}$ is then
a biholomorphism with its image.

(iii) It follows immediately from (i) and (ii), noticing that $\exp_{\tilde
p_0}(\ca D_{z_0})$ is simply connected because $\ca D_{z_0}$
is star-shaped with respect to the origin.\qedn

As a consequence, we have an interesting corollary:

\newthm Corollary \dqua: Let $\tilde\nabla$ be a holomorphic connection on 
$T\tilde S=\tilde S\times\C$, where $\tilde S\cong\C$ or $\Delta$ is a simply
connected Riemann surface. Then a metric $\tilde g$ adapted 
to~$\tilde\nabla$ is never complete, unless $\tilde S\cong\C$ and
$\tilde\nabla$ is the trivial connection.

\pf Assume $\tilde g$ complete. Then the geodesics are defined for all
times, and hence we have $\ca D_z=\C$ for all $z
\in\tilde S$. But then $\exp_z(\ca D_z)$ is a
copy of~$\C$ contained in~$\tilde S$; therefore $\tilde S=\C=
\exp_z(\ca D_z)$. In particular, $\exp_z$
must be affine linear, sending the origin to~$z$ and with 
derivative~1 at the origin, by \eqdtre; thus $\exp_z(v)=
v+z$. From Proposition~\dtre\ it follows that $J$ is
affine linear too; therefore $\exp(\tilde K)=J'$ is constant. Then
$\tilde K$ is constant, and hence $\tilde\eta=\de\tilde K\equiv 0$, that is
$\tilde\nabla$ is the trivial connection.\qedn

We end this section with a couple of remarks on the case $\tilde S=\C$. In this case,
if $\tilde\nabla$ is not trivial, $\exp_z$ cannot be surjective.
In fact, if $\exp_z$ is surjective then $J$ is globally
injective and thus (being $\tilde S=\C$) $J$ must be affine linear and,
as before, we find that $\tilde\nabla$ is trivial. Notice that in
general $J$ has an essential singularity at infinity; therefore
(open simply connected) sets where $J$ is injective tend to become very thin near infinity.

We have also seen that $\ca D_z$ cannot be~$\C$, unless 
$\tilde\nabla$ is trivial; so there are geodesics going to the
boundary in finite time. However, we cannot give a bound on this time.
More precisely, we have

\newthm Proposition \dsei: Let $\tilde\nabla$ be a holomorphic connection 
on $T\C$. Given $z_0\in\C$, let $\eps\colon S^1\to(0,+\infty]$
be defined by $\eps(e^{2\pi i\theta})=\eps_{e^{2\pi i\theta}}=
\sup\{t>0\mid te^{2\pi i\theta}\in\ca D_{z_0}\}$. Then $\eps$ is
unbounded on every interval of~$S^1$.

\pf Assume, by contradiction, that there are $0\le \theta_0<\theta_1<2\pi$
and $M>0$ such that $\eps(e^{2\pi i\theta})<M$ for all $\theta\in[\theta_0,
\theta_1]$. We know that $\exp_{z_0}$ is injective, and that
all geodesics tends to the boundary of~$\tilde S=\C$ (that is, to
infinity); therefore the geodesics issuing from~$z_0$ with
direction $e^{2\pi i\theta}$ with $\theta\in[\theta_1,\theta_2]$ swap
a wedge-like simply connected region $W\subset\C$ bounded by the
geodesics starting with direction $v_0=e^{2\pi i\theta_0}$
e $v_1=e^{2\pi i\theta_1}$. If $v=e^{2\pi i\theta}$ we have
$$
\bigl|J\bigl(\sigma_{v}(t)\bigr)\bigr|\le |J(z_0)|+
|J'(z_0)|\eps(v)<|J(z_0)|+|J'(z_0)|M\;.
$$
So $J$ is bounded on $W$; by a Phragmen-Lindel\"of argument (see
[S, Theorem I\negthinspace I\negthinspace I.3.4]), it follows that
$J\bigl(\sigma_{v_0}(t)\bigr)$ and $J\bigl(\sigma_{v_1}(t)
\bigr)$ must have the same limit. But this would imply $v_0=
v_1$, contradiction.\qedn

When $\tilde S=\Delta$, this argument just says that we cannot have an 
interval of geodesics all converging to the same boundary point of~$\Delta$
in finite time.

\smallsect 3. Multiply connected Riemann surfaces

Now let $S$ be any Riemann surface, and assume we have a holomorphic
connection $\nabla$ on $TS$; in this section we shall study the metric, horizontal
and geodesic foliations induced by $\nabla$ on $TS\setminus S$.

Let us begin with a few preliminaries. Let $\pi\colon\tilde S\to S$ be
the universal covering map. Since $S\ne\P^1(\C)$, the universal
covering space $\tilde S$ is biholomorphic either to~$\C$ or to~$\Delta$.
Let $\tilde\nabla=\pi^*\nabla$ be the holomorphic connection on $T\tilde S$
induced by~$\nabla$ via~$\pi$ (it is well-defined because $\pi$ is
locally invertible); it satisfies the equations
$$
d\pi\bigl(\tilde\nabla_{\tilde v}\tilde e\bigr)=\nabla_{d\pi(\tilde v)}
d\pi(\tilde e)\quad\Longleftrightarrow\quad (\id\otimes d\pi)\circ
\tilde\nabla=(\pi^*\otimes\id)\circ\nabla\circ d\pi\;.
\neweq\eqtuno
$$
Let $(U_\alpha,z_\alpha)$ be a chart of $S$, and $\eta_\alpha$ the local
holomorphic 1-form representing~$\nabla$ on $U_\alpha$. Denote by $w$
the coordinate on $\tilde S=\C$ or $\Delta$, and by $\tilde\eta$ the
global holomorphic 1-form representing $\tilde\nabla$. We  
define a local derivative $\pi'_\alpha\colon\pi^{-1}(U_\alpha)\to\C^*$ by 
$$
d\pi_w\left({\de\over\de w}\right)=\pi'_\alpha(w)
\left.{\de\over\de z_\alpha}\right|_{\pi(w)}\;.
\neweq\eqtumez
$$
Then using \eqtuno\ it is easy to see that $\eta_\alpha$ and $\tilde\eta$
are related by
$$
\tilde\eta=\pi^*\eta_\alpha+{1\over\pi'_\alpha}\,d\pi'_\alpha
\neweq\eqtdue
$$
over $\pi^{-1}(U_\alpha)$. 

As a first consequence we have

\newthm Proposition \tuno: Let $S$ be a (multiply connected) Riemann
surface, and $\nabla$ a holomorphic connection on~$TS$. Let $\tilde\pi
\colon\tilde S\to S$ be the universal covering map, and $\tilde\nabla$
the holomorphic connection on $T\tilde S$ induced by~$\nabla$ via~$\pi$.
Then:
\smallskip
\itm{(i)} $d\pi$ sends leaves of the metric (respectively, horizontal)
foliation in $T\tilde S$ onto leaves of the metric (respectively,
horizontal) foliation in $TS$;
\itm{(ii)} a curve $\tilde\sigma\colon I\to\tilde S$ is a geodesic for
$\tilde\nabla$ if and only if $\sigma=\pi\circ\tilde\sigma$ is a geodesic
for~$\nabla$.

\pf (i) Let $\omega=p^*\eta_\alpha+v_\alpha^{-1}\,dv_\alpha$ be the 
(global) 1-form generating the horizontal foliation on $TS\setminus S$, and
$\tilde\omega=\tilde p^*\tilde\eta+\tilde v^{-1}\,d\tilde v$ the corresponding form
generating the horizontal foliation on $T\tilde S\setminus\tilde S$,
where $p\colon TS\to S$ and $\tilde p\colon T\tilde 
S\to\tilde S$ are the projections.
In local coordinates, we can express $d\pi\colon T\tilde S\to TS$ by
$$
d\pi(w,\tilde v)=\bigl(\pi(w),\pi'_\alpha(w)\tilde v)\;,
$$
that is $v_\alpha\circ d\pi=(\pi'_\alpha\circ\tilde p)\tilde v$. 
Therefore 
$$
\eqalign{
(d\pi)^*\omega&=(d\pi)^*p^*\eta_\alpha+(d\pi)^*\left({1\over v_\alpha}\,d 
v_\alpha\right)=\tilde p^*\pi^*\eta_\alpha+
{1\over v_\alpha\circ d\pi}\,d(v_\alpha\circ d\pi)\cr
&=\tilde p^*\pi^*\eta_\alpha+{1\over\pi'_\alpha\circ\tilde p}\,d(\pi'_\alpha
\circ\tilde p)+{1\over\tilde v}\,d\tilde v=\tilde p^*\tilde\eta+{1\over
\tilde v}\,d\tilde v\cr
&=\tilde\omega\;,
\cr}
\neweq\eqtqua
$$
and this means exactly that $d\pi$ sends leaves of the horizontal 
foliation upstairs onto leaves of the horizontal foliation downstairs.

By Remark~1.5, the metric foliation downstairs (respectively, upstairs)
is generated by the local forms $\varpi_\alpha=|v_\alpha|^2\Re\omega$
(respectively, $\widetilde\varpi=|\tilde v|^2\Re\tilde\omega$). Then
$$
(d\pi)^*\varpi_\alpha=|v_\alpha\circ d\pi|^2\Re(d\pi)^*\omega=
|\pi'_\alpha\circ\tilde p|^2|\tilde v|^2\Re\tilde\omega=
|\pi'_\alpha\circ\tilde p|^2\widetilde\varpi\;,
$$  
and so $d\pi$ also sends leaves of the metric foliation upstairs onto
leaves of the metric foliation downstairs.

(ii) By definition we have
$$
\nabla_{\sigma'}\sigma'=\nabla_{d\pi(\tilde\sigma')}d\pi(\tilde\sigma')
=d\pi(\tilde\nabla_{\tilde \sigma'}\tilde\sigma')\;,
$$
and so $\nabla_{\sigma'}\sigma'\equiv O$ if and only if 
$\tilde\nabla_{\tilde \sigma'}\tilde\sigma'\equiv O$.\qedn

We shall need to know how $\tilde\eta$ behaves under the action of the
automorphism group $\Aut(\pi)$ of $\pi$. 

\newthm Lemma \monrep: Let $S$ be a (multiply connected) Riemann
surface, and $\nabla$ a holomorphic connection on~$TS$. Let $\tilde\pi
\colon\tilde S\to S$ be the universal covering map, and $\tilde\nabla$
the holomorphic connection on $T\tilde S$ induced by~$\nabla$ via~$\pi$.
Let $\tilde\eta$ be the holomorphic form representing $\tilde\nabla$, and 
$\tilde K$ a global primitive of $\tilde\eta$. Then
$$
\exp(\tilde K\circ\gamma)={\rho(\gamma)\over\gamma'}\exp(\tilde K)
\neweq\eqtsei
$$
for all $\gamma\in\Aut(\pi)$, where $\rho\colon\Aut(\pi)\to\C^*$ is the monodromy
representation of $\nabla$ (and we are identifying $\Aut(\pi)$ with the fundamental group 
of~$S$).

\pf Let $\gamma\in\Aut(\pi)$.
From $\pi\circ\gamma=\pi$ we get $\gamma\bigl(\pi^{-1}(U_\alpha)\bigr)
=\pi^{-1}(U_\alpha)$ and
$$
(\pi'_\alpha\circ\gamma)\gamma'=\pi'_\alpha\;.
$$
Therefore
$$
\eqalign{
\gamma^*\tilde\eta&=\gamma^*\pi^*\eta_\alpha+\gamma^*\left(
{1\over\pi'_\alpha}\,d\pi'_\alpha\right)=\pi^*\eta_\alpha+{1\over
\pi'_\alpha\circ\gamma}\,d(\pi'_\alpha\circ\gamma)=\pi^*\eta_\alpha+{\gamma'\over\pi'_\alpha}\,d\left({\pi'_\alpha\over
\gamma'}\right)\cr
&=\tilde\eta-{1\over\gamma'}\,d\gamma'\;.
\cr}
\neweq\eqtcin
$$
Let now $\tilde K$ be a holomorphic primitive of $\tilde\eta$. Then \eqtcin\
becomes
$$
d(\tilde K\circ\gamma-\tilde K)=-{1\over\gamma'}d\gamma'\;,
$$
and thus we can find a $\rho(\gamma)\in\C^*$ such that
$$
\exp(\tilde K\circ\gamma)={\rho(\gamma)\over\gamma'}\exp(\tilde K)\;.
$$
So we are left to proving that $\rho(\gamma)$ is given by the monodromy representation.

Choose an oper cover $\{(U_\alpha,z_\alpha)\}$ of $S$, where the $U_\alpha$ are simply connected (and all non-empty intersections $U_\alpha\cap U_\beta$ are connected).
For any $\alpha$, fix a connected component $\tilde U_{\alpha,\id}$ 
of~$\pi^{-1}(U_\alpha)$; setting $\tilde U_{\alpha,\gamma}=\gamma(\tilde U_{\alpha,\id})$,
varying $\gamma\in\Aut(\pi)$ we get all connected components of~$\pi^{-1}(U_\alpha)$.
In particular, $\{(\tilde U_{\alpha,\gamma},z_\alpha\circ\pi,\de/\de w)\}$ is an open
cover of~$\tilde S$ trivializing $T\tilde S$. By construction, the cocycle representing 
$T\tilde S$ with respect to this cover is trivial. 

Choose a holomorphic primitive $K_\alpha$ of $\eta_\alpha$ on $U_\alpha$. Then 
\eqtdue\ yields constants $c_{\alpha,\gamma}\in\C^*$ such that
$$
\exp(\tilde K)|_{\tilde U_{\alpha,\gamma}}=c_{\alpha,\gamma}\exp(K_\alpha\circ\pi)
\pi'_\alpha|_{\tilde U_{\alpha,\gamma}}\;.
$$
Assume that $U_\alpha\cap U_\beta\ne\void$. Then for every $\gamma\in\Aut(\pi)$
there is a unique $\gamma'\in\Aut(\pi)$ so that $\tilde U_{\alpha,\gamma}\cap
\tilde U_{\beta,\gamma'}\ne\void$. In this intersection we have
$$
1={\exp(\tilde K)\over\exp(\tilde K)}={c_{\alpha,\gamma}\exp(K_\alpha\circ\pi)
\pi'_\alpha\over c_{\beta,\gamma'}\exp(K_\beta\circ\pi)
\pi'_\beta}={c_{\alpha,\gamma}\over c_{\beta,\gamma'}}{\exp(K_\alpha\circ\pi)
\over\exp(K_\beta\circ\pi)}\left({\de z_\alpha\over\de z_\beta}\circ\pi
\right)=
{c_{\alpha,\gamma}\over c_{\beta,\gamma'}}\hat\psi_{\alpha\beta}\;,
\neweq\eqrep
$$
where $\{\hat\psi_{\alpha\beta}\}$ is the (locally constant) cocycle representing the monodromy representation of~$\nabla$; see \eqhatxi.

Now take $\gamma_0$, $\gamma\in\Aut(\pi)$. Then
$$
\eqalign{
\exp(\tilde K\circ\gamma)|_{\tilde U_{\alpha,\gamma_0}}&=
\exp(\tilde K)|_{\tilde U_{\alpha,\gamma\gamma_0}}\circ\gamma=
c_{\alpha,\gamma\gamma_0}\exp(K_\alpha\circ\pi\circ\gamma)(\pi'_\alpha\circ\gamma)
={c_{\alpha,\gamma\gamma_0}\over\gamma'}\exp(K_\alpha\circ\pi)\pi'_\alpha\cr
&={c_{\alpha,\gamma\gamma_0}\over c_{\alpha,\gamma}}{1\over\gamma'}
\exp(\tilde K)|_{\tilde U_{\alpha,\gamma_0}}\;.
\cr}
$$
Therefore $\rho(\gamma)=c_{\alpha,\gamma\gamma_0}/c_{\alpha,\gamma}$, and the 
assertion follows from \eqrep\ and the definition of the canonical isomorphism
between \v Cech cohomology and singular cohomology.\qedn
 
\newdef We shall denote by
$\rho(\pi)\subseteq\C^*$ the image of $\Aut(\pi)$ under~$\rho$, and by
$|\rho|(\pi)\subseteq\R^+$ the image of $\Aut(\pi)$ under~$|\rho|$; in
particular $\nabla$ has real periods if and only if $|\rho|(\pi)=\{1\}$.

Using the monodromy representation we can describe the
metric and horizontal foliations:

\newthm Theorem \tdue: Let $S$ be a (multiply connected) Riemann
surface, and $\nabla$ a holomorphic connection on~$TS$. Let $L$ be a
leaf of the metric foliation, and take $v_0\in T_{z_0}S\cap L$.
Then 
$$
L\cap T_{z_0}S=|\rho|(\pi)\cdot (S^1\cdot v_0)\qquad\hbox{and}\qquad
\bar{L}\cap T_{z_0}S=\bar{|\rho|(\pi)}\cdot (S^1\cdot v_0)\;.
\neweq\eqtgrp
$$
In particular, either
\smallskip
\itm{(i)} $\nabla$ has real periods, and in that case all leaves of the
metric foliation are closed in $TS$; or,
\itm{(ii)} $\nabla$ has not real periods, and in that case all leaves
of the metric foliation accumulate all 
points of the\break\indent zero section of~$TS$. 

\pf Clearly $S^1\cdot v_0\subset L\cap T_{z_0}S$. Let $\pi\colon\tilde S\to
S$ be the universal covering map; by Proposition~\tuno\ we can find
a leaf $\tilde L$ of the metric foliation upstairs so that $L=d\pi(\tilde L)$.
Fix $\tilde z_0\in\pi^{-1}(z_0)$, and let $\tilde v_0\in T_{\tilde z_0}
\tilde S$ such that $d\pi_{\tilde z_0}(\tilde v_0)=v_0$. Then $(\tilde z,
\tilde v)\in\tilde L$ if and only if 
$$
\exp\bigl(\Re\tilde K(\tilde z)\bigr)|\tilde v|=
\exp\bigl(\Re\tilde K(\tilde z_0)\bigr)|\tilde v_0|\;,
$$
where $\tilde K$ is a holomorphic primitive of the holomorphic 1-form $\tilde
\eta$ representing the holomorphic connection~$\tilde\nabla$ induced 
by~$\nabla$ via~$\pi$. In particular, if $\tilde z=\gamma(\tilde z_0)$ for
some $\gamma\in\Aut(\pi)$, then \eqtsei\ implies that 
$(\tilde z,\tilde v)\in\tilde L$ if and only if
$$
|\tilde v|={|\gamma'(\tilde z_0)|\over|\rho(\gamma)|}|\tilde v_0|\;.
\neweq\eqtotto
$$
In other words,
$$
\tilde L_{\gamma(\tilde z_0)}={|\gamma'(\tilde z_0)|\over|\rho(\gamma)|}
\cdot\tilde L_{\tilde z_0}\;,
$$
where we put $\tilde L_{\tilde z}=\tilde L\cap T_{\tilde z}\tilde S$. From
$d\pi_{\tilde z_0}=\gamma'(\tilde z_0)d\pi_{\gamma(\tilde z_0)}$ we then get
$$
d\pi_{\gamma(\tilde z_0)}(\tilde L_{\gamma(\tilde z_0)})=
{1\over|\rho(\gamma)|}\,d\pi_{\tilde z_0}(\tilde L_{\tilde z_0})
=|\rho(\gamma^{-1})|\,d\pi_{\tilde z_0}(\tilde L_{\tilde z_0})\;.
$$
But $d\pi_{\tilde z_0}(\tilde L_{\tilde z_0})=S^1\cdot v_0$;
hence
$$
L\cap T_{z_0}S=\bigcup_{\gamma\in\Aut(\pi)}d\pi_{\gamma(\tilde z_0)}
(\tilde L_{\gamma(\tilde z_0)})=|\rho|(\pi)\cdot(S^1\cdot v_0)\;.
$$
Taking the closure we get \eqtgrp. In particular, if $\nabla$ has not
real periods then $0\in\bar{|\rho|(\pi)}$, and (ii) follows.

Finally, assume that $\nabla$ has real periods.
Let $L\subset TS\setminus S$ be a leaf of the metric foliation, 
and $\{(z_k,v_k)\}\subset L$ with $z_k\to z_0\in S$ 
and $v_k\to v_0\in T_{z_0}S$; to get (i) we must prove that $(z_0,v_0)\in L$. 

Let $\tilde L$ be the leaf of the metric foliation upstairs such that
$L=d\pi(\tilde L)$, and take a sequence $\{(\tilde z_k,\tilde v_k)\}\subset\tilde L$ such that $\pi(\tilde z_k)=z_k$
and $d\pi_{\tilde z_k}(\tilde v_k)=v_k$.
Fix a point $\hat z_0\in\pi^{-1}(z_0)$. Since $z_k\to z_0$ in~$S$, 
we can find a sequence $\hat z_k\to \hat z_0$ in~$\tilde S$ and a
sequence $\{\gamma_k\}\subset\Aut(\pi)$ such that 
$\tilde z_k=\gamma_k(\hat z_k)$. Put $\hat v_k=\tilde v_k/\gamma'_k(\hat
z_k)$, so that $d\pi_{\hat z_k}(\hat v_k)=d\pi_{\tilde z_k}(\tilde v_k)=v_k$.
Furthermore, \eqtotto\ yields $(\hat z_k,\hat v_k)\in\tilde L$.

Now, since $\exp\bigl(\Re\tilde K(\hat z_k)\bigr)|\hat v_k|$ is a
non zero constant and $\exp\bigl(\Re\tilde K(\hat z_k)\bigr)\to
\exp\bigl(\Re\tilde K(\hat z_0)\bigr)\ne 0$, the sequence $\{\hat v_k\}$
is bounded; therefore, up to a subsequence, we can assume that
$\hat v_k\to\hat v_0\in\C$. Clearly, $(\hat z_0,\hat v_0)$ still belongs to
the leaf~$\tilde L$; hence 
$$
d\pi_{\hat z_0}(\hat v_0)=\lim_{k\to+\infty} d\pi_{\hat w_k}(\hat v_k)=
\lim_{k\to\infty} v_k=v_0
$$
belongs to~$L$, as claimed.\qedn

In a similar way we can prove the following

\newthm Theorem \tdueb: Let $S$ be a (multiply connected) Riemann
surface, and $\nabla$ a holomorphic connection on~$TS$, and let $L$ be a
leaf of the horizontal foliation. Then $p(L)=S$, where $p\colon TS\to S$ is the canonical projection. Furthermore, take 
any $z_0\in S$ and $v_0\in T_{z_0}S\cap L$.
Then 
$$
L\cap T_{z_0}S=\rho(\pi)\cdot v_0\qquad\hbox{and}\qquad
\bar{L}\cap T_{z_0}S=\bar{\rho(\pi)}\cdot v_0\;.
\neweq\eqtgrpb
$$
In particular, either
\smallskip
\itm{(i)} $\nabla$ has real periods, and in that case either all leaves of the
horizontal foliation are closed in $TS$ or any leaf of the horizontal
foliation is dense in the leaf of
the metric foliation containing it; or,
\itm{(ii)} $\nabla$ has not real periods, and in that case all leaves
of the horizontal foliation accumulate all 
points of\break\indent the zero section of~$TS$. 

\pf Let $\pi\colon\tilde S\to
S$ be the universal covering map; by Proposition~\tuno\ we can find
a leaf $\tilde L$ of the horizontal foliation upstairs so that $L=d\pi(\tilde L)$. Since \eqhfp\ implies $\tilde p(\tilde L)=\tilde S$, it follows 
immediately that $p(L)=S$.

Fix $\tilde z_0\in\pi^{-1}(z_0)$, and let $\tilde v_0\in T_{\tilde z_0}
\tilde S$ such that $d\pi_{\tilde z_0}(\tilde v_0)=v_0$. Then $(\tilde z,
\tilde v)\in\tilde L$ if and only if 
$$
\exp\bigl(\tilde K(\tilde z)\bigr)\tilde v=
\exp\bigl(\tilde K(\tilde z_0)\bigr)\tilde v_0\;,
$$
where $\tilde K$ is a holomorphic primitive of the holomorphic 1-form $\tilde
\eta$ representing the holomorphic connection~$\tilde\nabla$ induced 
by~$\nabla$ via~$\pi$. In particular, if $\tilde z=\gamma(\tilde z_0)$ for
some $\gamma\in\Aut(\pi)$, then \eqtsei\ implies that 
$(\tilde z,\tilde v)\in\tilde L$ if and only if
$$
\tilde v={\gamma'(\tilde z_0)\over\rho(\gamma)}\tilde v_0\;.
\neweq\eqtotto
$$
Notice that $\tilde L$ intersect each $T_{\tilde z}\tilde S$ in just one
point. From
$d\pi_{\tilde z_0}=\gamma'(\tilde z_0)d\pi_{\gamma(\tilde z_0)}$ we then get
$$
d\pi_{\gamma(\tilde z_0)}(\tilde v)=
{1\over\rho(\gamma)}\,d\pi_{\tilde z_0}(\tilde v_0)
=\rho(\gamma^{-1})\,d\pi_{\tilde z_0}(\tilde v_0)\;.
$$
Hence
$$
L\cap T_{z_0}S=\bigcup_{\gamma\in\Aut(\pi)}d\pi_{\gamma(\tilde z_0)}
(\tilde L\cap T_{\gamma(\tilde z_0)}\tilde S)=\rho(\pi)\cdot v_0\;.
$$
Taking the closure we get \eqtgrp. In particular, if $\nabla$ has not
real periods then $0\in\bar{\rho(\pi)}$, and (ii) follows.

Finally, if $\nabla$ has real periods then $\rho(\pi)$ is a subgroup of~$S^1$.
The subgroups of $S^1$ are either cyclic or dense; in the first case it
is easy to check that $L$ is closed, and in the second case the
assertion follows from \eqtgrpb\ and \eqtgrp.\qedn

The monodromy representation enters in another question: deciding when
the automorphisms of~$\pi$ are isometries for a Hermitian metric
adapted to~$\tilde\nabla$.

\newthm Proposition \tduemez: Let $S$ be a (multiply connected) Riemann
surface, and $\nabla$ a holomorphic connection on~$TS$. Let $\tilde\pi
\colon\tilde S\to S$ be the universal covering map, $\tilde\nabla$
the holomorphic connection on $T\tilde S$ induced by~$\nabla$ via~$\pi$,
and $\tilde g=\exp(\Re\tilde K)g_0$ a Hermitian metric adapted 
to~$\tilde\nabla$ (where $g_0$ is the Euclidean metric). 
Let $\rho\colon\Aut(\pi)\to\C^*$ be the monodromy representation of~$\nabla$. Then:
\smallskip
\itm{(i)} $\gamma\in\Aut(\pi)$ is an isometry of $\tilde g$ if and only if $|\rho(\gamma)|=1$; in particular, $\Aut(\pi)\subseteq\mathop{\rm Iso}(\tilde g)$ if and only if
$\nabla$ has real periods;
\itm{(ii)} every $\gamma\in\Aut(\pi)$ sends $\tilde\nabla$-geodesics in
$\tilde\nabla$-geodesics, and we have
$$
\gamma\circ\exp_{\tilde z_0}(\tilde v)=\exp_{\gamma(\tilde z_0)}
\bigl(\gamma'(\tilde z_0)\tilde v\bigr)
$$
for all $\tilde z_0\in\tilde S$ and $\tilde v\in\C$;
\itm{(iii)} if $J\colon\tilde S\to\C$ is a primitive of $\exp(\tilde K)$ then
$$
J\bigl(\gamma(\tilde z)\bigr)-J\bigl(\gamma(\tilde z_0)\bigr)
=\rho(\gamma)\bigl[J(\tilde z)-J(\tilde z_0)\bigr]
$$
\indent for all $\tilde z_0$, $\tilde z\in\tilde S$ and $\gamma\in\Aut(\pi)$.

\pf (i) Formula \eqtsei\ yields
$$
\exp\bigl(\Re\tilde K\circ\gamma)={|\rho(\gamma)|\over|\gamma'|}
\exp(\Re\tilde K)\;.
$$
Therefore 
$$
\eqalign{
\tilde g\bigl(\gamma(\tilde z);d\gamma_{\tilde z}(\tilde v)\bigr)&=\exp
\bigl(2\Re\tilde K\bigl(\gamma(\tilde z)\bigr)\bigr)|\gamma'(\tilde z)|^2
|\tilde v|^2\cr
&=|\rho(\gamma)|^2\exp\bigl(2\Re\tilde K(\tilde z)\bigr)|\tilde v|^2=
|\rho(\gamma)|^2\tilde g(\tilde z;\tilde v)\;,
\cr}
$$
and $\gamma\in\Aut(\pi)$ is an isometry if and only if $|\rho(\gamma)|=1$.

(ii) Remark~2.1 says that a curve
$\sigma\colon I\to\tilde S$ is a geodesic if and only if
$$
\sigma'=\exp(-\tilde K(\sigma)\bigr)\exp\bigl(\tilde K\bigl(\sigma(t_0)
\bigr)\bigr)\sigma'(t_0)
$$
for some (and hence all) $t_0\in I$. Now
$$
\eqalign{
\exp\bigl(-\tilde K(\gamma\circ\sigma)\bigr)&\exp\bigl(\tilde K\bigl(
\gamma\circ\sigma(t_0)\bigr)\bigr)(\gamma\circ\sigma)'(t_0)\cr
&=
{\gamma'\circ\sigma\over\rho(\gamma)}\exp(-\tilde K\circ\sigma)
{\rho(\gamma)\over\gamma'\bigl(\sigma(t_0)\bigr)}
\exp\bigl(\tilde K\bigl(\sigma(t_0)\bigr)\bigr)\gamma'\bigl(\sigma(t_0)
\bigr)\sigma'(t_0)\cr
&=(\gamma'\circ\sigma)\bigl[\exp(-\tilde K\circ\sigma)
\exp\bigl(\tilde K\bigl(\sigma(t_0)\bigr)\bigr)\sigma'(t_0)\bigr]\;,\cr}
$$
and thus $\gamma\circ\sigma$ is a geodesic if and only if $\sigma$ is,
even when $\gamma$ is not an isometry.

(iii) Formula \eqtsei\ yields
$$
\exp(\tilde K\circ\gamma)\gamma'=\rho(\gamma)\exp(\tilde K)\;.
$$
Integrating this from $\tilde z_0$ to $\tilde z$ we get the assertion.\qedn

%

In the next
result we show how to compute the monodromy representation when $S\subseteq\C$, that is
when $S$ is covered by a single chart. The most interesting case will be
when $S$ is the complement in $\P^1(\C)$ of a finite set of points.

\newthm Proposition \ttre: Let $S\subseteq\C$ be a (multiply
connected) domain, $\nabla$ a holomorphic connection on~$TS$, and $\eta$
the holomorphic $1$-form representing~$\nabla$. 
Then the monodromy representation $\rho\colon H_1(S,\Z)\to\C^*$ is given by
$$
\rho(\gamma)=\exp\left(\int_\gamma \eta\right)
$$
for all $\gamma\in H_1(S,\Z)$.

\pf Let $\pi\colon\tilde S\to S$ be the universal covering map of~$S$. 
Choose a $z_0\in S$, a $\tilde z_0\in\pi^{-1}(z_0)$, and a loop (still
denoted by~$\gamma$) based at~$z_0$ representing $\gamma\in H_1(S,\Z)$.
Let $\tilde\gamma$ be the lift of~$\gamma$ based at~$\tilde z_0$; then the
action on~$\tilde z_0$ of the element of $\Aut(\pi)$ corresponding to~$\gamma$
(again denoted by $\gamma$)
is given by $\tilde\gamma(1)$.

Let $\tilde\nabla$ be the holomorphic connection on $T\tilde S$ induced
by $\nabla$ via~$\pi$, and $\tilde\eta$ the holomorphic 1-form
representing~$\tilde\nabla$. Choose a holomorphic primitive $\tilde K$
of $\tilde\eta$, and a determination of $\log\pi'$. Then
$$
\eqalign{
\tilde K\bigl(\gamma(w_0)\bigr)-\tilde K(w_0)&=
\tilde K\bigl(\tilde\gamma(1)\bigr)-\tilde K\bigl(\tilde\gamma(0)\bigr)\cr
&=\int_{\tilde\gamma}\tilde\eta=\int_{\tilde\gamma}\left(\pi^*\eta+d\log\pi'
\right)=\int_{\tilde\gamma}\pi^*\eta+(\log\pi')\bigl(\gamma(\tilde z_0)\bigr)-
\log\pi'(\tilde z_0)\cr
&=\int_\gamma\eta+(\log\pi')\bigl(\gamma(\tilde z_0)\bigr)-
\log\pi'(\tilde z_0)\;.
\cr}
\neweq\eqtnove
$$
Therefore
$$
\exp\bigl[\tilde K\bigl(\gamma(\tilde z_0)\bigr)-\tilde K(\tilde z_0)\bigr]=
{\pi'\bigl(\gamma(\tilde z_0)\bigr)\over\pi'(\tilde z_0)}
\exp\left(\int_\gamma\eta\right)
={1\over\gamma'(\tilde z_0)}\exp\left(\int_\gamma\eta\right)\;,
$$
and the assertion follows from \eqtsei.\qedn

We have actually proved something more. Keeping the notations 
introduced in the previous proof, from $(\pi'\circ\gamma)\gamma'=\pi'$
we deduce that for each $\gamma\in\Aut(\pi)$ there is a unique
determination of the logarithm of $\gamma'$ such that
$$
\log\gamma'=\log\pi'-(\log\pi')\circ\gamma\;.
$$
Then \eqtnove\ becomes
$$
\tilde K\bigl(\gamma(w_0)\bigr)=\tilde K(w_0)-\log\gamma'(w_0)+\int_\gamma
\eta\;.
$$
Put 
$$
\rho_0(\gamma)={1\over2\pi i}\int_\gamma \eta\;;
$$
then it is easy to check that $\rho_0\colon H_1(S,\Z)\to\C$ is a
homomorphism of abelian groups, and we can write \eqtnove\ as
$$
\tilde K\circ\gamma=\tilde K-\log\gamma'+2\pi i
\rho_0(\gamma)\;.
\neweq\eqtuzero
$$

\newdef The homomorphism $\rho_0\colon H_1(S,\Z)\to\C$ just introduced is the {\sl period
map} associated to~$\nabla$.

Since $\rho=\exp(2\pi i\rho_0)$, the connection $\nabla$ has 
real periods if and only if the image of the period map is contained 
in~$\R$. In particular, Proposition~\udue\ yields:

\newthm Corollary \tsette: 
Let $S\subseteq\C$ be a domain in the plane, and $\nabla$ a holomorphic connection on~$TS$. Then there exists a
Hermitian metric adapted to~$\nabla$ if and only if the 
period map is real-valued.

\smallsect 4. Meromorphic connections

Let us now specialize to the case we are mostly interested in, that
is meromorphic connections on~$\P^1(\C)$. If $\nabla$ is a meromorphic
connection on $\P^1(\C)$, then we can consider it as a holomorphic
connection on $S=\P^1(\C)\setminus\{p_0,\ldots,p_r\}$, where $\{p_0,\ldots,
p_r\}$ are the poles of the meromorphic connection. Without loss
of generality, we shall always assume $p_0=\infty$, so that $S\subseteq\C$,
and thus we have the period map $\rho_0\colon H_1(S,\Z)\to\C$
associated to~$\nabla$. The homology group $H_1(S,\Z)$ is generated by
the counterclock-wise loops~$\gamma_1,\ldots,\gamma_r$ around,
respectively, $p_1,\ldots,p_r$; therefore $\rho_0(\gamma_j)$ is,
practically by definition, the {\sl residue} $\Res_{p_j}(\nabla)$:
$$
\Res_{p_j}(\nabla)=\rho_0(\gamma_j)={1\over2\pi i}\int_{\gamma_j}\eta\;,
$$
where $\eta$ is the holomorphic 1-form representing~$\nabla$ on~$S$. 
We also have the residue at~$\infty$, which is given by $\rho(\gamma_0)$,
where $\gamma_0$ is a clock-wise Jordan loop in $\C$ containing $\{p_1,
\ldots,p_r\}$ in its interior.
It is
useful to keep in mind that the classical residue theorem for 
meromorphic connections (see, e.g., [IY, Theorem~I\negthinspace 
I\negthinspace I.17.33]) says that
$$
\sum_{j=0}^r\Res_{p_j}(\nabla)=\deg T\P^1(\C)=-2\;.
\neweq\eqqtr
$$

The aim of this section is to describe the recurrence properties of the
geodesics on $\nabla$, and of the geodesic flow on $TS$, recurrence
properties that, as we shall see, are strikingly different from 
the recurrence properties of the metric and horizontal foliations
described by Theorems~\tdue\ and~\tdueb. 

To state the main technical tool for our study we need a definition.

\newdef Let $\nabla$ be a meromorphic connection on 
$\P^1(\C)$, with poles $\{p_0=\infty,p_1,\ldots,p_r\}$, and set
$S=\P^1(\C)\setminus\{p_0,\ldots,p_r\}\subseteq\C$. An {\sl $s$-sided
geodesic polygon} is a simply connected domain $R_0\subset
\P^1(\C)$ whose boundary 
is composed by $s\ge 1$ simple geodesics $\sigma_j\colon[0,\ell_j]\to S$ with
$\sigma_j(\ell_j)=\sigma_{j+1}(0)=z_{j+1}$ for $j=1,\ldots,s-1$, 
$\sigma_r(\ell_n)=\sigma_1(0)=z_1$ and no other intersections; the geodesics are listed so that $\de R_0$ is
positively oriented (that is $R_0$ is the interior of $\de R_0$). The points
$z_1,\ldots,z_s$ are the {\sl vertices} of $R_0$.

Then

\newthm Theorem \qzero: Let $\nabla$ be a meromorphic connection on 
$\P^1(\C)$, with poles $\{p_0=\infty,p_1,\ldots,p_r\}$, and set
$S=\P^1(\C)\setminus\{p_0,\ldots,p_r\}\subseteq\C$. Let $R_0\subset
\P^1(\C)$ be an $s$-sided geodesic polygon
with vertices $z_1,\ldots,z_s$. For $j=1,\ldots,s$ let $\eps_j\in(-\pi,\pi)$ be the external
angle in~$z_j$, and let $\{p_1,\ldots,p_g\}$ be the poles of~$\nabla$ contained in~$R_0$.
Then
$$
\sum_{j=1}^s\eps_j=2\pi\left(1+\sum_{j=1}^g \Re\Res_{p_j}(\nabla)\right)\;.
\neweq\eqqerrza
$$
In particular,
$$
\sum_{j=1}^g \Re\Res_{p_j}(\nabla)\in \left(-{s+2\over2},{s-2\over 2}\right)\;.
\neweq\eqqerra
$$

\pf The idea is to apply the Gauss-Bonnet theorem, even if we do not have
a global metric. What we do have is the local Gauss-Bonnet theorem, 
expressed in terms of a local metric adapted to~$\nabla$. We know
that the Gaussian curvature of all such metrics is identically zero (see
Remark~1.2); furthermore, any such metric is a positive multiple of any
other one, and thus the external angles (and the notion of orthogonal
parametrizations) are the same for all of them. It follows that also the
integral of the geodesic curvature does not depend on the chosen local
metric (and this can be verified directly too; see below), and hence the 
standard proof (see, e.g., [AT2, Theorem~6.3.9] or [dC, p.~274]) of
the global Gauss-Bonnet theorem based on the local Gauss-Bonnet theorem
still works.

We shall apply the Gauss-Bonnet theorem to the region $R$ obtained removing
from~$R_0$ small disks around~$p_1,\ldots,p_g$. Denoting by 
by $\tau^1,\ldots,\tau^g\colon
[0,2\pi]\to S$ the small clock-wise circles bounding the disks around $p_1,\ldots,p_g$ respectively, and by 
$\kappa_g^j$ the geodesic curvature of~$\tau^j$,
the Gauss-Bonnet theorem says that
$$
\sum_{j=1}^g\int_{\tau^j}\kappa_g^j\,ds+\sum_{j=1}^s\eps_j=2\pi(1-g)\;.
\neweq\eqqzero
$$ 
Since all the local metrics adapted to $\nabla$ are (non-constant) 
multiples of the euclidean metric, the standard real coordinates $(x,y)$
on~$S\subseteq\C$ are orthogonal, and we can use the formula
$$
\int_{\tau^j}\kappa_g^j=\int_0^{2\pi}\left[{d\theta^j\over dt}+
{1\over2\sqrt{(EG)\circ\tau^j}}\left((\Im\tau^j)'{\de G\over\de x}\circ
\tau^j-(\Re\tau^j)'{\de E\over\de y}\circ\tau^j
\right)\right]\,dt\;,
\neweq\eqquno
$$
where $E=G=\exp(2\Re K)$ and $\theta^j$ is the angle between $\de/\de x$
and $(\tau^j)'$. Here $K$ is any local holomorphic primitive of the
form $\eta$ representing $\nabla$ in~$S$; since two such local primitives
differ only by an additive constant, the integrand of \eqquno\ does not
depend on the choice of $K$.

We have $\tau^j(t)=p_j+re^{-it}$ for $r>0$ small enough; hence $d\theta^j/dt
\equiv-1$. The Cauchy-Riemann equations say that
$$
{\de\Re K\over\de x}=\Re{\de K\over\de z}\quad\hbox{and}\quad
{\de\Re K\over\de y}=-\Im{\de K\over\de z}\;;
$$ 
therefore
$$
{1\over 2\sqrt{EG}}{\de G\over\de x}=\Re{\de K\over\de z}\quad\hbox{and}
\quad -{1\over 2\sqrt{EG}}{\de E\over\de y}=\Im{\de K\over\de z}\;.
$$
Hence \eqquno\ becomes
$$
\int_{\tau^j}\kappa_g^j=-2\pi+\int_0^{2\pi}\Im{d\over dt}(K\circ\tau^j)
\,dt\;.
$$

Let $\pi\colon\tilde S\to S$ be the universal covering map, and 
$\tilde\tau^j$ a lifting of~$\tau^j$. If $\gamma_j\in\Aut(\pi)$ is
the generator associated to the counterclock-wise loop around~$p_j$, we have
$\tilde\tau^j(2\pi)=\gamma_j^{-1}\bigl(\tilde\tau^j(0)\bigr)$.
Let now $\tilde K$ be a holomorphic primitive of the holomorphic
1-form $\tilde\eta$ representing the holomorphic connection $\tilde\nabla$
induced by $\nabla$ via~$\pi$. Then, choosing a determination for the
logarithm of~$\pi'$, we have
$$
d\tilde K=d(K\circ\pi+\log\pi')\;.
$$
Therefore
$$
{d\over dt}(K\circ\tau^j)={d\over dt}(K\circ\pi\circ\tilde\tau^j)
={d\over dt}\bigl((\tilde K-\log\pi')\circ\tilde\tau^j\bigr)\;.
$$
Hence, using \eqtuzero\ and $(\pi'\circ\gamma)\gamma'=\pi'$ for all
$\gamma\in\Aut(\pi)$, we get 
$$
\eqalign{
\int_{\tau^j}\kappa_g^j&=-2\pi+\Im\Bigl[\tilde K\bigl(\gamma_j^{-1}\bigl(\tilde\tau^j(0)\bigr)\bigr)-\tilde K\bigl(\tilde\tau^j(0)\bigr)+
\log\pi'\bigl(\tilde\tau^j(0)\bigr)-\log\pi'\bigl(\gamma_j^{-1}\bigl(\tilde\tau^j(0)\bigr)\bigr)\Bigr]\cr
&=-2\pi-\Im\bigl(2\pi i\rho_0(\gamma_j)\bigr)=-2\pi\bigl[1+\Re\Res_{p_j}
(\nabla)\bigr]\;,
\cr}
$$
and so 
$$
\sum_{j=1}^g\int_{\tau^j}\kappa^j_g=-2\pi g-2\pi\sum_{j=1}^g\Re\Res_{p_j}(\nabla)\;.
\neweq\eqqtre
$$
Putting this into \eqqzero\ we get \eqqerrza. Finally, \eqqerra\ follows from
\eqqerrza\ and the fact that the sum of the external angles belongs to the
interval $(-s\pi,s\pi)$.\qedn

\newthm Corollary \quno: Let $\nabla$ be a meromorphic connection on 
$\P^1(\C)$, with poles $\{p_0=\infty,p_1,\ldots,p_r\}$, and set
$S=\P^1(\C)\setminus\{p_0,\ldots,p_r\}\subseteq\C$. Let $\sigma\colon
[0,\ell]\to S$ be a geodesic with $\sigma(0)=\sigma(\ell)$ and no
other self-intersections; in particular, $\sigma$ is an oriented
Jordan curve. Let $\{p_1,\ldots,p_g\}$ be the poles of $\nabla$
contained in the interior of~$\sigma$, and $\eps\in(-\pi,\pi)$ the external
angle at~$\sigma(0)$. Then
$$
\eps=2\pi\left(1+\sum_{j=1}^g \Re\Res_{p_j}(\nabla)\right)\;,
\neweq\eqqerrz
$$
and hence
$$
\sum_{j=1}^g \Re\Res_{p_j}(\nabla)\in (-3/2,-1/2)\;.
\neweq\eqqerr
$$

\pf It follows from Theorem~\qzero\ with $s=1$.\qedn

\newthm Corollary \qdue: Let $\nabla$ be a meromorphic connection on 
$\P^1(\C)$, with poles $\{p_0=\infty,p_1,\ldots,p_r\}$, and set
$S=\P^1(\C)\setminus\{p_0,\ldots,p_r\}\subseteq\C$. Let $\sigma_0\colon[0,\ell_0]\to S$ 
and $\sigma_1\colon[0,\ell_1]\to S$ be two distinct geodesics with
$\sigma_0(0)=z_0=\sigma_1(0)$ and $\sigma_0(\ell_0)=z_1=\sigma_1(\ell_1)$
and not intersecting elsewhere. Let $\{p_1,\ldots,p_g\}$ be the poles of $\nabla$
contained in the simply connected domain~$R_0$ bounded by $\sigma_0$ and~$\sigma_1$, and 
$\eps_j\in(-\pi,\pi)$ the external
angle at~$z_j$, for $j=1$,~$2$. Then
$$
\eps_0+\eps_1=2\pi\left(1+\sum_{j=1}^g \Re\Res_{p_j}(\nabla)\right)\;,
\neweq\eqqerrzb
$$
and hence
$$
\sum_{j=1}^g \Re\Res_{p_j}(\nabla)\in (-2,0)\;.
\neweq\eqqerrb
$$

\pf It follows from Theorem~\qzero\ with $s=2$.\qedn

To prove the next corollary we need a lemma and a definition.

\newthm Lemma \ops: Let $\nabla$ be a meromorphic connection on 
$\P^1(\C)$, with poles $\{p_0=\infty,p_1,\ldots,p_r\}$, and set
$S=\P^1(\C)\setminus\{p_0,\ldots,p_r\}\subseteq\C$. Let $R_0\subset
\P^1(\C)$ be an $s$-sided geodesic polygon
with vertices $z_1,\ldots,z_s$. 
Let $\sigma_j\colon[0,\ell_j]\to S$ be the geodesics composing the
boundary of~$R_0$, with $\sigma_j(0)=z_j$ for $j=1,\ldots,s$,
and let $\{p_1,\ldots,p_g\}$ be the poles of~$\nabla$ contained in~$R_0$. Then
$$
\prod_{j=1}^s\sigma'_j(\ell_j)=
\exp\left(-2\pi i\sum_{j=1}^g\Res_{p_j}(\nabla)\right)
\prod_{j=1}^s\sigma'_j(0)\;.
\neweq\eqqqua
$$

\pf Choose a point 
$\tilde z_1\in\pi^{-1}(z_1)$, and let $\tilde\sigma_1\colon[0,\ell_1]\to\tilde S$ be the
lifting of~$\sigma_1$ with $\tilde\sigma_1(0)=\tilde z_1$. Then recursively choose 
the lifting~$\tilde\sigma_j\colon[0,\ell_j]\to\tilde S$ of~$\sigma_j$ with $\tilde\sigma_j(0)
=\tilde\sigma_{j-1}(\ell_{j-1})$. In particular, $\tilde\sigma_s(\ell_s)=\gamma(w_1)$, where 
$\gamma\in\Aut(\pi)$ is the element associated to the class of~$\de R_0$ in~$\pi_1(S,z_1)$.

By Proposition~\tuno.(ii), the $\tilde\sigma_j$ 
are geodesics for $\tilde\nabla$; hence (Remark~2.1)
$$
\tilde\sigma'_j(\ell_j)=\exp\bigl(-\tilde K\bigl(\tilde\sigma_j(\ell_j)\bigr)\bigr)
\exp\bigl(-\tilde K\bigl(\tilde\sigma_j(0)\bigr)\bigr)\tilde\sigma'_j(0)\;.
$$
Recalling that $\tilde\sigma_j(0)=\tilde\sigma_{j-1}(\ell_{j-1})$ 
and using \eqtuzero\ we get
$$
\eqalign{
\prod_{j=1}^s\tilde\sigma'_j(\ell_j)&=\exp\bigl(-\tilde K\bigl(\gamma(\tilde z_1)\bigr)\bigr)
\exp\bigl(\tilde K(\tilde z_1)\bigr)\prod_{j=1}^s\tilde\sigma'_j(0)\cr
&=\exp\bigl(-2\pi i\rho_0(\gamma)\bigr)\gamma'(\tilde z_1)\prod_{j=1}^s\tilde\sigma'_j(0)\;.
\cr}
$$
Now, $\sigma'_j=(\pi'\circ\tilde\sigma_j)\tilde\sigma'_j$; therefore
$$
\eqalign{
\prod_{j=1}^s\sigma'_j(\ell_j)&=\prod_{j=1}^s\pi'\bigl(\tilde\sigma_j(\ell_j)\bigr)
\prod_{j=1}^s\tilde\sigma'_j(\ell_j)=\exp\bigl(-2\pi i\rho_0(\gamma)\bigr)\gamma'(\tilde z_1)
\prod_{j=1}^s\pi'\bigl(\tilde\sigma_j(\ell_j)\bigr)\prod_{j=1}^s\tilde\sigma'_j(0)\cr
&=\exp\bigl(-2\pi i\rho_0(\gamma)\bigr)\gamma'(\tilde z_1)
{\prod_{j=1}^s\pi'\bigl(\tilde\sigma_j(\ell_j)\bigr)\over
\prod_{j=1}^s\pi'\bigl(\tilde\sigma_j(0)\bigr)}\prod_{j=1}^s\sigma'_j(0)\cr
&=\exp\bigl(-2\pi i\rho_0(\gamma)\bigr){\pi'\bigl(\gamma(\tilde z_1)\bigr)\gamma'(\tilde z_1)
\over\pi'(\tilde z_1)}\prod_{j=1}^s\sigma'_j(0)\cr
&=\exp\bigl(-2\pi i\rho_0(\gamma)\bigr)\prod_{j=1}^s\sigma'_j(0)\;.
\cr}
$$
Now by construction, $[\gamma]=[\gamma_1]+\cdots+[\gamma_g]$ in $H_1(S,\Z)$,
where $\gamma_j$ is a counterclock-wise loop around~$p_j$; 
therefore
$$
\rho_0(\gamma)=\rho_0(\gamma_1)+\cdots+\rho_0(\gamma_g)=\sum_{j=1}^g
\Res_{p_j}(\nabla)\;,
\neweq\eqqcin
$$
and we are done. \qedn

%
%

\newdef A geodesic $\sigma\colon[0,\ell]\to S$ is {\sl closed}
if $\sigma(\ell)=\sigma(0)$ and $\sigma'(\ell)$ is a positive 
multiple of $\sigma'(0)$; it is {\sl periodic} if $\sigma(\ell)=\sigma(0)$
and $\sigma'(\ell)=\sigma'(0)$.

\newrem Contrarily to the case of Riemannian geodesics,
closed geodesics are not necessarily periodic; see Example~6.1.

\newthm Corollary \qunomez: Let $\nabla$ be a meromorphic connection on 
$\P^1(\C)$, with poles $\{p_0=\infty,p_1,\ldots,p_r\}$, and set
$S=\P^1(\C)\setminus\{p_0,\ldots,p_r\}\subseteq\C$. Let $\sigma\colon
[0,\ell]\to S$ be a geodesic with $\sigma(0)=\sigma(\ell)$ and no
other self-intersections; in particular, $\sigma$ is an oriented
Jordan curve. Let $\{p_1,\ldots,p_g\}$ be the poles of $\nabla$
contained in the interior of~$\sigma$. Then $\sigma$ is a closed geodesic if and only if
$$
\sum_{j=1}^g\Re\Res_{p_j}(\nabla)=-1\;,
$$
and it is a periodic geodesic if and only if
$$
\sum_{j=1}^g\Res_{p_j}(\nabla)=-1\;.
$$
If $\sigma$ is closed, it can be extended to an
infinite length geodesic $\sigma\colon J\to S$, where $J$ is a half-line
(possibly $J=\R$). Moreover,
\smallskip
\itm{(i)}if $\sum\limits_{j=1}^g \Im\Res_{p_j}(\nabla)<0$ then 
$\sigma'(t)\to O$ as $t\to+\infty$ and $|\sigma'(t)|\to+\infty$ as 
$t$ tends to the other end of~$J$ 
\itm{(ii)}if $\sum\limits_{j=1}^g \Im\Res_{p_j}(\nabla)>0$ then
$\sigma'(t)\to O$ as $t\to-\infty$ and 
$|\sigma'(t)|\to+\infty$ as $t$ tends to the other end of~$J$.

\newrem It might actually happen that the a closed geodesic blows up at finite 
time (that is $J\ne\R$); see Example~6.1.

\pf A self-intersecting geodesic $\sigma$ is closed if and only if the
external angle at the intersection point is~0; therefore the first assertion
follows from Corollary~\quno.

Let $\pi\colon\tilde S\to S$ be the universal covering map, and 
let $\tilde\sigma$ be a lifting of~$\sigma$,
In this case \eqqqua\ becomes
$$
\sigma'(\ell)=\exp\left(-2\pi i\sum_{j=1}^g\Res_{p_j}(\nabla)\right)\sigma'(0)\;.
$$
So $\sigma$ is periodic if and only if the sum of the residues
is an integer contained (by Corollary~\quno) in the interval~$(-3/2,-1/2)$, 
and the second assertion follows.

Finally, when $\sigma$ is closed at every turn the velocity vector is multiplied
by $\exp\left(2\pi\sum_{j=1}^g \Im\Res_{p_j}(\nabla)\right)$,
and (i) and (ii) follow.\qedn

To state our main theorem we need two more definitions.

\newdef Let $\sigma\colon I\to S$ be a curve in $S=\P^1(\C)\setminus\{p_0,
\ldots,
p_r\}$. A {\sl simple loop in~$\sigma$} is the restriction 
of $\sigma$ to a closed interval $[t_0,t_1]\subseteq I$ such that 
$\sigma|_{[t_0,t_1]}$ is a simple loop~$\gamma$. If $p_1,\ldots, p_g$
are the poles of~$\nabla$ contained in the interior of~$\gamma$,
we shall say that $\gamma$ {\sl surrounds}~$p_1,\ldots,p_g$.

\newdef A {\sl saddle connection} for a meromorphic connection $\nabla$ 
on $\P^1(\C)$ with poles $\{p_0,\ldots,p_r\}$ is a maximal geodesic
$\sigma\colon(\eps_-,\eps_+)\to S=\P^1(\C)\setminus\{p_0,\ldots,p_r\}$
(with $\eps_-\in[-\infty,0)$ and $\eps_+\in(0,+\infty]$) such that
$\sigma(t)$ tends to a pole of~$\nabla$ both when $t\uparrow\eps_+$ and
when $t\downarrow\eps_-$. A {\sl simple cycle of saddle connections} is a
Jordan piecewise smooth curve in $\P^1(\C)$ made up of saddle connections.
Again, we shall say that a simple cycle of saddle connections
{\sl surrounds} the poles of~$\nabla$ contained in its interior.

We can now prove a
Poincar\'e-Bendixson theorem for meromorphic connections on $\P^1(\C)$:

\newthm Theorem \qqua: Let $\sigma\colon[0,\eps_0)\to S$ be a maximal
geodesic for a meromorphic connection~$\nabla$ on $\P^1(\C)$, where
$S=\P^1(\C)\setminus\{p_0,\ldots,p_r\}$ and $p_0,\ldots,p_r$ are the
poles of~$\nabla$. Then either
\smallskip
\itm{(i)} $\sigma(t)$ tends to a pole of~$\nabla$ as $t\to\eps_0$; or
\itm{(ii)} $\sigma$ is closed, and then surrounds poles $p_1,\ldots,p_g$
 with $\sum\limits_{j=1}^g\Re\Res_{p_j}(\nabla)=-1$; or
\item{\rm (iii)} the $\omega$-limit set of $\sigma$ in $\P^1(\C)$ is given by the support of a closed geodesic surrounding poles
$p_1,\ldots,p_g$ with $\sum\limits_{j=1}^g\Re\Res_{p_j}(\nabla)=-1$; or
\item{\rm (iv)} the $\omega$-limit set of $\sigma$ in $\P^1(\C)$ is a simple cycle of saddle connections surrounding
poles $p_1,\ldots,p_g$ with $\sum\limits_{j=1}^g\Re\Res_{p_j}(\nabla)=-1$; or
\itm{(v)} $\sigma$ intersects itself infinitely many times, and in this
case every simple loop of $\sigma$ surrounds a set of poles whose
sum of residues has real part belonging to $(-3/2,-1)\cup(-1,-1/2)$.
\smallskip
\noindent In particular, a recurrent geodesic either intersects itself
infinitely many times or is closed.

\pf Assume that $\sigma$ is not closed, nor intersect itself 
infinitely many times (the condition on the residues in these
cases follows from Corollaries~\quno\ and \qunomez). Then up to 
changing the starting point we can assume that $\sigma$ does not
intersect itself. Let $W$ be the $\omega$-limit set of~$\sigma$ 
in~$\P^1(\C)$. Since $W$ is
connected, to end the proof it suffices to show that if $W$ 
contains a point $z_0\in S$ then we are in cases (iii) or (iv).

Let $\pi\colon\tilde S\to S$ be the universal covering map,
and $\tilde\nabla$ the holomorphic connection on~$T\tilde S$ induced
by~$\nabla$ via~$\pi$. Choose $\tilde z_0\in\pi^{-1}(z_0)$, a simply 
connected neighbourhood~$U\subset S$ of~$z_0$, and let $\tilde U$
be the connected component of~$\pi^{-1}(U)$ containing~$\tilde z_0$. 
By Proposition~\tuno.(ii), the segments of $\nabla$-geodesic contained 
in~$U$ are exactly the images under~$\pi$ of the segments
of $\tilde\nabla$-geodesic contained in~$\tilde U$. Furthermore,
by Proposition~\dzero\ up to shrinking~$U$ and $\tilde U$ we can find
an isometry $J$ between $\tilde U$ endowed with the Hermitian metric adapted
to~$\tilde\nabla$ and an open set of~$\C$ endowed with the euclidean 
metric. Therefore $\pi\circ J^{-1}$ sends line segments into segments of
$\nabla$-geodesic in~$U$, and conversely every segment of $\nabla$-geodesic
in~$U$ is image of a line segment via $\pi\circ J^{-1}$. 

So the geometry of the geodesics in a neighbourhood~$U$ of~$z_0$ is the same 
as the geometry of line segments; in particular, we can find in $U$ a (simple) geodesic $\tau$ issuing
from $z_0$ and intersecting $\sigma$ in infinitely many points
converging to~$z_0$ in $U$. Notice that all intersections are transversal
because $\sigma$ does not intersect itself.

Let $z_1$ be an intersection point between $\sigma$ and $\tau$. Following
$\sigma$ from $z_1$, let $z'_1$ be the first intersection point
between $\sigma$ and $\tau$ closer to~$z_0$ than $z_1$; let $R$ be the
Jordan domain bounded by the segments of $\sigma$ and $\tau$ between
$z_1$ and $z'_1$. By Corollary~\qdue, this domain must contain at least
one pole of~$\nabla$. 

If the two external angles of $R$ have opposite signs, then we set~$R_1=R$.
If not, we follow $\sigma$ after $z'_1$ until it meets again~$\tau$ in
a point $z''_1$ between $z_1$ and $z_0$. If $z''_1$ is between $z_1$ and
$z'_1$, then the Jordan domain bounded by $\sigma$ and $\tau$ between
$z_1$ and $z''_1$ has external angles with opposite signs, and we call
$R_1$ this domain. If instead $z''_1$ is between $z'_1$ and $z_0$
we have two possibilities. If the Jordan domain bounded by $\sigma$ and
$\tau$ between $z'_1$ and $z''_1$ has external angles with opposite signs,
we call $R_1$ this domain. If not, since by construction this domain~$R'$
is disjoint from~$R$, the poles inside $R'$ must be disjoint from the
poles inside $R$. Since the number of poles is finite, repeating
this construction sooner or later we get a Jordan domain $R_1$ whose
external angles have opposite signs.

Thus in this way we can build a sequence $\{R_j\}$ of disjoint 2-sided
geodesic polygons whose external angles have opposite signs bounded
by a segment of $\sigma$ and a segment of $\tau$ in such a way that
both vertices converge to~$z_0$. Every $R_j$ must contain poles; 
since there are only finitely many poles, up to a subsequence we can
assume they all contain the same poles $p_1,\ldots,p_g$. Since their
boundaries are disjoint, they are nested; so up to a subsequence
we can also assume that either $R_{j+1}\subset R_j$ for all $j$ or
$R_{j+1}\supset R_j$ for all $j$. 

Up to a subsequence, we can
also assume that the direction of~$\sigma$ at the vertex closest to~$z_0$ along~$\tau$
is converging to a given direction~$v_0$ in~$S^1$. Since $\sigma$ does not self-intersect,
the local geometry of the geodesics near~$z_0$ implies that the direction of~$\sigma$
at the other vertex must also converge to~$v_0$, and thus the sum of the
external angles must converge to~0. But, by Corollary~\qdue, the 
sum of the external angles is constant; so it must be zero. This
means that in each~$R_j$ the two intersections of $\sigma$ with $\tau$
are parallel, and that 
$$
\sum_{j=1}^g \Re\Res_{p_j}(\nabla)=-1\;.
$$

We have a decreasing or increasing sequence of 2-sided geodesic polygons,
with opposite external angles. If $z\in S$ is any point accumulated
by this sequence, again using the fact that $\sigma$ does not
self-intersect and the local geometrical structure of the geodesics,
we see that this sequence actually accumulates the support of a 
geodesic issuing from~$z$. Therefore $W\cap S$ is the union of
supports of finitely many disjoint geodesics. So if $W\not\subset S$ then we are in case (iv); 
if instead $W\subset S$ (for instance if the sequence $\{R_j\}$ is decreasing),
then $W$ is the support of a self-intersecting 
geodesic~$\sigma_0$
surrounding $p_1,\ldots,p_g$; hence, by Corollary~\qunomez, $\sigma_0$
is closed, and we are done.

Finally, the $\omega$-limit set of a recurrent geodesic intersects the support
of the geodesic, and the last assertion follows.\qedn

\newrem We have examples (see Examples~6.1, 8.1 and~8.2) of cases (i), (ii),
(iii) and (v), but no examples
yet of case (iv).

Using these methods we can also say something about 
self-intersecting geodesics. For instance, we can prove
the following:

\newthm Proposition \infgeo: Let $\sigma\colon[0,\eps_0)\to S$ be a 
geodesic for a meromorphic connection $\nabla$ on $\P^1(\C)$,
where $S$ is the complement in~$\P^1(\C)$ of the poles of~$\nabla$.
Assume that $\sigma$ contains two distinct simple loops 
$\sigma|_{[t_0,t_1]}$
and $\sigma|_{[t'_0,t'_1]}$ based on the same point $z_0=\sigma(t_0)=
\sigma(t_1)=\sigma(t'_0)=\sigma(t'_1)$ and 
representing the
same class $[\gamma_0]\in \pi_1(S,z_0)$. Then $\sigma$ is closed.

\pf Let $\pi\colon\tilde S\to S$ be
the universal covering map, and $\tilde\sigma\colon[0,\eps_0)\to\tilde S$
a lifting of~$\sigma$. Since
$\sigma|_{[t_0,t_1]}$ is a simple loop representing $[\gamma_0]$, we
must have $\tilde\sigma(t_1)=\gamma_0\bigl(\tilde\sigma(t_0)\bigr)$,
where $\gamma_0$ is the element of~$\Aut(\pi)$ corresponding 
to $[\gamma_0]\in\pi_1(S,z_0)$.
Then Propositions~\duno\ and \tduemez.(iii) yield
$$
\eqalign{
\rho(\gamma_0)c_0 t_0&=\rho(\gamma_0)\bigl[J\bigl(\tilde\sigma(t_0)\bigr)
-J\bigl(\tilde\sigma(0)\bigr)\bigr]=J\bigl(\gamma_0\bigl(\tilde\sigma(t_0)
\bigr)\bigr)-J\bigl(\gamma_0\bigl(\tilde\sigma(0)\bigr)\bigr)\cr
&=J\bigl(\tilde\sigma(t_1)\bigr)-J\bigl(\gamma_0\bigl(\tilde\sigma(0)\bigr)
\bigr)=c_0t_1+J\bigl(\tilde\sigma(0)\bigr)-J\bigl(\gamma_0\bigl(
\tilde\sigma(0)\bigr)\bigr)\;,
\cr}
$$
that is
$$
t_1-\rho(\gamma_0)t_0={1\over c_0}\Bigl[J\bigl(\gamma_0\bigl(\tilde\sigma
(0)
\bigr)\bigr)-J\bigl(\tilde\sigma(0)\bigr)\Bigr]\;,
$$
for a suitable $c_0\ne 0$. Repeating this argument for 
$\sigma|_{[t'_0,t'_1]}$, we get 
$$
t'_1-\rho(\gamma_0)t'_0=t_1-\rho(\gamma_0)t_0\;.
$$
So, taking the imaginary part and recalling that $t'_0\ne t_0$,
we get $\Im\rho(\gamma_0)=0$.
Since $\rho(\gamma_0)=\exp(2\pi i\rho_0(\gamma_0)\bigr)$, this implies
$\sin\bigl(2\pi\Re\rho_0(\gamma_0)\bigr)=0$, that is 
$2\pi\Re\rho_0(\gamma_0)=k\pi$ for a suitable $k\in\Z$. By \eqqcin\
and Corollary~\quno\ the only possibility is $k=-2$, that is
$\Re\rho_0(\gamma_0)=-1$. But then, by Corollary~\qunomez, $\sigma$
is necessarily closed, and we are done.\qedn


\smallsect 5. Holomorphic self-maps, homogeneous vector fields and meromorphic connections

We start this section adapting concepts introduced in [ABT1] to our
situation.

\newdef Let $f\colon M\to M$ be a holomorphic self-map of a complex $n$-dimensional
manifold~$M$, and assume that
$f$ leaves a smooth hypersurface~$S\subset M$ pointwise fixed; we write
$f\in\End(M,S)$, and always assume that $f\not\equiv\id_M$. We shall
say that a local chart $(U,z)$ of $M$, with $z=(z^1,\ldots,z^n)$, is {\sl adapted}
to~$S$ if $S\cap U=\{z^1=0\}$.

{\sc Example 5.1:} A particularly interesting example of map $f\in\End(M,S)$ is obtained blowing up a map tangent to the identity. Let $f_o$ be a (germ of) holomorphic self-map of~$\C^n$ fixing the origin and {\sl tangent to the identity,} that is such that $d(f_o)_O=\id$. If $\pi\colon M\to\C^n$ denotes the blow-up of the origin, let
$S=\pi^{-1}(O)=\P^1(\C)$ be the exceptional divisor; then the
lifting $f$ of~$f_o$, that is the unique holomorphic self-map of~$M$ such that
$f_o\circ\pi=\pi\circ f$ (see, e.g.,~[A1] for details), belongs to~$\End(M,S)$. 
\bigbreak

We denote by $N_S=TM|_S/TS$ the normal bundle of~$S$ into~$M$, by
$\ca N_S$ the sheaf of germs of holomorphic sections of~$N_S$, by $\ca T_M$ the sheaf of germs of holomorphic sections of~$TM$, and we put $\ca T_{M,S}=\ca T_M\otimes\ca O_S$, where $\ca O_S$ is the structure sheaf of~$S$. More generally, given a complex vector bundle (e.g., $E$), we shall denote by the corresponding calligraphic letter (e.g., $\ca E$) the sheaf of germs
of its holomorphic sections. 

Let $f\in\End(M,S)$ and 
take $p\in S$. Then for every
$h\in\ca O_{M,p}$ (where $\ca O_M$ is the structure sheaf of~$M$) the germ $h\circ f$ is
well-defined, and we have $h\circ f-h\in\ca I_{S,p}$, where $\ca I_S$ is the ideal sheaf 
of~$S$. 

\newdef The {\sl $f$-order of vanishing} at~$p$ of~$h\in\ca O_{M,p}$ is
$$
\nu_f(h;p)=\max\{\mu\in\N\mid h\circ f-h\in\ca I_{S,p}^\mu\}\;,
$$
and the {\sl order of contact}~$\nu_f$ of~$f$ with~$S$ is
$$
\nu_f=\min\{\nu_f(h;p)\mid h\in\ca O_{M,p}\}\;.
$$

In [ABT1] we proved that $\nu_f$ does not depend on $p$, and that 
$$
\nu_f=\min\{\nu_f(z^1;p),\ldots,\nu_f(z^n;p)\}\;,
$$
where $(U,z)$ is any local chart centered at~$p\in S$ and $z=(z^1,\ldots,z^n)$. 

\newdef A map $f\in\End(M,S)$ is {\sl tangential} to~$S$ if 
$$
\min\bigl\{\nu_f(h;p)\mid h\in\ca I_{S,p}\bigr\}>\nu_f
$$
for some (and hence any) point $p\in S$. 


Let $p\in S$, and take a chart $(U,z)$ adapted to~$S$ and centered at~$p$. If $f\in\End(M,S)$ and $f^j=z^j\circ f$, we can then write 
$$
f^j-z^j= (z^1)^{\nu_f}\,g^j\;,
\neweq\eqXfprep
$$
where $g^1,\ldots,g^n$ are holomorphic and do not all vanish when 
restricted to~$S$. They in general depend on the chosen chart; however, in [ABT1] we proved that setting
$$
\ca X_f= \sum_{j=1}^n g^j
{\de\over\de z^j}\otimes (dz^1)^{\otimes\nu_f}
\neweq\eqcaXf
$$
then $\ca X_f|_{U\cap S}$ defines a {\it global} section $X_f$ of the bundle $TM|_S\otimes (N_S^*)^{\otimes\nu_f}$, where $N_S^*$ is the conormal bundle 
of~$S$ into~$M$. 
The bundle $TM|_S\otimes (N_S^*)^{\otimes\nu_f}$ is canonically isomorphic to the bundle~$\Hom(N_S^{\otimes\nu_f},TM|_S)$; therefore the section~$X_f$ induces a morphism from $N_S^{\otimes\nu_f}$ to $TM|_S$, still denoted by $X_f$.

\newdef The morphism $X_f\colon N_S^{\otimes\nu_f}\to TM|_S$ just defined
is the {\sl canonical morphism} associated to $f\in\End(M,S)$.

It is easy to check (see [ABT1]) that $f$ is tangential if and only if the image of $X_f$ is contained in~$TS$, which amounts to saying that
$g^1|_{U\cap S}\equiv O$ for any local chart adapted to~$S$.

\newdef Assume that $f\in\End(M,S)$ is tangential.
We shall say that $p\in S$ is a {\sl singular point for~$f$} if $X_f$ vanishes at~$p$. We shall denote by $\Sing(f)$ the set of singular points for~$f$, and 
by $S^o=S\setminus\Sing(f)$ the subset of regular points. Since
$N_S^{\otimes\nu_f}$ is a line bundle, $X_f$ is a injective
on $N^{\otimes\nu_f}_{S^o}$. In particular, $X_f$ defines a rank 1
singular holomorphic foliation~$\ca F_f$ of~$S$, regular on~$S^o$. 


\newdef Assume we have a complex vector bundle $\pi_F\colon F\to S^o$ on
a complex manifold~$S^o$, and a morphism $X\colon F\to TS^o$.
A {\sl partial $X$-holomorphic connection} on a complex vector bundle $\pi_E\colon E\to S^o$ is a $\C$-linear map $\nabla\colon\ca E\to
\ca F^*\otimes\ca E$ such that
$$
\nabla_u(gs)=X(u)(g)s+g\nabla_u s
$$
for all $s\in\ca E$, $u\in\ca F$ and $g\in\ca O_{S^o}$. Clearly, if $X$ is injective we can identify $F$ with its image in~$TS^o$; in that case
we shall talk of a {\sl partial holomorphic connection along $X(F)\subset TS^o$.} Finally, if both $E$, $F$ and $X$ extend to a larger manifold~$S$, with
$S^o$ dense in $S$ and $X$ injective in~$S^o$ but not necessarily in~$S$,
we shall sometimes say that $\nabla$ is a {\sl partial meromorphic connection along $X$} 
on~$E$.

When $f$ is tangential, in [ABT1] we introduced a partial meromorphic
connection~$\nabla$ along~$X_f$ on~$N_S$ by setting
$$
\nabla_u(s)=\pi([\ca X_f(\tilde u),\tilde s]|_S)\;,
\neweq\eqCSuno
$$
where: $s\in\ca N_{S^o}$; $u\in\ca N_{S^o}^{\otimes\nu_f}$; $\pi\colon\ca T_{M,S^o}\to\ca N_{S^o}$ is the canonical projection; $\tilde s$ is any element in $\ca T_{M,S^o}$ such
that $\pi(\tilde s|_{S^o})=s$; $\tilde u$ is any
element in~$\ca T_{M,S^o}^{\otimes\nu_f}$ such that $\pi(\tilde
u|_S)=u$; and $\ca X_f$ is locally given by \eqcaXf. In a chart $(U,z)$ adapted to~$S$, a local generator of $N_{S^o}$ is $\de_1=\pi(\de/\de z^1)$, and a local generator of~$N_{S^o}^{\otimes\nu_f}$ is~$\de_1^{\otimes\nu_f}=\de_1\otimes\cdots\otimes\de_1$. 
Therefore using $\tilde u=(\de/\de z^1)^{\otimes\nu_f}$ as extension of~$\de_1^{\otimes\nu_f}$, and $\tilde s=\de/\de z^1$ as extension of~$\de_1$ we get
$$
\nabla_{\de_1^{\otimes\nu_f}}\de_1=-\left.{\de g^1\over\de z^1}
\right|_{U\cap S^o}\de_1\;.
$$

Up to here we limited ourselves to summarize [ABT1]; let us now introduce
new ideas.

A partial meromorphic connection~$\nabla$ along $X_f$ on~$N_S$ canonically induces a partial meromorphic connection (still denoted by~$\nabla$) along~$X_f$ on~$N_S^{\otimes\nu_f}$ by setting
$$
\nabla(s_1\otimes\cdots\otimes s_{\nu_f})=\sum_{j=1}^{\nu_f} s_1\otimes\cdots\otimes
\nabla s_j\otimes\cdots\otimes s_{\nu_f}\;.
$$
In particular we get 
$$
\nabla_{\de_1^{\otimes\nu_f}}\de_1^{\otimes\nu_f}=-\nu_f\left.
{\de g^1\over\de z^1}\right|_{U\cap S^o}\de_1^{\otimes\nu_f}\;.
$$

As remarked before, the morphism $X_f$ defines a rank 1 singular holomorphic foliation $\ca F_f$ on~$S$, locally generated by
$$
v_o=X_f(\de_1^{\otimes\nu_f})=\sum_{p=2}^n g^p|_{U\cap S}{\de\over\de z^p}
\;.
\neweq\eqvzero
$$
We can then define a partial meromorphic connection $\nabla^o\colon\ca F_f
\to\ca F_f^*\otimes\ca F_f$ along the identity on~$\ca F_f$, holomorphic
on~$S^o$, by setting
$$
\nabla^o_v s=X_f\bigl(\nabla_{X_f^{-1}(v)}X_f^{-1}(s)\bigr)\;.
$$
Notice that, by construction, {\it $\nabla^o$ induces a (standard) holomorphic connection on each leaf of the foliation~$\ca F_f$;} so the
geodesics we shall introduce momentarily will be {\it geodesics for a
holomorphic connection on a Riemann surface,} that is exactly of the kind
we have studied in the first part of this paper.

\newrem When $n=2$, the morphism $X_f$ is an isomorphism between $N_{S^o}^{\otimes\nu_f}$
and $TS^o$; so $\nabla^o$ is a standard meromorphic connection on~$S$.
In particular, locally we have $\de/\de z^2={1\over g^2}v_0$, and thus
$\nabla^o$ is represented by the 1-form
$$
\eta=-\left.\nu_f{1\over g^2}{\de g^1\over\de z^1}\right|_{U\cap S^o}dz^2\;.
$$

\newdef A {\sl $\nabla^o$-geodesic} is a (real) curve $\sigma\colon I\to S^o$ such that $\sigma'(t)\in(\ca F_f)_{\sigma(t)}$ for all $t\in I$ (that is,
the image of~$\sigma$ is contained in a leaf of the foliation) and
$\nabla^o_{\sigma'}\sigma'\equiv O$.  

In local coordinates $(U,z)$ adapted to~$S$, saying that the image of 
$\sigma$ is contained in a leaf of the foliation~$\ca F_f$ is 
equivalent to saying that $\sigma'$ is a multiple of the generator~$v_0$ introduced in~\eqvzero. In other words, writing $\sigma^j=z^j\circ\sigma$
and denoting, with a slight abuse of notation, by $\sigma'_o$ this multiple, we have that the image of $\sigma$ is contained in a leaf of the foliation 
if and only if 
$$
(\sigma^j)'=\sigma'_o\cdot (g^j\circ\sigma)
\neweq\eqcunob
$$
for $j=1,\ldots,n$ (in particular, $\sigma^1\equiv 0$). Furthermore, $\sigma$ is a
$\nabla^o$-geodesic if and only $X_f^{-1}(\sigma')=\sigma'_o\de_1^{\otimes\nu_f}$ with
$$
(\sigma'_o)'-\nu_f\left({\de g^1\over\de z^1}\circ\sigma\right)(\sigma'_o)^2=0\;.
\neweq\eqcdueb
$$
This suggests to introduce a holomorphic vector field $G$ defined on the
total space of $p\colon N^{\otimes\nu_f}_S\to S$ by setting
$$
G|_{p^{-1}(U)}=\sum_{p=2}^n g^p|_{U\cap S}\,v{\de\over\de z^p}+\nu_f
\left.{\de g^1\over\de z^1}\right|_{U\cap S}v^2{\de\over\de v}\;,
\neweq\eqctreb
$$
where $v$ is the fiber coordinate corresponding to the generator~$\de_1^{\otimes\nu_f}$
(and $v^2$ is the square of the coordinate $v$). 

\newthm Proposition \geofb: Let $f\in\End(M,S)$ be tangential. Then:
\smallskip
\itm{(i)} the formula $\eqctreb$ defines a global holomorphic vector field on the total space of $N^{\otimes\nu_f}_S$;
\itm{(ii)} a curve $\sigma\colon I\to S^o$ is a $\nabla^o$-geodesic if and
only if the image of $\sigma$ is contained in a leaf of~$\ca F_f$ and
$X_f^{-1}(\sigma')$\break\indent is an integral curve of~$G$. 

\pf (i) follows from a not too difficult computation (using, e.g., [ABT1, (3.6) and (4.2)]), while (ii) follows from \eqcunob\ and \eqcdueb.\qedn

As mentioned before, our main example is when $S$ is the exceptional divisor of the
blow-up of the origin in $\C^n$, and $f$ is the lifting of a germ tangent to the identity.
Let us now discuss some peculiar features of this case.

Let $\pi\colon M\to\C^n$ be the blow-up of the origin in~$\C^n$, and $S=\pi^{-1}(O)=\P^{n-1}(\C)$ the exceptional divisor. Let $w=(w^1,\ldots,w^n)$ denote coordinates in~$\C^n$, and set $H_j=\{w\in\C^n\mid w^j\ne 0\}\subset\C^n$
for $j=1,\ldots,n$. We can cover $M$ with $n$ charts $(U_j,z_j)$, where
$U_j=\pi^{-1}(H_j)$ for $j=1,\ldots,n$; the chart $(U_j,z_j)$ is
centered in $[0:\cdots:1:\cdots:0]\in\P^{n-1}(\C)$, and $U_j\cap S=\{z^j_j=0\}$, where $z_j=(z_j^1,\ldots,z_j^n)$. The projection~$\pi$ on~$U_j$
is given by
$$
\pi\bigl(z_j(p)\bigr)=(z^1_jz^j_j,\ldots,z^j_j,\ldots,z_j^nz_j^j)\;;
$$
and the coordinate changes in $U_i\cap U_j=\{z_i^j,z_j^i\ne 0\}$ are 
given by
$$
z_j^h=\cases{z^i_iz^j_i&for $h=j\;$;\cr
1/z^j_i& for $h=i\;$;\cr
z^h_i/z^j_i& for $h\ne i$,~$j\;$;\cr}
$$
see [A1] for details. It follows that tangent vectors and covectors change 
according to the following rules:
$$
dz^h_j=\cases{z^i_i\,dz^j_i+z^j_i\,dz^i_i&for $h=j\;$;\cr
-{1\over(z_i^j)^2}\,dz^j_i&for $h=i\;$;\cr
{1\over z^j_i}\,dz^h_i-{z^h_i\over(z_i^j)^2}\,dz^j_i&for $h\ne i$,~$j\;$;
\cr}
$$
and
$$
{\de\over\de z^h_j}=\cases{\displaystyle{1\over z^j_i}{\de\over\de z^i_i}&for $h=j\;$;\cr
\displaystyle z^j_i\left(2z^i_i{\de\over\de z^i_i}-\sum_{k=1}^n z_i^k
{\de\over\de z^k_i}\right)&for $h=i\;$;\cr
\displaystyle z_i^j{\de\over\de z^h_i}&for $h\ne i$,~$j\;$.
\cr}
$$
We shall denote by $(\zeta_j,v_j)$ the induced coordinates on~$N_S^{\otimes\nu}$, where
$$
\zeta_j=(\zeta_j^1,\ldots,\zeta_j^{n-1})=(z_j^1,\ldots,\widehat{z_j^j},
\ldots,z_j^n)\in\C^{n-1}\;.
$$
The coordinate changes in $N_S^{\otimes\nu}$ are then given by
$$
\zeta_j^h=\cases{\zeta^h_i/\zeta^j_i& for $1\le h\le j-1$ and $i\le h\le n-1\;$,\cr
\zeta^{h+1}_i/\zeta^j_i&for $j\le h\le i-2\;$,\cr
1/\zeta^j_i&for $h=i-1\;$,
\cr}\quad\hbox{and}\quad v_j=(\zeta_i^j)^\nu v_i
\neweq\eqcnoveb
$$
when $j<i$, and by similar formulas when $j>i$.

The first consequence of these formulas is the following

\newthm Proposition \cdue: Let $\pi\colon M\to\C^n$ be the blow-up of the origin in~$\C^2$, and let $S=\pi^{-1}(O)=\P^{n-1}(\C)$ be the exceptional divisor. Then for any~$\nu\in\N^*$ we can define a $\nu$-to-$1$
holomorphic covering map $\chi_\nu\colon\C^n\setminus\{O\}\to N^{\otimes\nu}_S\setminus S$ by setting
$$
\zeta_j(w)=\left({w^1\over w^j},\ldots,\widehat{{w^j\over w^j}},\ldots,
{w^n\over w^j}\right)\qquad\hbox{and}\qquad v_j(w)=(w_j)^\nu
$$
for all $w\in H_j$ and $j=1,\ldots,n$. In particular,
$p\circ\chi_\nu
\colon\C^n\setminus\{O\}\to\P^{n-1}(\C)$ is the canonical projection, where
$p\colon N_S^{\otimes\nu}\to S=\P^{n-1}(\C)$ is the projection.

\pf It suffices to remark that $\chi_\nu$ is well-defined, thanks to 
\eqcnoveb.\qedn

\newdef We shall call the map $\chi_\nu\colon\C^n\setminus\{O\}\to
N_{\P^{n-1}(\C)}^{\otimes\nu}\setminus \P^{n-1}(\C)$ just defined the {\sl $\nu$-polar coordinates} of~$\C^n$. Furthermore, we shall denote by
$[\cdot]\colon\C^n\setminus\{O\}\to\P^{n-1}(\C)$ the canonical projection; so
$$
p\circ\chi_\nu(w)=[w]
$$
for all $w\in\C^n\setminus\{O\}$, where $p\colon N_{\P^{n-1}(\C)}^{\otimes\nu}
\to\P^{n-1}(\C)$ is the projection.

%

Let us now $f\in\End(M,S)$ be obtained blowing-up a germ~$f_o\in\End(\C^n,O)$ tangent to the identity, as in Example~5.1. Write
$$
f_o^j(w)=w^j+\sum_{h\ge\nu+1}Q_h^j(w)\;,
\neweq\eqcfz
$$
where $Q_h^j$ is a homogeneous polynomial of degree~$h$, and $\nu+1\ge 2$ is the {\sl order}~$\nu(f_o)$ of~$f_o$, chosen so that $(Q_{\nu+1}^1,\ldots,
Q_{\nu+1}^n)\not\equiv O$. We associate to~$f_o$ the homogeneous
vector field of degree $\nu+1$
$$
Q_{f_o}=\sum_{j=1}^n Q^j_{\nu+1}{\de\over\de w^j}\;.
$$

\newdef We say that a homogeneous vector field $Q$ is {\sl dicritical}
if it is a multiple of the {\sl radial vector field}
$$
\sum_{j=1}^n w^j{\de\over\de w^j}\;.
$$
In other words, $Q=\sum_j Q^j{\de\over\de w^j}$ is dicritical
if and only if
$$
w^hQ^k\equiv w^kQ^h
$$
for all $h$, $k=1,\ldots, n$.
A map~$f_o$ tangent to the identity is {\sl dicritical} if $Q_{f_o}$ is.

Let $\pi\colon M\to\C^n$ be the blow-up of the origin, and $f\in\End\bigl(M,
\P^{n-1}(\C)\bigr)$ the lift of~$f_o$ to th blow-up. Then in the 
chart $(U_1,z_1)$ introduced before setting $f_1=z_1\circ f$ we get
$$
f_1^j(z_1)=\cases{\displaystyle z_1^1+(z_1^1)^\nu\sum_{h\ge\nu+1}
(z_1^1)^{h-\nu-1}Q_h^1(1,\zeta_1)&for $j=1\;$,\cr
z_1^j+(z_1^1)^\nu\,{\sum\limits_{h\ge\nu+1}(z_1^1)^{h-\nu-1}
\bigl[Q_h^j(1,\zeta_1)-z_1^jQ_h^1(1,\zeta_1)\bigr]\over
1+(z_1^1)^\nu\sum\limits_{h\ge\nu+1}(z_1^1)^{h-\nu-1}Q_h^1(1,\zeta_1)}
&for $j>1\;$;
\cr}
$$
similar formulas hold in the other charts. In particular, it follows that
\smallskip
\item{(i)} if $f_o$ is non-dicritical then $f$ is tangential to the exceptional divisor $\P^{n-1}(\C)$ and $\nu_f=\nu(f_o)-1$;
\item{(ii)} if $f_o$ is dicritical then $f$ is not tangential and $\nu_f=\nu(f_o)$. 
\smallskip
\noindent Thus {\it most maps constructed with this procedure are tangential.}

Assume then that $f_o$ is non-dicritical, so that $f$ is tangential and
$\nu_f=\nu(f_o)-1=\nu$. Then in the canonical chart $(U_1,z_1)$ we have
$$
\displaylines{
\left.{\de g^1_1\over\de z^1_1}\right|_{U_1\cap S}=Q^1_{\nu+1}(1,\zeta_1)\;,
\qquad g^p_1|_{U_1\cap S}=Q^p_{\nu+1}(1,\zeta_1)-\zeta_1^{p-1}
Q^1_{\nu+1}(1,\zeta_1)\hbox{ for $p=2,\ldots,n\;$,}\cr
G|_{p^{-1}(U_1\cap S)}=\sum_{h=1}^{n-1}\bigl[Q_{\nu+1}^{h+1}(1,\zeta_1)
-\zeta_1^hQ_{\nu+1}^1(1,\zeta_1)\bigr]v_1{\de\over\de\zeta_1^h}
+\nu Q^1_{\nu+1}(1,\zeta_1)v_1^2{\de\over\de v_1}\;,\cr}
$$
and similar formulas hold in the other charts. In particular, it follows
that the canonical morphism~$X_f$ and the
connection~$\nabla$ (and hence the connection $\nabla^o$ and the $\nabla^o$-geodesics) depend only on the homogeneous vector field~$Q_{f_o}$. Thus
{\it we can use the same formulas to associate to {\rm any} non-dicritical
homogeneous vector field~$Q$ of degree~$\nu+1$ the canonical morphism
$X_Q\colon N^{\otimes\nu}_{\P^{n-1}(\C)}\to T\P^{n-1}(\C)$, the meromorphic
connection $\nabla$ on $N^{\otimes\nu}_{\P^{n-1}(\C)}$, and the geodesic field $G$.} In other words,
we associate to $Q=\sum_j Q^j{\de\over\de z^j}$ all the objects we would get 
starting from the time-1 map~$f_Q$ of $Q$, which is of the form
$$
f_Q^j(z)=z^1+Q^j(z)+O(\|z\|^{\nu+2})
$$
for $j=1,\ldots,n$.

To describe another consequence of these formulas we need another definition.

\newdef A {\sl characteristic direction} of a homogeneous vector field
$Q=\sum_j Q^j{\de\over\de w^j}$ is a direction $[v]\in\P^{n-1}(\C)$
such that the line $L_v=\C v$, where $v\in\C^n\setminus\{O\}$ is any representative of~$[v]$, is $Q$-invariant; in this case we shall say that
$L_v$ is a {\sl characteristic leaf} of~$Q$. If moreover $Q|_{L_v}\equiv O$
we say that $[v]$ is {\sl degenerate;} otherwise, it is {\sl non-degenerate.} If $S^o\subset\P^{n-1}(\C)$ is the complement of the 
characteristic directions of~$Q$, we shall write $\hat S_Q=\{w\in\C^n
\setminus\{O\}\mid [w]\in S^o\}$.

It is easy to see that $[v]\in\P^{n-1}(\C)$ is characteristic if and only if
$\bigl(Q^1(v),\ldots,Q^n(v)\bigr)=\lambda v$ for some $\lambda\in\C$ (clearly depending on the representative $v$ of~$[v]$), and that $[v]$
is non-degenerate if and only if $\lambda\ne 0$. Then it is clear that
{\it the singular points of $X_Q$ are exactly the characteristic directions
of~$Q$.}

This is just the first signal that we can relate the dynamics of~$Q$ to
the geodesic flow of~$\nabla^o$. And indeed we have the following:

\newthm Theorem \ctreb: Let $Q$ be a non-dicritical homogeneous vector field of degree~$\nu+1\ge 2$ in~$\C^n$. Then
$$
d\chi_\nu(Q)=G\;,
\neweq\eqcstarc
$$
where $G$ is the geodesic field on $N^{\otimes\nu}_{\P^{n-1}(\C)}$ 
associated to~$Q$.
In particular
for a real curve
$\gamma\colon I\to\hat S_Q$ the following are equivalent:
\itemitem{\rm(i)} $\gamma$ is an integral curve of~$Q$;
\itemitem{\rm(ii)} $\chi_\nu\circ\gamma$ is an integral curve of~$G$;
\itemitem{\rm(iii)} $[\gamma]$ is a $\nabla^o$-geodesic.

\pf \eqcstarc\ follows from the homogeneity of~$Q$ and the formula
$$
d\chi_\nu\left({\de\over\de w^h}\right)=\cases{\displaystyle
{1\over w^j}{\de\over\de\zeta_j^h}&if $h<j\;$,\cr
\displaystyle{1\over w^j}{\de\over\de\zeta_j^{h-1}}& if $h>j\;$,\cr
\displaystyle-{1\over(w^j)^2}\left[\sum_{k=1}^{j-1}w^k{\de\over\de\zeta_j^k}
+\sum_{k=j}^{n-1}w^{k+1}{\de\over\de\zeta_j^k}\right]+\nu(w^j)^{\nu-1}
{\de\over\de v_j}&if $h=j\;$.
\cr}
$$
The last assertion then follows from Proposition~\geofb.\qedn

In particular, if $\gamma$ is an integral curve of~$Q$ than $[\gamma]$ 
necessarily belongs to a leaf $L$ of the holomorphic singular
foliation~$\ca F_Q$ of $\P^{n-1}(\C)$ induced by the canonical morphism~$X_Q$, and it is a geodesic for a meromorphic connection on~$L$. Thus
the study of the dynamics of~$Q$ boils down to the study of the singular
foliation~$\ca F_Q$ of $\P^{n-1}(\C)$ and of the geodesic flow of
meromorphic connections on Riemann surfaces. 

To show the power of this
approach, from the next section on we shall discuss what happens in 
dimension~2, where the foliation~$\ca F_Q$ is trivial; but we end 
this section describing the dynamics of dicritical vector fields, the
only case our approach does not work. 

Actually, the dynamics of a dicritical vector field is very easy to study, 
because all directions are characteristic and the dynamics inside
a characteristic leaf is one-dimensional, as shown by the following

\newthm Lemma \dcl: Let $L_v=\C v$ be a characteristic leaf of a
homogeneous vector field $Q$ of degree $\nu+1\ge 2$ in~$\C^n$. Then:
\smallskip
\itm{(i)} if $[v]\in\P^1(\C)$ is a degenerate characteristic direction
then the dynamics of $Q$ on $L_v$ is trivial;
\itm{(ii)} if $[v]\in\P^1(\C)$ is a non-degenerate characteristic 
direction, then the integral curve of $Q$ issuing from $\zeta_0 v\in L_v$
is given by
$$
\gamma_{\zeta_0 v}(t)={\zeta_0 v\over(1-\lambda_0\zeta_0^\nu \nu t)^{1/\nu}}\;,
\neweq\eqintcd
$$
where $\lambda_0\in\C^*$ is such that $Q^j(v)=\lambda_0 v^j$ for $j=1,\ldots,n$.
In particular no (non-constant) integral curve is recurrent, and we have:
\smallskip
\itemitem{\rm (a)} if $\lambda_0\zeta_0^\nu\notin\R^+$ then $\lim\limits_{t\to+\infty}
\gamma_{\zeta_0 v}(t)=O$;
\itemitem{\rm (b)} if $\lambda_0\zeta_0^\nu\in\R^+$ then $\lim\limits_{t\to(\lambda_0\zeta_0^\nu\nu)^{-1}}\|\gamma_{\zeta v}(t)\|=+\infty$.

\pf Part (i) is clear. For part (ii), let
$\phe(\zeta)=\zeta v$ be a parametrization of~$L_v$. Then
$$
d\phe^{-1}(Q|_{L_v})=\lambda_0 \zeta^{\nu+1}{\de\over\de\zeta}\;.
$$
The integral curves of this one-dimensional vector field are
$$
\zeta(t)={\zeta_0\over(1-\lambda_0\nu\zeta_0^\nu t)^{1/\nu}}\;,
$$
where the determination of the $\nu$-th root is chosen so that $\zeta(0)=
\zeta_0$; therefore the integral curve of $Q$ issuing from $\zeta_0 v$ is
given by \eqintcd.\qedn

\smallsect 6. Global dynamics in dimension~2

Let $f\in\End(M,S)$ be tangential, where $M$ is a 2-dimensional complex manifold
and $S$ a 1-dimensional complex submanifold of~$M$. Then
the partial connection $\nabla$ introduced in the previous section is a 
{\it bona fide} holomorphic connection on $N_{S^o}^{\otimes\nu_f}$, and
we have a isomorphism $X_f\colon N_{S^o}^{\otimes\nu_f}\to TS^o$. 
We then are in the situation described in Section~1, with $E=N_{S^o}^{\otimes\nu_f}$ and $X=X_f$. In particular, we get
\smallskip
\item{--} the metric foliation on $N_{S^o}^{\otimes\nu_f}\setminus S^o$, a real non-singular foliation of real rank~3;  
\item{--} the horizontal foliation on $N_{S^o}^{\otimes\nu_f}$, a complex non-singular foliation of complex rank~1;
\item{--} the geodesic foliation on $N_{S^o}^{\otimes\nu_f}$, a real foliation of real rank~1, singular only on the zero section.
\smallskip
\noindent Furthermore, the holomorphic connection~$\nabla^o$ 
induced by $\nabla$ via $X_f$ on $TS^o$ is a meromorphic connection on~$S$,
and the geodesic field $G$ introduced in the previous section coincides
with the geodesic field introduced in Section~1.  
Locally, if $(U_\alpha,z_\alpha)$ is a chart adapted to~$S$, then $z^2_\alpha$ is a local coordinate on~$S$, and $\de/\de z_\alpha^2$ is a local generator of~$TS$; hence we have
$$
X_\alpha=g^2_\alpha|_{U_\alpha\cap S^o}\;,\quad
\nabla^o_{\de/\de z^2_\alpha}{\de\over\de z^2_\alpha}=-
\nu_f{1\over g^2_\alpha}\left.{\de g^1_\alpha\over\de z^1_\alpha}\right|_{U_\alpha\cap S^o}\quad\hbox{and}\quad \eta_\alpha=-\left.{\nu_f\over g^2_\alpha}\,{\de g^1_\alpha\over\de z^1_\alpha}\right|_{U_\alpha\cap S^o}\,dz^2_\alpha\;,
\neweq\eqkalfa
$$
and
$$
G=X_\alpha v_\alpha H_\alpha=g^2_\alpha|_{U_\alpha\cap S}\, v_\alpha\de_\alpha+\nu_f\left.{\de g^1_\alpha\over\de z^1_\alpha}\right|_{U_\alpha\cap S} v_\alpha^2\,{\de\over\de v_\alpha}\;.
\neweq\eqGext
$$
In particular, $G$ vanishes only on the zero section and on the fibers over those singular points where $\de g^1_\alpha/\de z^1_\alpha|_S$ vanishes too. 

In particular, $G$ defines a singular extension of the horizontal foliation to the whole of~$N_S^{\otimes\nu_f}$. We can use \eqGext\ to study the associated saturated foliation. Indeed, assume that the chart $(U_\alpha,z_\alpha)$ is centered at a singular point~$p\in U_\alpha\cap S$. Then on $U_\alpha\cap S$ we can write
$$
X_\alpha v_\alpha H_\alpha=v_\alpha(z^2_\alpha)^\mu\left[h^2_\alpha\, \de_\alpha+\nu_f h^1_\alpha v_\alpha\,{\de\over\de v_\alpha}\right]\;,
$$
with $h^1_\alpha$, $h^2_\alpha$ holomorphic functions on~$U_\alpha\cap S$ not both vanishing at~$p$. Then the section in square brackets generates the saturation of the horizontal foliation, and it is vanishing only when $h^2_\alpha(p)=0$ and $v_\alpha=0$, that is in the singular points for $f$ in the zero section where $g^2_\alpha$ vanishes at a higher order than $\de g^1_\alpha/\de z^1_\alpha$.

\newdef Let $f\in\End(M,S)$ be tangential, where $M$ is a 2-dimensional complex manifold
and $S$ a 1-dimensional complex submanifold of~$M$. The {\sl order}~$\mu_p$ of a point
$p\in S$ is the order of vanishing of~$X_f$ at~$p$. In a local
chart $(U_\alpha,z_\alpha)$ adapted to~$S$, we have
$$
\mu_p=\ord_p(g^2_\alpha|_S)\;.
$$
In particular, $p\in\Sing(f)$ if and only if~$\mu_p\ge 1$.
We say that $p\in \Sing(f)$ is an 
{\sl apparent singularity} if it is not a pole of~$\nabla$, that is if
$\mu_p\le\ord_p(\de g^1_\alpha/\de z^1_\alpha)$ for a
(and hence any) local chart $(U_\alpha,z_\alpha)$ adapted to~$S$. 
Furthermore, we shall say that $p$ is a
{\sl Fuchsian singularity} if $\mu_p=\ord_p(\de g^1_\alpha/\de z^1_\alpha)+1$, and that is an {\sl irregular singularity} if 
$\mu_p>\ord_p(\de g^1_\alpha/\de z^1_\alpha)+1$. Finally,
we say that $p$ is a {\sl strictly fixed point} if it is a pole of~$\nabla$
(that is, a Fuchsian or irregular singularity), and that it is
a {\sl degenerate  singularity} if $\mu_p$,~$\ord_p(\de g^1_\alpha/\de 
z^1_\alpha)\ge 1$. Finally, the {\sl index} of $p\in S$ is given by
$$
\iota_p(f,S)=-{1\over\nu_f}\Res_p(\nabla)\;.
$$
It is not difficult to check (see [A2, ABT1] and Section~7)
that these definitions do not depend on the adapted local chart. 

\newrem The index defined here is the same one introduced in [A2] and [ABT1]. In particular, there we proved that if $S$ is compact then
$$
\sum_{p\in S}\iota_p(f,S)=\int_S c_1(N_S)\;,
$$
where $c_1(N_S)$ is the first Chern class of the normal bundle~$N_S$.

Using these definitions we can say that
{\it the geodesic field $G$ extends to a holomorphic section of $T(N_S^{\otimes\nu_f})$, vanishing only on the zero section and on the fibers over degenerate singularities; and
the horizontal foliation extends to a saturated rank~$1$ complex foliation of ~$N_S^{\otimes\nu_f}$, whose singular points are the strictly fixed points for $f$ in the zero section.}

\newrem If we replace the map $f$ by its inverse $f^{-1}$, it is easy to see that $X_\alpha$ changes sign whereas $\eta_\alpha$ does not change. Therefore $f$ and $f^{-1}$ induce the same metric and horizontal foliations, but the geodesic field of~$f$ is the opposite of the geodesic field of~$f^{-1}$ --- and thus the orientation of the leaves of the geodesic foliations of~$f$ is opposite to the orientation of the leaves of the geodesic foliations 
of~$f^{-1}$, even though the leaves are the same.

Let $\pi\colon M\to\C^2$ be the blow-up of the origin in~$\C^2$, and $S=\pi^{-1}(O)=\P^1(\C)$ the exceptional divisor.
In this case we have an atlas of~$M$ adapted to~$S$ composed by two charts $(U_0,z_0)$
and~$(U_\infty,z_\infty)$, with $z_0(U_0)=z_\infty(U_\infty)=\C^2$. The chart~$(U_0,z_0)$ is centered at the
point~$0=[1:0]\in\P^1(\C)$ of the exceptional divisor, while the chart~$(U_\infty,z_\infty)$
is centered at the point~$\infty=[0:1]\in\P^1(\C)$ of the exceptional divisor. If we denote
by~$(w^1,w^2)$ the coordinates of~$\C^2$, the projection~$\pi$ from~$M$ onto~$\C^2$ is given by
$$
\cases{w^1=z^1_0\;,\cr
w^2=z^1_0z^2_0\;,\cr}\ \hbox{on $U_0\;$,\qquad and by\qquad}
\cases{w^1=z^1_\infty z^2_\infty\;,\cr
w^2=z^1_\infty\;,\cr}\ \hbox{on $U_\infty\;$.}
$$
In particular, $\pi(U_0)=\C^2\setminus\{(0,w^2)\mid w^2\ne0\}$ and $\pi(U_\infty)=\C^2
\setminus\{(w^1,0)\mid w^1\ne0\}$, and the cocycles of $TS$ and $N_S$ are respectively given by
$$
\psi_{0\infty}(z^2_0)=-(z_0^2)^2\qquad\hbox{and}\qquad \xi_{0\infty}(z^2_0)={1\over
z_0^2}\;.
\neweq\eqtdue
$$
Furthermore, if $(\zeta_0,v_{0,\nu})$ and $(\zeta_\infty,v_{\infty,\nu})$ are the local coordinates on the total space of~$N_S^{\otimes\nu}$ for $\nu\in\N^*$ induced by the canonical charts~$(U_0,z_0)$ and $(U_\infty,z_\infty)$ of~$M$, we have
$$
\zeta_\infty={1\over\zeta_0}\qquad\hbox{and}\qquad v_{\infty,\nu}={1\over\xi_{0\infty}(\zeta_0)^\nu}\, v_{0,\nu}=\zeta_0{}^\nu v_{0,\nu}\;.
\neweq\eqctre
$$
%
%

Let now $Q=Q^1{\de\over\de w^1}+Q^2{\de\over\de w^2}$ be a non-dicritical homogeneous 
vector field of degree~$\nu+1\ge 2$. Using as~$f\in\End(M,S)$ the 
blow-up of the time-1 map~$f_Q$ of~$Q$, and 
%
%
%
writing for simplicity $v_0$ and $\zeta_0$ instead of $v_{0,\nu}$ and $z_0^2$, and $v_\infty$ and $\zeta_\infty$ instead of $v_{\infty,\nu}$ and $z_\infty^2$, we get
$$
\displaylines{
\left.{\de g^1_0\over\de z^1_0}\right|_{U_0\cap S}=Q^1(1,\zeta_0)\;,\qquad
g^2_0|_{U_0\cap S}=Q^2(1,\zeta_0)-\zeta_0 Q^1(1,\zeta_0)\;,\cr
X_0=Q^2(1,\zeta_0)-\zeta_0 Q^1(1,\zeta_0)\;,\qquad 
\eta_0=-{\nu\, Q^1(1,\zeta_0)\over Q^2(1,\zeta_0)-\zeta_0 Q^1(1,\zeta_0)}\,
d\zeta_0\;,\cr
\omega|_{p^{-1}(U_0\cap S^o)}={-\nu\,Q^1(1,\zeta_0)\over
Q^2(1,\zeta_0)-\zeta_0 Q^1(1,\zeta_0)}\,d\zeta_0+{1\over v_0}\,dv_0\;,\cr
H_0={\de\over\de\zeta_0}+{\nu\,Q^1(1,\zeta_0)\over 
Q^2(1,\zeta_0)-\zeta_0 Q^1(1,\zeta_0)}\,v_0{\de\over\de v_0}\;,\cr
G|_{p^{-1}(U_0\cap S^o)}=\bigl(Q^2(1,\zeta_0)-\zeta_0 Q^1(1,\zeta_0)\bigr)v_0
{\de\over\de\zeta_0}+\nu\,Q^1(1,\zeta_0)(v_0)^2{\de\over\de v_0}\;,
\cr}
$$
and
$$
\displaylines{
\left.{\de g^1_\infty\over\de z^1_0}\right|_{U_\infty\cap S}=Q^2(\zeta_\infty,1)\;,\qquad
g^2_\infty|_{U_\infty\cap S}=Q^1(\zeta_\infty,1)-\zeta_\infty Q^2(\zeta_\infty,1)\;,\cr
X_\infty=Q^1(\zeta_\infty,1)-\zeta_\infty Q^2(\zeta_\infty,1)\;,\qquad
\eta_\infty=-{\nu\,Q^2(\zeta_\infty,1)\over 
Q^1(\zeta_\infty,1)-\zeta_\infty Q^2(\zeta_\infty,1)}\,d\zeta_\infty\;,\cr
\omega|_{p^{-1}(U_\infty\cap S^o)}={-\nu\,Q^2(\zeta_\infty,1)\over 
Q^1(\zeta_\infty,1)-\zeta_\infty Q^2(\zeta_\infty,1)}\,d\zeta_\infty+{1\over v_\infty}\,dv_\infty\;,\cr 
H_\infty={\de\over\de\zeta_\infty}+{\nu\,Q^2(\zeta_\infty,1)\over Q^1(\zeta_\infty,1)-\zeta_\infty Q^2_(\zeta_\infty,1)}\,v_\infty{\de\over\de v_\infty}\;,\cr
G|_{p^{-1}(U_\infty\cap S^o)}=\bigl(Q^1(\zeta_\infty,1)-\zeta_\infty Q^2(\zeta_\infty,1)\bigr)v_\infty
{\de\over\de\zeta_\infty}+\nu Q^2(\zeta_\infty,1)(v_\infty)^2{\de\over\de v_\infty}\;.
\cr}
$$

If $[v]\in\P^1(\C)$, the {\sl index}~$\iota_{[v]}
(Q)$ of~$Q$ at $[v]$ is the index introduced in Definition~6.1, that is
$$
\iota_{[v]}(Q)=-{1\over\nu}\Res_{[v]}(\nabla)\;.
$$
In particular, if $[v]=[1:v_0]$ then 
$$
\iota_{[v]}(Q)=-{1\over\nu}\Res_{[1:v_0]}(\eta_0)=\Res_{v_0}\left(
{Q^1(1,\zeta)\over Q^2(1,\zeta)-\zeta Q^1(1,\zeta)}\right)\;;
\neweq\eqindex
$$
a similar formula using $\eta_\infty$ yields the index of $Q$ at $[0:1]$.
Clearly, the index at $[v]$ can be different from zero only if $[v]$ is
a characteristic direction.

\newrem If $[v]=[1:v_0]$ is a non-degenerate characteristic direction, a related 
number is the
{\sl director}~$\delta_{[v]}(Q)$ of~$[v]$ given by (see [\'E1--4] and [H1])
$$
\delta_{[v]}(Q)={1\over Q^1(1,v_0)}{\de\bigl(Q^2(1,\zeta)-\zeta Q^1(1,\zeta)\bigr)\over\de \zeta}(v_0)\;.
$$
A similar formula yields the director of $\infty=[0:1]$ when the latter is a
characteristic direction. 

\newrem A homogeneous vector field of degree $\nu+1$ is dicritical if
and only if every direction is characteristic. If $Q$ is non-dicritical
then it has $\nu+2$ characteristic directions, counted with 
multiplicity (which is just the order introduced in Definition~6.1; see [AT1]). 

\newrem
Comparing definitions we immediately see that:
\smallskip
\item{(a)} $[v]\in\P^1(\C)$ is a characteristic direction for $Q$ if 
and only if it is a singular point of $X_Q$;
\item{(b)} $[v]\in\P^1(\C)$ is a degenerate characteristic direction
for $Q$ if and only if it is a degenerate singularity of~$X_Q$;
\item{(c)} if $[v]\in\P^1(\C)$ is a non-degenerate characteristic direction
then it is a strictly fixed point of $X_Q$;
\item{(d)} $[v]\in\P^1(\C)$ is a non-degenerate characteristic direction
with non-vanishing director if and only if it is a Fuchsian singularity of~$X_Q$ of order $\mu_{[v]}=1$, and then 
$$
\delta_{[v]}(Q)={-\nu\over\Res_{[v]}(\nabla)}={1\over\iota_{[v]}
(Q)}\;.
$$

We can now study the complex foliation induced by~$Q$ in~$\C^2$ by means
of the horizontal foliation of~$N^{\otimes\nu}_S$. The link between the two
is provided by the following analogous of Theorem~\ctreb:

\newthm Proposition \ctre: Let $Q$ be a non-dicritical homogeneous vector field of degree~$\nu+1\ge 2$ in~$\C^2$, and let $S^o\subset\P^1(\C)$ be the complement
in $\P^1(\C)$ of the characteristic directions of~$Q$.
Then 
$$
\chi_\nu^*\omega={\nu\over w^2Q^1(w)-w^1 Q^2(w)}
\bigl[-Q^2(w)\,dw^1+Q^1(w)\,dw^2\bigr]
$$
on $S^o$ where $\omega$ is the holomorphic $1$-form representing the 
horizontal foliation on the total space of~$N^{\otimes\nu}_{S^o}$. In particular
a complex curve $\Lambda\subset\widehat S_Q$ is a (complex) leaf of
the holomorphic foliation induced by $Q$ if and only if $\chi_\nu(\Lambda)$ is a leaf of
the horizontal foliation of $N_{S^o}^{\otimes\nu}$.

\pf It is an easy computation.\qedn

So to study the holomorphic foliation induced by $Q$ is equivalent to
studying the horizontal foliation of the total space of $N_S^{\otimes\nu}$.
For instance, we can get a description of the closure of the leaves
completely analogous to the one given in Theorem~\tdueb. To state it we need
a couple of definitions.

\newdef Let $[v_1],\ldots,[v_g]\in\P^1(\C)$ be the characteristic 
directions of a non-dicritical homogeneous vector field $Q$ of 
degree $\nu+1\ge 2$ in~$\C^2$. Then the {\sl monodromy group} associated to~$Q$ is
the subgroup~$G(Q)$ of $\C^*$ given by
$$
G(Q)=\exp\left[2\pi i\left({1\over\nu}\Z\oplus\bigoplus_{h=1}^g\Z\iota_{[v_j]}(Q)
\right)\right]\subset\C^*\;.
$$ 
Notice that $G(Q)\subset S^1$ if and only if all indeces of~$Q$ are 
real numbers; in this case we shall say that $Q$ has {\sl real periods.}

\newrem It is not difficult to check that $G(Q)$ is a finite cyclic subgroup
of~$S^1$
if and only if there is $\ell\in\N^*$ such that
$$
\nu\,\iota_{[v_h]}(Q)\in{1\over\ell}\,\Z 
$$
for all $h=1,\ldots,g$. 

\newdef Let $Q=Q^1{\de\over\de z^1}+Q^2{\de\over\de z^2}$ be a non-dicritical homogeneous vector field of degree
$\nu+1\ge 2$. The {\sl metric foliation} of $\hat S_Q$ induced
by~$Q$ is the real rank-3 foliation given by the form 
$$
\Re\left({\nu\over w^2Q^1(w)-w^1 Q^2(w)}
\bigl[-Q^2(w)\,dw^1+Q^1(w)\,dw^2\bigr]\right)=\Re(\chi^*_\nu\omega)\;.
$$
In other words, it is the foliation induced by the metric foliation of
$N_S^{\otimes\nu}$ via $\chi_\nu$.

Then we have the following 

\newthm Theorem \clvf: Let $\Lambda\subset\C^2\setminus\{O\}$ be a 
non-characteristic leaf of the foliation induced by a non-dicritical homogeneous vector field 
$Q$ of degree $\nu+1\ge 2$ in $\C^2$. Then $[\Lambda]\subset\P^1(\C)$ is the complement~$S^o$ in~$\P^1(\C)$ of the characteristic directions of~$Q$.
Furthermore, take $[v]\in S^o$ and $z_0\in\Lambda\cap\C^* v$.
Then
$$
\Lambda\cap \C^*v=G(Q)\cdot z_0\qquad\hbox{and}\qquad
\bar\Lambda\cap\C^* v=\bar{G(Q)}\cdot z_0\;.
\neweq\eqclhf
$$
In particular, either
\smallskip
\itm{(i)} $Q$ has real periods, and then either all non-characteristic leaves of $Q$ are
closed in $\widehat S_Q$ (and this happens if and only if $G(Q)$ is a 
finite cyclic group) or any non-characteristic leaf is dense in the leaf of the metric foliation containing it (which is necessarily closed
in $\widehat S_Q$); or
\itm{(ii)} $Q$ does not have real periods, and then all non-characteristic leaves of $Q$
accumulate both the origin\break\indent and infinity in all directions.

\pf A non-characteristic leaf cannot intersect a characteristic one, and so $[\Lambda]\subseteq S^o$. By Proposition~\ctre,
$\chi_\nu(\Lambda)$ is then a $\nu$-to-1 cover of a leaf $L$ of the
horizontal foliation of $N_{S^o}^{\otimes\nu}$; in particular,
since $p(L)=S^o$ we get $[\Lambda]=S^o$.

Theorem~\tdueb\ says that
$$
\chi_\nu(\Lambda\cap\C^*v)=L\cap N^{\otimes\nu}_{[v]}=\rho(\pi)\cdot
\chi_\nu(z_0)\;,
$$
where $\rho(\pi)\subset\C^*$ is the image of the monodromy representation of the holomorphic
connection~$\nabla$ induced by $Q$ on~$N_{S^o}^{\otimes\nu}$. Therefore
$$
\Lambda\cap\C^*v=\rho(\pi)^{1/\nu}\cdot z_0\qquad\hbox{and}\qquad
\bar\Lambda\cap\C^*v=\bar{\rho(\pi)^{1/\nu}}\cdot z_0\;,
$$
where $\rho(\pi)^{1/\nu}=\{\zeta\in\C\mid\zeta^\nu\in\rho(\pi)\}$.
Hence to prove \eqclhf\ it suffices to show that $\rho(\pi)^{1/\nu}=G(Q)$.

Let $[v_1],\ldots,[v_g]\in\P^1(\C)$ be the characteristic directions of~$Q$.
By Proposition~\ttre, 
$$
\rho(\pi)=\exp\left[2\pi i\bigoplus_{h=1}^g\Z r_h\right]\;,
$$
where $r_h\in\C$ is the residue in~$[v_h]$ of the meromorphic connection~$\nabla^o$ induced by~$\nabla$ on~$TS^o$ via~$X_Q$. Therefore
$$
\rho(\pi)^{1/\nu}=\exp\left[2\pi i\left(\Z{1\over\nu}\oplus
\bigoplus_{h=1}^g\Z {r_h\over\nu}\right)\right]\;.
$$
Now, it is easy to check that
the meromorphic form $\eta^o_\alpha$ representing $\nabla^o$ in a canonical chart
is related to the
form $\eta_\alpha$ representing $\nabla$ in the corresponding chart by
$$
\eta^o_\alpha=\eta_\alpha-{1\over X_\alpha}\,dX_\alpha\;;
$$
therefore
$$
r_h=\Res_{[v_h]}(\nabla^o)=\Res_{[v_h]}(\nabla)-\ord_{[v_h]}(X_Q)\;.
\neweq\eqinducres
$$
Since $\ord_{[v_h]}(X_Q)\in\N$ we get $\rho(\pi)^{1/\nu}=G(Q)$, as claimed.

If $Q$ does not have real periods, then $G(Q)\subset\C^*$ accumulates 
both~0 and~$\infty$, and (ii) follows. If $Q$ has real periods, then
$G(Q)$ is a subgroup of~$S^1$, and hence it is either cyclic or dense,
and (i) follows from Theorem~\tdueb.(i).\qedn

To real flow of~$Q$ is instead described by the geodesic
foliation of~$N^{\otimes\nu}_{\P^1(\C)}$, that we studied in the first
part of this paper.
%
%
%
%
As a first consequence, we get the following 
Poincar\'e-Bendixson theorem for homogeneous vector fields:

\newthm Theorem \cqua: Let $Q$ be a homogeneous holomorphic vector
field on $\C^2$ of degree~$\nu+1\ge 2$, and let $\gamma\colon[0,\eps_0)\to\C^2$
be a recurrent maximal integral curve of $Q$. Then $\gamma$ is periodic or 
$[\gamma]\colon[0,\eps_0)\to\P^1(\C)$ intersects itself infinitely many times.

\pf If $Q$ is dicritical, Lemma~\dcl\ implies that 
no (non-constant) integral curve is recurrent.

Assume then that $Q$ is not dicritical. Then $\gamma$, being a recurrent integral curve, cannot intersect a characteristic leaf, because in that case it would be contained in it and, again
by Lemma~\dcl, no (non-constant) integral curve in a characteristic
leaf is recurrent. 

Since $\gamma$ is recurrent in $\C^2\setminus
\{O\}$, then $\chi_\nu\circ\gamma$ is recurrent in~$N^{\otimes\nu}_{\P^1(\C)}\setminus\P^1(\C)$. More precisely,
if we denote by $S$ the complement in~$\P^1(\C)$ of the characteristic
directions of~$Q$, we know that $\chi_\nu\circ\gamma$ is 
recurrent in $N^{\otimes\nu}_S\setminus S$ because the support of~$\gamma$
is contained in $\widehat S_Q$. Then $p\circ\chi_\nu
\circ\gamma$ is recurrent in~$S$; by Theorem~\qqua, this 
implies that $p\circ\chi_\nu\circ\gamma=[\gamma]$ is closed or intersect
itself infinitely many times. In the latter case we are done. In the former
case, $[\gamma]$ must satisfy the conditions described in Corollary~\qunomez; in
particular, the only way for $\chi_\nu\circ\gamma$ to be recurrent is
to be periodic. Since $\chi_\nu$ is a finite-to-one map,
it follows that $\gamma$ is periodic, as claimed.\qedn

\newthm Corollary \ccin: Let $\gamma\colon\R\to\C^2\setminus\{O\}$ be a
non-constant periodic integral curve of a
homogeneous vector field $Q$ of degree~$\nu+1\ge 2$. Then the
characteristic directions $[v_1],\ldots,[v_g]\in\P^1(\C)$ surrounded
by~$[\gamma]$ satisfy
$$
\sum_{j=1}^g \Res_{[v_j]}(Q)=-1+\sum_{j=1}^g \ord_{[v_j]}(X_Q)\;.
$$

\pf It follows immediately from Corollary~\qunomez\
and \eqinducres.\qedn

We can clearly say more along this line; but to get a better
understanding of the dynamics of $Q$ we need to know something about
the behavior of the geodesic field near the singularities. 
The next two sections are devoted to this task; 
we instead end this section describing a homogeneous vector field actually
having a periodic integral curve and, more generally,
examples of meromorphic connections on $\P^1(\C)$ having periodic geodesics,
closed geodesics, non-closed geodesics accumulating a closed one, and
geodesic converging to a pole, thus displaying the behavior
described in Theorem~\qqua.(i), (ii) and (iii).
\medbreak
{\sc Example 6.1:} Take
$$
Q=i\gamma (w^1)^2{\de\over\de w^1}+(1+i\gamma) w^1w^2{\de\over\de w^2}\;,
$$
where $\gamma\in\R$. In this case $\nu=1$; hence $\chi_\nu$ is a biholomorphism, and thus
$Q$ is biholomorphically conjugated to the geodesic field $G$.

The field $Q$ is non dicritical. It has two characteristic
directions: $0=[1:0]$, which is non-degenerate for $\gamma\ne 0$, and $\infty=[0:1]$, which is always
degenerate. In the chart centered at $0$ we have
$$
\displaylines{
X_0=\zeta_0\;,\qquad \eta_0=-{i\gamma\over\zeta_0}\,d\zeta_0\;,\qquad
G=\zeta_0 v_0{\de\over\de\zeta_0}+i\gamma(v_0)^2{\de\over\de v_0}\;;
\cr}
$$
in particular, $0$ is a Fuchsian singularity and, denoting by $\nabla^o$ the connection induced by $\nabla$
on~$\P^1(\C)$ via~$X_Q$, we have
$$
\mu_0=\ord_0(X_0)=1\;,\qquad\Res_0(\nabla)=-i\gamma\;,\qquad \Res_0(\nabla^o)
=-1-i\gamma\;.
$$

If $\gamma\ne 0$ the integral curves of $G$ in the standard chart
centered at~0 are of the form
$\sigma(t)=\bigl(\zeta(t),v(t)\bigr)$, with 
$$
\zeta(t)=\zeta_0\exp\left[{i\over\gamma}\log(1-i\gamma v_0 t)
\right]\;,\qquad v(t)={v_0\over 1-i\gamma v_0 t}\;,
$$
where we have chosen the determination of the logarithm so that
$\zeta(0)=\zeta_0$. The 
curve $\zeta$ is a geodesic for the meromorphic connection $\nabla^o$,
and we can write 
$$
\zeta(t)=\zeta_0\exp\bigl[i\gamma^{-1}\log|1-i\gamma v_0 t|-
\gamma^{-1}\arg(1-i\gamma v_0 t)\bigr]\;.
$$
If $i\gamma v_0\in\R^*$ then 
$\arg(1-i\gamma v_0 t)$ is equal
either to 0 or to~$\pi$
depending on the sign of $1-i\gamma v_0 t$. In both cases 
$\arg(1-i\gamma v_0 t)$ is constant, and so $\zeta$ is a closed
geodesic. Notice that if $i\gamma v_0\in\R^+$ then $\zeta$ is defined
only on the half-line $(-\infty,(i\gamma v_0)^{-1})$. 

If instead $i\gamma v_0\notin\R$ then $\arg(1-i\gamma v_0t)\to\arg
(-i\gamma v_0)$ as $t\to+\infty$, and thus the geodesic $\zeta$ 
accumulates a circumference, which is easily seen to be the support of a
closed geodesic.

In all these cases $v(t)\to O$ as $t\to+\infty$ (except when $i\gamma v_0\in
\R^+$; in this case $|v(t)|\to+\infty$ as $t\to(i\gamma v_0)^{-1}$), which
means that the integral curves of~$Q$ are converging to the origin
(or escaping to infinity in the exceptional case) without being tangent
to any direction (because $\zeta(t)$ has no limit); in particular,
they are not periodic.

However, if $\gamma=0$ we have $\Res_0(\nabla^o)=-1$,
and so we expect to find periodic integral curves of $Q$. Indeed,
if $\gamma=0$ the integral curves of $G$ in the standard chart 
centered at~0 are given by 
$$
\zeta(t)=\zeta_0\exp(v_0 t)\;,\qquad v(t)\equiv v_0\;.
$$
Therefore if $v_0\in i\R$ we get
a periodic integral curve, as desired. If instead $\Re v_0<0$ 
we have $\zeta(t)\to O$, which means that the geodesic tends to~$[1:0]$ and that the integral curve of~$Q$
issuing from $(v_0,\zeta_0 v_0)$ tends to a non-zero element of the 
characteristic leaf~$L_{(1,0)}\subset\C^2$. Finally, if $\Re v_0>0$ we have
$|\zeta(t)|\to+\infty$, which means that the geodesic tends to $[0:1]$ and
the integral curve of~$Q$
issuing from $(v_0,\zeta_0 v_0)$ escapes to infinity.

\smallsect 7. Local study of the singularities

The next step consists in the study of the geodesic flow nearby singular points.
We shall work in the following setting: we have a line bundle $E$ on
a Riemann surface $S$ (e.g., $S=\P^1(\C)$ and $E=N^{\otimes\nu}_S$);
a morphism $X\colon E\to TS$ which is an isomorphism on $S^o=X\setminus\Sing(X)$,
where $\Sing(X)$ contains only isolated points (e.g., $X=X_Q$); and a meromorphic connection
$\nabla$ on $E$, holomorphic on $E|_{S^o}$, such that the geodesic field $G$
extends holomorphically from $E|_{S^o}$ to the whole of~$E$.

Fix a local chart $(U_\alpha,z_\alpha,e_\alpha)$ trivializing~$E$ and centered at a 
singular point~$p_0$ of~$X$. 
In local coordinates we can write
$$
X(e_\alpha)=X_\alpha{\de\over\de z_\alpha}
$$
where $X_\alpha\colon U_\alpha\to\C$ is holomorphic with $X(p_0)=0$. 
The holomorphic 1-form $\eta_\alpha$ representing $\nabla$ in these coordinates is of the form $\eta_\alpha=k_\alpha\,dz_\alpha$, where
$k_\alpha$ is meromorphic with possibly a pole in $p_0$. The geodesic field $G$ in these
coordinates is given by
$$
G=X_\alpha v_\alpha\de_\alpha-(X_\alpha k_\alpha) v_\alpha^2{\de\over\de v_\alpha}\;;
$$
saying that $G$ extends holomorphically over the singular point means that 
$Y_\alpha=X_\alpha k_\alpha$ is holomorphic in~$U_\alpha$. Therefore we shall
write
$$
\eta_\alpha={Y_\alpha\over X_\alpha}\,dz_\alpha\quad\hbox{and}\quad
G=X_\alpha v_\alpha\de_\alpha-Y_\alpha v_\alpha^2{\de\over\de v_\alpha}\;,
\neweq\eqszero
$$
with $X_\alpha$, $Y_\alpha\in\ca O(U_\alpha)$ and $X_\alpha(p_0)=0$. 

We shall denote by $\nabla^o$ the meromorphic connection on $TS$ induced by $\nabla$
via~$X$. We have already remarked that the meromorphic form $\eta^o_\alpha$ representing 
$\nabla^o$ in the chart $(U_\alpha,z_\alpha,\de/\de z_\alpha)$ is given by
$$
\eta^o_\alpha=\eta_\alpha-{1\over X_\alpha}\,dX_\alpha\;.
$$
In particular,
$$
\Res_{p_0}(\nabla^o)=\Res_{p_0}(\nabla)-\ord_{p_0}(X_\alpha)\;.
\neweq\eqsuno
$$
We shall write
$$
X_\alpha=z_\alpha^{\mu_X}h^X_\alpha\qquad\hbox{and}\qquad
Y_\alpha=z_\alpha^{\mu_{\alpha,Y}}h^Y_\alpha\;,
$$
with $\mu_X=\ord_{p_0}(X_\alpha)$, $\mu_{\alpha,Y}=\ord_{p_0}(Y_\alpha)$ and 
$h_\alpha^X(p_0)$, $h^Y_\alpha(p_0)\ne 0$. Notice that $\mu_X$ does not depend on the
coordinates, by \eqsuno, whereas $\mu_{\alpha,Y}$ in general does. Let us recast
Definition~6.1 in this context:

\newdef The {\sl order}~$\mu_X\ge 1$ of $p_0\in\Sing(X)$ is the
order of vanishing of~$X$ at~$p_0$, that is $\mu_X=\ord_{p_0}(X)$. Choosing
a local chart $(U_\alpha,z_\alpha,e_\alpha)$ and
writing $G$ as in \eqszero, we say that $p_0$ is
\smallskip
\item{--} an {\sl apparent singularity} if $\mu_X\le\mu_{\alpha,Y}$;
\item{--} a {\sl Fuchsian singularity} if $\mu_X=\mu_{\alpha,Y}+1$;
\item{--} an {\sl irregular singularity} if $\mu_X>\mu_{\alpha,Y}+1$; in this case
$m=\mu_X-\mu_{\alpha,Y}>1$ is the {\sl irregularity} of~$p_0$;
\item{--} a {\sl degenerate singularity} if $\mu_{\alpha, Y}\ge 1$.
\smallskip
\noindent In Remark~7.3 we shall see that these notions do not depend on the chosen chart.
Finally, the {\sl residue} of~$p_0$ is $\Res_{p_0}(\nabla)$, and the {\sl induced residue}
of~$p_0$ is $\Res^o_{p_0}(X)=\Res_{p_0}(\nabla)-\mu_X$. 

\newrem With these notations, we have
$$
\eta^o_\alpha=\left({h_\alpha^Y\over z_\alpha^{\mu_X-\mu_{\alpha,Y}}h^X_\alpha}
-{\mu_X\over z_\alpha}-{(h^X_\alpha)'\over h_\alpha^X}\right)\,dz_\alpha\;.
$$
In particular, an apparent singularity is a pole of~$\nabla^o$ without being
a pole of~$\nabla$. Furthermore, a Fuchsian singularity of $X$ is a Fuchsian singularity
of~$\nabla^o$ in the classical sense (that is a pole of order~1) unless the induced residue
of~$p_0$ vanishes. This is the only possibility for $p_0$ being a pole of~$\nabla$
but not of~$\nabla^o$; indeed, 
is easy to see that this happens if and only if
$$
\ord_{p_0}(z_\alpha Y_\alpha-\mu_X X_\alpha)>\mu_X=\ord_{p_0}(X_\alpha)\;,
$$
which is equivalent to
$$
\ord_{p_0}(Y_\alpha)=\ord_{p_0}(X_\alpha)-1\quad\hbox{and}\quad \Res_{p_0}(\nabla)=\mu_X\ge 1\;,
$$
that is to $p_0$ being Fuchsian with vanishing induced residue. 

\newrem If $p_0$ is a degenerate singularity, that is $Y_\alpha(p_0)=0$,
then the geodesic field $G$ vanishes identically when restricted to the fiber
$E_{p_0}$ over the singularity.\hfill\break\indent
When instead $\ord_{p_0}(Y_\alpha)=0$, that is $a_0=Y_\alpha(p_0)\ne 0$, the geodesic field restricted to~$E_{p_0}$ is given by $-a_0v_\alpha^2{\de\over\de v_\alpha}$. In particular, $E_{p_0}$ is $G$-invariant, and the integral curve
of~$G$ issuing from $v_0 e_\alpha(p_0)\in E_{p_0}\setminus\{O\}$ is given by
$$
v(t)={v_0\over 1+a_0v_0 t}\,e_\alpha(p_0)\;.
$$
Notice that $\lim\limits_{t\to\pm\infty} v(t)=O$, and that $v(t)$ is defined for all~$t\in\R$
unless $av_0\in\R^-$.

Our aim is to simplify as much as possible the expression of $G$ by changing the local chart
and the local generator of~$E$. Thus we shall use changes of the form
$$
(z_\beta,v_\beta)=\phe(z_\alpha,v_\alpha)=\bigl(\psi(z_\alpha),\xi(z_\alpha)v_\alpha
\bigr)\;,
\neweq\eqallcc
$$
with $\psi(0)=0$ and $\psi'(0)$, $\xi(0)\ne 0$ (in the notations of Section~1, we have
$\psi'=\psi_{\beta\alpha}$ and $\xi=\xi_{\beta\alpha}$). A quick computation gives
$$
X_\beta\circ\psi={\psi' X_\alpha\over\xi}\qquad\hbox{and}\qquad
Y_\beta\circ\psi={1\over\xi}Y_\alpha-{\xi'\over\xi^2}X_\alpha\;.
\neweq\eqsdue
$$

\newrem In particular, if $\mu_X>\mu_{\alpha,Y}$ then $\mu_{\beta,Y}=\mu_{\alpha,Y}$; if $\mu_X\le\mu_{\alpha, Y}$ then $\mu_X\le\mu_{\beta,Y}$; and
$\mu_{\alpha,Y}\ge 1$ implies $\mu_{\beta,Y}\ge 1$.
So Definition~7.1 is well-posed.

\newrem Since we only allow changes of coordinates of the form \eqallcc, the normal forms we are going to obtain for $G$ are (related to but) different from the usual normal forms for
singular holomorphic vector fields in $\C^2$ with respect to unrestricted (formal or) biholomorphic changes of coordinates. Furthermore, as it will
become apparent with Theorem~8.1, our classification is also different
from the classical classification of meromorphic connections on~$\P^1(\C)$.

Let us choose $\psi=\id$ and $\xi=h^X_\alpha$. We find
$$
X_\beta=z_\beta^{\mu_X}\qquad\hbox{and}\qquad 
Y_\beta={1\over h_\alpha^X}\bigl(z_\beta^{\mu_{\alpha,Y}}h^Y_\alpha-z_\beta^{\mu_X}
(h^X_\alpha)'\bigr)=z_\beta^{\mu_{\beta,Y}}h^Y_\beta\;,
$$
where
\smallskip
\item{--} if $\mu_{\alpha,Y}<\mu_X$ then $\mu_{\beta,Y}=\mu_{\alpha,Y}$ and 
$h_\beta^Y=\bigl(h^Y_\alpha-z_\beta^{\mu_X-\mu_{\alpha,Y}}(h_\alpha^X)'\bigr)/
h_\alpha^X$;
\item{--} if $\mu_{\alpha,Y}\ge\mu_X$ then $\mu_{\beta,Y}=\mu_X+\ord_{p_0}
\bigl((h_\alpha^X)'-z_\beta^{\mu_{\alpha,Y}-\mu_X}h_\alpha^Y\bigr)\ge\mu_X$.
\smallskip
\noindent So we can assume $h^X_\alpha\equiv 1$. 

The first result of this section is a complete holomorphic classification of apparent singularities:

\newthm Proposition \suno: Let $E$ be a line bundle on
a Riemann surface $S$, and assume we have  
a morphism $X\colon E\to TS$ which is an isomorphism on $S^o=X\setminus\Sing(X)$,
and a meromorphic connection
$\nabla$ on $E$, holomorphic on $E|_{S^o}$, such that the geodesic field $G$
extends holomorphically from $E|_{S^o}$ to the whole of~$E$. Let $p_0\in\Sing(X)$
be an apparent singularity of order~$\mu$. Then there exists a chart 
$(U,z,e)$ centered at~$p_0$ such that $G$ in this chart is given by
$$
G=\cases{\displaystyle
z v {\de\over\de z}&if $\mu=1$,\cr
\displaystyle z^\mu(1+a z^{\mu-1})v{\de\over\de z}&for some $a\in\C$ if $\mu>1$.\cr}
$$
Furthermore, if $\mu>1$ then $a\in\C$ is a (holomorphic and formal) invariant.

\pf We have already remarked that we can find a chart $(U_\alpha,z_\alpha,e_\alpha)$ such
that $h^X_\alpha\equiv 1$. 
Furthermore, $\mu_{\alpha,Y}\ge\mu_X$ because $p_0$ is an apparent singularity.

As a first step, take $\psi=\id$ and $\xi$ solving the Cauchy problem
$$
\cases{\displaystyle {1\over\xi}z_\alpha^{\mu_{\alpha,Y}}h^Y_\alpha-{\xi'\over\xi^2}
z_\alpha^\mu=0\;,\cr
\noalign{\smallskip}
\xi(p_0)=1\;;
\cr}
$$
this problem has a solution holomorphic in a neighbourhood of~$p_0$ because we
can rewrite the differential equation in the form
$$
\xi'=z^{\mu_{\alpha,Y}-\mu}h_\alpha^Y\xi\;.
$$
Then \eqsdue\ says that in the new coordinates $(U_\beta,z_\beta,e_\beta)$ we have
$$
X_\beta=z_\beta^\mu h^X_\beta\qquad\hbox{and}\qquad Y_\beta\equiv 0\;,
$$
with $h_\beta^X(p_0)=1$. 

Now, it is known (see, e.g., [IY, Theorem~5.25]) that we can find a local 
change of variable $\psi$ fixing~$p_0$ such that 
$$
(\psi' X_\beta)\bigl(\psi^{-1}(z)\bigr)=\cases{z &if $\mu=1$,\cr
z^\mu(1+a z^{\mu-1})&for some $a\in\C$ if $\mu>1$;\cr}
$$
hence taking $\xi\equiv 1$ we bring $G$ in the
required form. Finally, the last assertion follows from the corresponding fact for
one variable vector fields (see, e.g., [IY, Section~$6{\bf B}_2$]).\qedn

\newdef The invariant $a$ introduced in the latter proposition is the {\sl apparent
index} of the apparent singularity.

As a consequence, we can describe the behavior of the geodesic flow nearby
an apparent singularity:

\newthm Corollary \sdue: Let $E$ be a line bundle on
a Riemann surface $S$, and assume we have  
a morphism $X\colon E\to TS$ which is an isomorphism on $S^o=X\setminus\Sing(X)$, and a meromorphic connection
$\nabla$ on $E$, holomorphic on $E|_{S^o}$, such that the geodesic field $G$
extends holomorphically from $E|_{S^o}$ to the whole of~$E$. Let $p_0\in\Sing(X)$
be an apparent singularity of order~$\mu$ and apparent index $a\in\C$ if $\mu>1$.
Then if $\sigma\colon[0,\eps)\to S^o$ is a geodesic
with $\sigma(t)\to p_0$ as $t\to\eps\in(0,+\infty]$, then $\sigma'(t)\to O_{p_0}$ and 
$X^{-1}\bigl(\sigma'(t)\bigr)$ tends to a non-zero element of~$E_{p_0}$.\hfill\break
\indent
Furthermore, there is a neighbourhood $U\subseteq S$ of $p_0$ such that if $z_0\in U
\setminus\Sing(X)$, and $\sigma_v\colon[0,\eps_v)\to S^o$ denotes the maximal geodesic
issuing from~$z_0$ in the direction~$X(v)\in T_{z_0}S$, then
\itm{(i)} if $\mu=1$ then there is a non-zero direction $v_0\in E_{z_0}$ such that
\itemitem{--} if $v=\zeta v_0\in E_{z_0}$ with $\Re\zeta<0$ then $\sigma_v(t)\to p_0$;
\itemitem{--} if $v=\zeta v_0\in E_{z_0}$ with $\Re\zeta>0$ then $\sigma_v$
escapes from $U$;
\itemitem{--} if $v=\zeta v_0\in E_{z_0}$ with $\Re\zeta=0$ then $\sigma_v$ is a
periodic geodesic surrounding~$p_0$;
\itm{(ii)} if $\mu>1$ and $a=0$ then there are $\mu-1$ non-zero 
directions $v_1,\ldots, v_{\mu-1}\in E_{z_0}$ such that
\itemitem{--} if $v\notin\R v_1\cup\cdots\cup\R v_{\mu-1}$ then $\sigma_v(t)\to p_0$;
\itemitem{--} if $v\in\R v_1\cup\cdots\cup\R v_{\mu-1}$ then $\sigma_v(t)$ escapes
from~$U$.
\itm{(iii)} if $\mu>1$ and $a\ne 0$ then there is a non-zero direction
$v_0\in E_{z_0}$ such that
\itemitem{--} if $v=\zeta v_0\in E_{z_0}$ with $\Re(\zeta/a)> 0$ then 
$\sigma_v(t)\to p_0$;
\itemitem{--} if $v=\zeta v_0\in E_{z_0}$ with $\Re(\zeta/a)=0$ then 
either $\sigma_v(t)\to p_0$, or $\sigma_v$ is a periodic geodesic surrounding~$p_0$, or $\sigma_v$
escapes from~$U$;
\itemitem{--} if $v=\zeta v_0\in E_{z_0}$ with $\Re(\zeta/a)<0$ then either
$\sigma_v(t)\to p_0$ or $\sigma_v$ escapes from~$U$.

\pf By Proposition~\suno, we can find a chart $(U,z,e)$ centered at~$p_0$ such that
a curve $\sigma\colon[0,\eps)\to U$ is a geodesic if and only if 
$X^{-1}\bigl(\sigma'(t)\bigr)=\bigl(z(t),v(t)\bigr)$ satisfies 
$$
\cases{z'=z^\mu(1+a z^{\mu-1}) v\;,\cr
v'=0\;,
\cr}
\neweq\eqstre
$$
with $a=0$ if $\mu=1$.
In particular, $v(t)\equiv v(0)$, and thus clearly has a finite non-zero limit
as $t\to\eps$. 
Moreover, $\sigma(t)\to p_0$ if and only if $z(t)\to 0$; hence $z(t)^\mu\bigl(1+az(t)^{\mu-1}\bigr) v(t)\to0$ as $t\to\eps$. This means exactly
that $\sigma'(t)\to O_{p_0}$, and the first assertion is proved.

Assume now $\mu=1$. Then solving \eqstre\ we find
$$
z(t)=z_0\exp(v t)\;,
$$
where $z_0=z(0)$ and $v(t)\equiv v$. In particular, $\sigma(t)\to p_0$ as $t\to\eps$ if and
only if $\Re v<0$. If $\Re v=0$ we get a 
periodic geodesic around $p_0$ (and indeed the induced residue of~$p_0$ is~$-1$).
Finally, if $\Re v>0$ then $\sigma(t)$ escapes from $U$. So we have proved (i),
with $v_0=e(z_0)$.  

If instead $\mu>1$ and $a=0$ solving \eqstre\ we get
$$
z(t)=z_0\left(1-{v z_0^{\mu-1}\over\mu-1}t\right)^{-1/(\mu-1)}\;,
$$
where the determination of the root is chosen so that $z(0)=z_0$, and $v(t)\equiv v$.
In particular $\sigma(t)\to p_0$ as $t\to+\infty$ unless $v z_0^{\mu-1}\in\R^+$,
and we get (ii).

If $\mu>1$ and $a\ne 0$ the solution of \eqstre\ satisfies $v(t)\equiv v$
and
$$
\left(1+{1\over az(t)^{\mu-1}}\right)\exp\left[-\left(1+{1\over az(t)^{\mu-1}}\right)\right]
=c_0\exp\left[{(\mu-1)v\over a}\,t\right]\;,
\neweq\eqscinbb
$$
with
$$
c_0=\left(1+{1\over az_0^{\mu-1}}\right)\exp\left(-\left(1+{1\over az_0^{\mu-1}}\right)
\right)\;.
$$
In particular, if $z_0$ is one of the $\mu-1$ roots of
$az_0^{\mu-1}=-1$, then $c_0=0$ and so $z(t)\equiv z_0$. We assume that $U$
is small enough to exclude these points; in particular, we have $c_0\ne 0$.

Assume $\Re(v/a)>0$. Then the modulus of the right-hand side of \eqscinbb\
goes to~$+\infty$ as $t\to+\infty$. Hence the modulus of the left-hand side of
\eqscinbb, given by $|w(t)|\exp\bigl(-\Re w(t)\bigr)$ where
$$
w(t)=1+{1\over a z(t)^{\mu-1}}\;,
$$
goes to $+\infty$ too. This forces $|w(t)|\to+\infty$, and hence $z(t)\to 0$, as required.

If $\Re(v/a)=0$, then the left-hand side of \eqscinbb\ has 
constant modulus~$|c_0|$, and argument going to $\pm\infty$. Looking at the level sets of the function
$|w e^{-w}|$ one sees that there is a critical value $\hat c$ so that
\smallskip
\item{--} if $|c_0|\ge\hat c$ then necessarily $\arg w(t)$ is bounded and $|\Im w(t)|\to+\infty$;
\item{--} if $|c_0|<\hat c$ then, depending on which connected component of
the level set contains~$w(0)$, either $\arg w(t)$ is bounded and $|\Im w(t)|\to+\infty$, or $\arg w(t)$ is unbounded and $|w(t)|$ is bounded
and bounded away from zero.
\smallskip\noindent
If $|\Im w(t)|\to+\infty$, then $|w(t)|\to+\infty$, 
and hence $z(t)\to 0$ as before. If instead $|w(t)|$ is bounded, then
$w(t)$ must periodically trace the bounded connected component of the
level set of $|w e^{-w}|$, and hence $z(t)$ is a periodic geodesic
(or escapes $U$ if $U$ is too small). 

Finally, if $\Re(v/a)<0$ then the modulus of the right-hand side of \eqscinbb\ goes to~$0$ as $t\to+\infty$, and so either $w(t)\to 0$
or $\Re w(t)\to+\infty$. In the former case, $z(t)$ must tend to one
of the excluded points, and so $\sigma_v$ escapes. In the latter case
$|w(t)|\to+\infty$, and hence $z(t)\to 0$ as usual.\qedn

The translation of this corollary for homogeneous vector fields is the
following:

\newthm Corollary \stre: Let $Q$ be a homogeneous holomorphic vector
field on $\C^2$ of degree~$\nu+1\ge 2$. 
Let $[v_0]\in\P^1(\C)$ be an apparent singularity 
of~$X_Q$ of order~$\mu\ge 1$ (and apparent index $a\in\C$ if $\mu>1$). Then:
\smallskip
\itm{(i)} if the direction $[\gamma(t)]\in\P^1(\C)$ of an integral curve
$\gamma\colon[0,\eps)\to\C^2\setminus\{O\}$ of $Q$ tends to $[v_0]$
as $t\to\eps$ then $\gamma(t)$ tends to a non-zero point
of the characteristic leaf $L_{v_0}\subset\C^2$;
\itm{(ii)} no integral curve of $Q$ tends to the origin tangent to~$[v_0]$;
\itm{(iii)} there is an
open set of initial conditions whose integral curves tend to a non-zero 
point of~$L_{v_0}$;
\itm{(iv)} if $\mu=1$ or $\mu>1$ and $a\ne 0$ then $Q$ admits periodic orbits of
arbitrarily long periods accumulating at\break\indent the origin.

\pf Notice that $Q$ is identically zero on~$L_{v_0}$; therefore either
an integral curve is a constant point of~$L_{v_0}$ (and then all the
assertions are trivial) or does not intersect $L_{v_0}$. Furthermore,
since characteristic leaves are $Q$-invariant, we are interested only
in integral curves $\gamma$ contained in $\hat S_Q$.

Assume that $[\gamma(t)]$ converges to~$[v_0]$. Set $\sigma=[\gamma]=p\circ\chi_\nu\circ\gamma$; then $\sigma$ is
a geodesic converging to~$[v_0]$, and $X^{-1}(\sigma')=\chi_\nu\circ\gamma$.
Then Corollary~\sdue\ says that $\chi_\nu\circ\gamma$ tends to a
non-zero element of~$E_{v_0}$, where $E=N^{\otimes\nu}_{\P^1(\C)}$. Now,
$E_{v_0}=\chi_\nu(L_{v_0})$, and $\chi_\nu$ is a $\nu$-to-1 map; since
the set of accumulation points of~$\gamma$ is connected and contained 
in~$L_{v_0}$, it follows that $\gamma(t)$ tends to a non-zero element
di~$L_{v_0}$, and (i) is proved. (ii) follows from (i), and (iii)
follows from Corollary~\sdue.(i), (ii) and (iii).

Finally, the periodc geodesics given by Corollary~\sdue.(i) and (iii) yields 
periodic integral curves accumulating the origin of arbitrarily long period
(because the period is inversally proportional to the modulus of the vector $v$ giving
the initial condition of the geodesic), and we get (iv).\qedn

The next section is devoted to the classification (formal and, when
possible, holomorphic) of Fuchsian and irregular singularities,
and its dynamical consequences. We end this section with a
preliminary result on the dynamics of Fuchsian singularities that
shall be useful later to deal with resonances.

\newthm Proposition \squa: Let $E$ be a line bundle on
a Riemann surface $S$, and assume we have  
a morphism $X\colon E\to TS$ which is an isomorphism on $S^o=X\setminus\Sing(X)$, and a meromorphic connection
$\nabla$ on $E$, holomorphic on $E|_{S^o}$, such that the geodesic field $G$
extends holomorphically from $E|_{S^o}$ to the whole of~$E$. Let $p_0\in\Sing(X)$ be a Fuchsian
singularity, and assume that $\Res_{p_0}(\nabla)
\in\R^*$. Then:
\smallskip
\itm{(i)} if $\Res_{p_0}(\nabla)<0$ then all leaves of the metric
foliation over $p\in S^o$ tend to the zero section as~$p\to p_0$;
\itm{(ii)} if $\Res_{p_0}(\nabla)>0$ then all leaves of the metric
foliation over $p\in S^o$ tend to infinity as~$p\to p_0$.
\smallskip
\noindent In particular, if $\sigma\colon[0,\eps)\to S^o$ is a geodesic
with $\sigma(t)\to p_0$ as~$t\to\eps$, then 
\smallskip
\itm{(a)} if $\Res_{p_0}(\nabla)<0$ then $X^{-1}\bigl(\sigma'(t)\bigr)$
tends to $O_{p_0}$ as $t\to\eps$;
\itm{(b)} if $\Res_{p_0}(\nabla)>0$ then $|X^{-1}\bigl(\sigma'(t)\bigr)|\to
+\infty$ as $t\to\eps$.

\pf Saying that $p_0$ is Fuchsian means that $p_0$ is a pole of
order~1 of~$\nabla$. Therefore in local coordinates $(U_\alpha,z_\alpha,
e_\alpha)$ centered at~$p_0$
we can write $\eta_\alpha=k_\alpha\,dz_\alpha$ with
$$
k_\alpha={\rho\over z_\alpha}+k^*_\alpha\;,
$$
where $\rho=\Res_{p_0}(\nabla)$ and $k^*_\alpha$ is holomorphic in a neighbourhood
of~$p_0$. 
 
In Remark~1.1 we noticed that, since the residue $\rho$ is real, 
we can find a metric $g_\alpha$ adapted to~$\nabla$ 
in~$U_\alpha\setminus\{p_0\}$ by setting
$$
g_\alpha(z_\alpha;v_\alpha)=\exp\bigl(2\Re K^*_\alpha(z_\alpha)\bigr)
|z_\alpha|^{2\rho} |v_\alpha|^2\;,
\neweq\eqsqua
$$
where $K^*_\alpha\in\ca O(U_\alpha)$ is the holomorphic primitive 
of~$k^*_\alpha$ with $K^*_\alpha(p_0)=0$. The leaves of the metric
foliation over $U_\alpha\setminus\{p_0\}$ are the level sets of~$g_\alpha$,
and then (i) and (ii) clearly follow. Assertions (a) and (b) are then
consequences of the fact that the geodesic foliation is contained in
the metric foliation.\qedn

\newthm Corollary \scin: Let $Q$ be a homogeneous holomorphic vector
field on $\C^2$ of degree~$\nu+1\ge 2$. 
Let $[v_0]\in\P^1(\C)$ be a Fuchsian characteristic direction 
of~$Q$ with real residue~$\rho\in\R^*$. Assume that the direction 
$[\gamma(t)]\in\P^1(\C)$ of an integral curve
$\gamma\colon[0,\eps)\to\C^2\setminus\{O\}$ of $Q$ tends to $[v_0]$
as $t\to\eps\in(0,+\infty]$. Then $\gamma(t)$ tends to the origin if
$\rho<0$, and to infinity if $\rho>0$. In particular, if $\rho>0$ no
integral curve outside the characteristic leaf $L_{v_0}$ can tend to the origin tangent to~$[v_0]$.

\pf It follows from Proposition~\squa, as usual.\qedn

\newrem The residue in a Fuchsian singularity~$p_0$ is necessarily
different from zero (otherwise $p_0$ would not be a pole of the connection).

\newrem In Corollary~8.5 we shall see that the same conclusions can be
inferred considering the real part of the residue when the residue is
not real; but the proof in the resonant case shall depend on Corollary~\scin.

\newrem In the irregular case, \eqsqua\ still holds, but now 
$K^*_\alpha$ is meromorphic in~$U_\alpha$, with a pole at~$p_0$. 
Therefore the behavior of the leaves of the metric foliation might
depend on the way $p$ approaches~$p_0$.

\smallsect 8. Classification of singularities

The next result provides the formal classification of Fuchsian and 
irregular singularities. Contrarily to the classical case of 
meromorphic connections on~$\P^1(\C)$, in the Fuchsian case we might have resonances.

\newthm Theorem \stuno: Let $E$ be a line bundle on
a Riemann surface $S$, and assume we have  
a morphism $X\colon E\to TS$ which is an isomorphism on $S^o=X\setminus\Sing(X)$,
and a meromorphic connection
$\nabla$ on $E$, holomorphic on~$E|_{S^o}$, such that the geodesic field $G$
extends holomorphically from $E|_{S^o}$ to the whole of~$E$. 
Let $p_0\in\Sing(X)$ be a Fuchsian or irregular singularity, and in a chart
$(U_\alpha,z_\alpha,e_\alpha)$ centered at~$z_0$ write
$$
G=z_\alpha^{\mu_X}(a_0+a_1 z_\alpha+\cdots)v_\alpha\de_\alpha-
z_\alpha^{\mu_Y}(b_0+b_1 z_\alpha+\cdots)v_\alpha^2{\de\over\de v_\alpha}\;,
$$
with $\mu_X>\mu_Y$ and $a_0$, $b_0\ne 0$. Put $\rho=b_0/a_0\ne 0$; if $p_0$ is
Fuchsian then $\rho=\Res_{p_0}(\nabla)$. Then:
\smallskip
\itm{(i)} if $p_0$ is Fuchsian then 
\itemitem{\rm(a)} if $\mu_Y-\rho\notin\N^*$ then $G$ is
formally conjugated to
$$
z^{\mu_X-1}\left(z v\de-\rho v^2{\de\over\de v}\right)\;;
$$
\itemitem{\rm (b)} if $\mu_Y-\rho=n\in\N^*$ then $G$ is
formally conjugated to
$$
z^{\mu_X-1}\left(zv\de-\rho v^2(1+a z^n){\de\over\de v}\right)
$$
for a suitable $a\in\C$ which is a formal invariant;
\itm{(ii)} if $p_0$ is irregular then
$G$ is formally conjugated to
$$
z^{\mu_X-m}\left(z^m v\de-\rho v^2(1+a z^{m-1}){\de\over\de v}\right)\;,
$$
\indent where $m=\mu_X-\mu_Y>1$ is the irregularity and $a=\Res_{p_0}(\nabla)/\rho$
is a formal invariant.

\pf The proof follows the usual formal Poincar\'e-Dulac paradigm. Write
$$
X_\alpha=z_\alpha^{\mu_X}\sum_{j=0}^{+\infty}a_j z_\alpha^j\quad\hbox{and}
\quad Y_\alpha=z_\alpha^{\mu_Y}\sum_{j=0}^{+\infty}b_j z_\alpha^j\;;
$$
given $n\in\N^*$, we start computing the action on~$X_\alpha$
and $Y_\alpha$ of a change of coordinates of the form 
$$
(z_\beta,v_\beta)=\phe(z_\alpha,v_\alpha)=(z_\alpha+c_1 z_\alpha^{n+1},
v_\alpha+c_2 z_\alpha^n v_\alpha)\;,
$$
with $(c_1,c_2)\in\C^2\setminus\{O\}$. Using \eqsdue\ it is easy to see that
$$
X_\beta=z_\beta^{\mu_X}\left[\sum_{j=0}^{n-1} a_j z_\beta^j+
\bigl[a_n+a_0\bigl((n+1-\mu_X)c_1-c_2\bigr)\bigr]z_\beta^n+o(z_\beta^n)\right]
$$
and
$$
Y_\beta=\cases{\displaystyle
z_\beta^{\mu_Y}\left[\sum_{j=0}^{n-1} b_j z_\beta^j+
\bigl[b_n-(\mu_Yc_1+c_2)b_0-nc_2a_0\bigr]z_\beta^n+o(z_\beta^n)\right]&if
$m=\mu_X-\mu_Y=1$,\cr
\displaystyle
z_\beta^{\mu_Y}\left[\sum_{j=0}^{n-1} b_j z_\beta^j+
\bigl[b_n-(\mu_Yc_1+c_2)b_0\bigr]z_\beta^n+o(z_\beta^n)\right]&if
$m=\mu_X-\mu_Y>1$.\cr}
$$
So such a change of coordinates does not modify the terms of degree
less than~$n$, and acts in the specified way on the terms of degree $n$. In
particular, to get $X_\beta$ and $Y_\beta$ without terms of degree~$n$
we must choose $c_1$ and $c_2$ so that
$$
\cases{a_0(\mu_X-n-1)c_1+a_0c_2=a_n\;,\cr
\mu_Y b_0 c_1+(na_0+b_0)c_2=b_n\;,
\cr}
\neweq\eqstuno
$$
if $m=1$, or so that
$$
\cases{a_0(\mu_X-n-1)c_1+a_0c_2=a_n\;,\cr
\mu_Y b_0 c_1+b_0c_2=b_n\;,
\cr}
\neweq\eqstdue
$$
if $m>1$.

In case (i), the determinant of the system \eqstuno\ vanishes if and only if
$$
n=\mu_Y-\rho\;.
$$
If this happens, the first equation in \eqstuno\ becomes $a_0(\rho c_1+c_2)
=a_n$, and so it can anyway be solved.

In case (ii), the determinant of the system \eqstdue\ vanishes if and 
only if
$$
n+1=m\;.
$$
If this happens, the first equation in \eqstdue\ can anyway be solved.

Summing up, in case (i) we can always kill the term of degree $n$ in 
$X_\beta$, whereas we can kill the term of degree $n$ in $Y_\beta$
if $n\ne \mu_Y-\rho$. In case (ii), we can always kill the term of 
degree $n$ in $X_\beta$, whereas we can kill the term of degree $n$ in 
$Y_\beta$ if $n\ne m-1$. Furthermore, in both cases a quick inspection of
\eqstuno\ and \eqstdue\ shows that (as soon as we have killed the terms below
the resonance level) we cannot anymore modify the coefficient of the resonant term.
Therefore proceeding by induction on $n$
we get the assertion.\qedn

\newdef  The formal invariant $a\in\C$ is called {\sl resonant index.}

Our next aim is to prove that, for Fuchsian singularities, the formal normal forms given
by Theorem~\stuno.(i) are actually holomorphic normal forms, that is that
we can find a holomorphic change of coordinates of the form \eqallcc\ bringing $G$ in the given
normal form. To do so, we shall adapt the holomorphic Poincar\'e-Dulac paradigm as described in [IY, Section~I.5].

\newdef Given $r>0$, the {\sl majorant $r$-norm} of a formal power
series in $\C[\![z]\!]$ is defined by
$$
\left\|\sum_{j=0}^\infty a_j z^j\right\|_r=\sum_{j=0}^\infty |a_j|r^j\;.
$$
When $h=(h_1,h_2)\in\C[\![z]\!]^2$ we set
$$
\|h\|_r=\|h_1\|_r+\|h_2\|_r\;.
$$
We shall denote by $\ca B_r$ (respectively, $\ca B_r^2$) the space of $h\in\C[\![z]\!]$
(respectively,  of $h\in\C[\![z]\!]^2$) with finite
majorant $r$-norm, and by $\ca B_{r,\ell}$ (respectively, $\ca B_{r,\ell}^2$) the subspace of
elements of order at least~$\ell$. It is easy to check (see [IY, 
Proposition~5.8]) that $\ca B_r$, $\ca B_{r,\ell}$, $\ca B_r^2$ and $\ca B_{r,\ell}^2$ are Banach spaces.

Majorant norms are multiplicative, that is
$$
\|fg\|_r\le\|f\|_r\cdot\|g\|_r
\neweq\eqsttre
$$
for all $f$, $g\in\ca B_r$; see [IY, Lemma~5.10].

Clearly if $r'\le r$ we have $\|h\|_{r'}\le\|h\|_r$, and so the natural
inclusion $\ca B_{r,\ell}^2\hookrightarrow\ca B_{r',\ell}^2$ is continuous.

\newdef Let $S\colon\ca B_{r,\ell}^2\to\ca B_{r,\ell}^2$ be an operator
defined on $\ca B_{r,\ell}^2$ for fixed $\ell\in\N$ and all $r$ small enough (and commuting with the inclusions $\ca B_{r,\ell}^2\hookrightarrow\ca 
B_{r',\ell}^2$). We shall say that $S$ is
{\sl strongly contracting} if
\smallskip
\item{(i)} $\|S(0)\|_r=O(r^2)$, and
\item{(ii)} $S$ is Lipschitz on the ball $B_{r,\ell}^2=\{h\in\ca B_{r,\ell}^2
\mid \|h\|_r\le r\}$, with Lipschitz constant no greater than $O(r)$ as~$r\to 0$.

The point of strongly contracting operators is that, for $r$ small
enough, they are a contraction of $B_{r,\ell}^2$ into itself, and hence 
admit a unique fixed point there. 

The next lemma contains examples of strongly contracting operators
we shall need later on.

\newthm Lemma \stdue: Given $\ell\in\N$, the following operators are strongly contracting:
\smallskip
\itm{(i)} $P_j\colon\ca B_{r,\ell}^2\to\ca B_{r,\ell}^2$ given by
$P_j(h_1,h_2)=(h_1^j,0)$, if $j\ge 2$;
\itm{(ii)} $Q_j\colon\ca B_{r,\ell}^2\to\ca B_{r,\ell}^2$ given by
$Q_j(h_1,h_2)=(h_2h_1^j,0)$, if $j\ge 1$;
\itm{(iii)} $R_{A,t}\colon\ca B_{r,\ell}^2\to\ca B_{r,\ell}^2$ given by
$R_{A,t}(h)=z^t A(h)$, where $t\ge 1$ and $A\colon\ca B_{r,\ell}^2\to
\ca B_{r,\ell}^2$ satisfies $\|A(O)\|_r=O(r)$ and $\|A(h)-A(k)\|_r\le O(1)\|h-k\|_r$;
\itm{(iv)} $S_g\colon\ca B_{r,\ell}^2\to\ca B_{r,\ell}^2$ given by
$S_g(h_1,h_2)=\bigl(g_1\circ(\id +h_1),g_2\circ(\id+h_1)\bigr)$, where 
$\ell\ge 2$ and $g=(g_1,g_2)\in\ca B_{s,\ell}^2$ for $s$ small enough;
\itm{(v)} $T_{p_1,p_2,S}\colon\ca B_{r,\ell}^2\to\ca B_{r,\ell}^2$ given by
$T_{p_1,p_2,S}(h_1,h_2)=\bigl(p_1(h_1,h_2)S_1(h_1,h_2),p_2(h_1,h_2)
S_2(h_1,h_2)\bigr)$, where $p_1$, $p_2$ are polynomials, and 
$S=(S_1,S_2)\colon\ca B_{r,\ell}^2\to\ca B_{r,\ell}^2$ is strongly
contracting;
\itm{(vi)} linear combinations of strongly contracting operators;
\itm{(vii)} operators of the form $A\circ 
S\colon\ca B_{r,\ell}^2\to\ca B_{r,\ell}^2$, with 
$S\colon\ca B_{r,\ell}^2\to\ca B_{r,\ell}^2$ strongly contracting 
and $A\colon\ca B_{r,\ell}^2\to\ca B_{r,\ell}^2$\break\indent
linear (commuting with the inclusions) so that $\|Ah\|_r\le O(1)\|h\|_r$.

\pf (i) It suffices to check the Lipschitz constant of $P_j$ on 
$B_{r,\ell}^2$. Using \eqsttre\ we have
$$
\|P_j(h_1,h_2)-P_j(k_1,k_2)\|_r=\|h_1^j-k_1^j\|_r\le
\|h_1-k_1\|_r\sum_{s=0}^{j-1}\|h_1\|^s_r\|k_1\|^{j-s-1}_r
\le jr^{j-1}\|h_1-k_1\|_r\;,
$$
and so the Lipschitz constant of~$P_j$ is~$O(r^{j-1})$ on $B_{r,\ell}^2$.

(ii) Analogously, if $(h_1,h_2)$, $(k_1,k_2)\in B_{r,\ell}^2$ we get
$$
\eqalign{\|Q_j(h_1,h_2)-Q_j(k_1,k_2)\|_r&=\|h_2h_1^j-k_2k_1^j\|_r\le
\|h_2\|_r\|h_1^j-k_1^j\|_r+\|h_2-k_2\|_r\|k_1\|_r^j\cr
&\le jr^j\|h_1-k_1\|_r+r^j\|h_2-k_2\|_r\le jr^j\|(h_1,h_2)-(k_1,k_2)\|_r\;,
\cr}
$$
and so the Lipschitz constant of~$Q_j$ is~$O(r^j)$ on $B_{r,\ell}^2$.

(iii) It follows immediately from $\|z^t h\|_r=r^t\|h\|_r$.

(iv) See [IY, Lemma~5.14].

(v) First of all,
$$
\|T_{p_1,p_2,S}(O)\|_r=\|p_1(O)S_1(O)\|_r+\|p_2(O)S_2(O)\|_r\le 
\max\{|p_1(O)|,|p_2(O)|\}\|S(O)\|_r=O(r^2)\;.
$$
Furthermore for $j=1$, 2 we have
$$
\eqalign{\|p_j(h)S_j(h)-p_j(k)S_j(k)\|_r&
\le\|p_j(h)\|_r\|S_j(h)-S_j(k)\|_r+\|p_j(h)-p_j(k)\|_r\|S_j(k)\|_r\cr
&\le O(1)\|S_j(h)-S_j(k)\|_r\cr
&\quad+\|p_j(h)-p_j(k)\|_r\|S_j(k)-S_j(O)\|_r
+\|p_j(h)-p_j(k)\|_r\|S_j(O)\|_r\cr
&\le O(r)\|h-k\|_r+O(1)\|h-k\|_r O(r)\|k\|_r+O(1)\|h-k\|_rO(r^2)\cr
&\le O(r)\|h-k\|_r\;,
\cr}
$$
where we used the fact that $\|p_j(h)\|_r$ is bounded on $B_{r,\ell}^2$,
and that 
$$
p_j(h)-p_j(k)=(h_1-k_1)q_1(h,k)+(h_2-k_2)q_2(h,k)
$$
for suitable polynomials $q_1$ and $q_2$, and the assertion follows.

(vi) and (vii) are obvious.\qedn

We are now ready to prove the main theorem of this section:

\newthm Theorem \sttre: Let $E$ be a line bundle on
a Riemann surface $S$, and assume we have  
a morphism $X\colon E\to TS$ which is an isomorphism on $S^o=X\setminus\Sing(X)$,
and a meromorphic connection
$\nabla$ on $E$, holomorphic on~$E|_{S^o}$, such that the geodesic field $G$
extends holomorphically from $E|_{S^o}$ to the whole of~$E$. 
Let $p_0\in\Sing(X)$ be a Fuchsian singularity, and in a chart
$(U_\alpha,z_\alpha,e_\alpha)$ centered at~$z_0$ write
$$
G=z_\alpha^{\mu_X}(a_0+a_1 z_\alpha+\cdots)v_\alpha\de_\alpha-
z_\alpha^{\mu_Y}(b_0+b_1 z_\alpha+\cdots)v_\alpha^2{\de\over\de v_\alpha}\;,
$$
with $\mu_X=\mu_Y+1$ and $a_0$, $b_0\ne 0$. Put $\rho=b_0/a_0=\Res_{p_0}(\nabla)$. Then we can
find a chart $(U,z,e)$ centered in~$p_0$ in which $G$ is
given by
$$
z^{\mu_X-1}\left(z v\de-\rho v^2{\de\over\de v}\right)
$$
if $\mu_Y-\rho\notin\N^*$, or by
$$
z^{\mu_X-1}\left(zv\de-\rho v^2(1+a z^n){\de\over\de v}\right)
$$
for a suitable $a\in\C$ if $n=\mu_Y-\rho\in\N^*$.

\pf By (the proof of) Theorem~\stuno, we can find a chart
$(U_\beta,z_\beta,e_\beta)$ centered in~$p_0$ where $G$ is given by
$G=X_\beta v_\beta\de_\beta-Y_\beta v_\beta^2{\de\over\de v_\beta}$ with
$$
X_\beta=z^{\mu_X}_\beta\bigl(1+g_1(z_\beta)\bigr)\quad
\hbox{and}\quad Y_\beta=\rho z_\beta^{\mu_X-1}\bigl(1+a z^n_\beta+g_2(z_\beta)\bigr)\;,
$$
where $\ord_{p_0}(g_1)$, $\ord_{p_0}(g_2)>n$ and $a\ne 0$ only if
$n=\mu_X-1-\rho\in\N^*$. Recalling \eqsdue, we need to find $\psi$
and $\xi$ with $\psi(0)=0$, $\psi'(0)=\xi(0)=1$ such that
$$
X_\beta\bigl(\psi(z)\bigr)={\psi'(z)\over\xi(z)}\,z^{\mu_X}
\quad\hbox{and}\quad
Y_\beta\bigl(\psi(z)\bigr)={1\over\xi(z)}\rho z^{\mu_X-1}(1+a z^n)-{\xi'(z)\over
\xi(z)^2}z^{\mu_X}\;.
\neweq\eqstqua
$$
Writing $\psi(z)=z+zh_1(z)$ and $\xi(z)=1+h_2(z)$ with $(h_1,h_2)\in\ca
B_{r,\ell}^2$ for $r$ small enough and $\ell\ge 1$ to be determined later,
we can reformulate \eqstqua\ as follows:
$$
\cases{\displaystyle
(1+h_2)z^{\mu_X}(1+h_1)^{\mu_X}\bigl[1+g_1(z+zh_1)\bigr]=
(1+zh_1'+h_1)z^{\mu_X}\;,\cr\displaystyle
(1+h_2)^2\rho z^{\mu_X-1}(1+h_1)^{\mu_X-1}\bigl[1+az^n(1+h_1)^n+g_2(z+zh_1)\bigr]=
(1+h_2)\rho z^{\mu_X-1}(1+az^n)-h'_2 z^{\mu_X}\;.
\cr}
$$
Simplifying the powers of~$z$, expanding the powers of $1+h_1$ and $1+h_2$ 
and rearranging the terms we get
$$
\cases{\displaystyle
zh'_1+(1-\mu_X)h_1-h_2=\sum_{j=2}^{\mu_X}{\mu_X\choose j}h_1^j+h_2
\sum_{j=1}^{\mu_X}{\mu_X\choose j}h_1^j+(1+h_2)(1+h_1)^{\mu_X}g_1(z+zh_1)\;,\cr
\noalign{\smallskip}\displaystyle
\eqalign{-zh'_2+\rho(1&-\mu_X)h_1-\rho h_2\cr
&=\rho\sum_{j=2}^{\mu_X-1}
{\mu_X-1\choose j}h_1^j+2\rho h_2\sum_{j=1}^{\mu_X-1}{\mu_X-1\choose j}h_1^j
\cr
&\quad+\rho a z^n\left[(\mu_X+n-1)h_1+h_2+\sum_{j=2}^{\mu_X+n-1}{\mu_X+n-1\choose j}
h_1^j+2h_2\!\sum_{j=1}^{\mu_X+n-1}{\mu_X+n-1\choose j}h_1^j\right]
\cr
&\quad+\rho(1+h_1)^{\mu_X-1}\bigl(1+az^n(1+h_1)^n\bigr)h_2^2+\rho(1+h_2)^2(1+h_1)^{\mu_X-1}g_2(z+zh_1)\;.
\cr}
\cr}
$$
We can write this in a more compact form as
$$
Ah=B_1h+B_2h+Ch+Dh+Eh\;,
$$
where $A$, $B_1$, $B_2$, $C$, $D$, and $E$ are operators on $\ca 
B^2_{r,\ell}$ respectively given by
$$
\displaylines{
A(h_1,h_2)=\bigl(zh'_1+(1-\mu_X)h_1-h_2,-zh'_2+\rho(1-\mu_X)h_1-\rho h_2
\bigr)\;,\cr
B_1(h_1,h_2)=\left(\sum_{j=2}^{\mu_X}{\mu_X\choose j}h_1^j,
\rho\sum_{j=2}^{\mu_X-1}
{\mu_X-1\choose j}h_1^j\right)\;,\cr
B_2(h_1,h_2)=\left(h_2\sum_{j=1}^{\mu_X}{\mu_X\choose j}h_1^j,
2\rho h_2\sum_{j=1}^{\mu_X-1}{\mu_X-1\choose j}h_1^j\right)\;,\cr
C(h_1,h_2)=\left(0,\rho az^n\left[(\mu_X+n-1)h_1+h_2+
\sum_{j=2}^{\mu_X+n-1}{\mu_X+n-1\choose j}h_1^j+h_2\sum_{j=1}^{\mu_X+n-1}
{\mu_X+n-1\choose j}h_1^j\right]\right)\;,\cr
D(h_1,h_2)=\left(0,\rho(1+h_1)^{\mu_X-1}\bigl(1+az^n(1+h_1)^n\bigr)h_2^2\right)\;,\cr
E(h_1,h_2)=\bigl((1+h_2)(1+h_1)^{\mu_X}g_1(z+zh_1),\rho
(1+h_2)^2(1+h_1)^{\mu_X-1}g_2(z+zh_1)\bigr)\;.
\cr}
$$
The operators $B_1$, $B_2$, $C$, $D$ and $E$ are linear combinations
of strongly contracting operators, and hence are strongly 
contracting (by Lemma~\stdue). 
The operator $A$ is linear, and it preserves the degrees: we have
$$
A(c_1 z^d, c_2 z^d)=\bigl(\bigl[(d+1-\mu_X)c_1-c_2\bigr]z^d,
\bigl[\rho(1-\mu_X)c_1-(d+\rho)c_2\bigr]z^d\bigr)\;.
$$
So if $d\ne\mu_X-1-\rho$ we get
$$
A^{-1}(c_1 z^d,c_2 z^d)={1\over d(d+\rho+1-\mu_X)}
\Bigl(\bigl[(d+\rho)c_1-c_2\bigr]z^d,\bigl[\rho(1-\mu_X)c_1-(d+1-\mu_X)c_2)
\bigr]z^d\Bigr)\;.
$$
In particular,
$$
\|A^{-1}(c_1 z^d,c_2 z^d)\|_r\le O(1/d) \|(c_1 z^d,c_2z^d)\|_r\;,
\neweq\eqstcin
$$
uniformly on~$r$ and $d$. So $A$ is invertible on 
$\ca B^2_{r,\ell}$ as soon as $\ell>\mu_X-1-\rho=n$ (or $\ell>1$ if 
$\mu_X-1-\rho\notin\N$), and $\|A^{-1}\|_r=O(1)$
uniformly on~$r$.

Summing up, to solve \eqstqua\ is enough to solve the fixed point
problem
$$
h=A^{-1}B_1h+A^{-1}B_2h+A^{-1}Ch+A^{-1}Dh+A^{-1}Eh
$$
on $B^2_{r,\ell}$ for $r$ small enough and $\ell>n$ (or $\ell>1$ when
$a=0$). Since $A^{-1}B_1$, $A^{-1}B_2$, $A^{-1}C$, $A^{-1}D$ and
$A^{-1}E$ are strongly contracting by Lemma~\stdue.(vii), this is possible, and we are done.\qedn


\newrem This approach does not work in the irregular case.
Working as in the previous proof we reduce the problem to solving 
an equation of the form $Ah=B_1h+B_2h+Ch+Dh+Eh+Fh$, where
$B_1$, $B_2$, $C$, $D$ and $E$ are as before (just replace $n$ by
$m-1$ and $\mu_X-1$ by $\mu_X-m$ everywhere), whereas $A$ and $F$
are given by
$$
\displaylines{
A(h_1,h_2)=\bigl(zh'_1+(1-\mu_X)h_1-h_2,\rho(m-\mu_X)h_1-\rho h_2
\bigr)\;,\cr
F(h_1,h_2)=(0,z^m h'_2)\;.
\cr}
$$
In particular
$$
A^{-1}(c_1 z^d,c_2 z^d)={1\over \rho(d+1-m)}
\Bigl(\bigl[\rho c_1-c_2\bigr]z^d,\bigl[\rho(m-\mu_X)c_1-(d+1-\mu_X)c_2)
\bigr]z^d\Bigr)
$$
and
$$
A^{-1}F(c_1 z^d,c_2 z^d)={-c_2\over \rho}
\left(z^{d+m-1},(d-\mu_Y)z^{d+m-1}\right)\;,
$$
and hence $A^{-1}F$ is not strongly contracting.

As a corollary, we can describe the behavior of the geodesics nearby
non-resonant Fuchsian singularities:

\newthm Proposition \stqua: Let $E$ be a line bundle on
a Riemann surface $S$, and assume we have  
a morphism $X\colon E\to TS$ which is an isomorphism on $S^o=X\setminus\Sing(X)$,
and a meromorphic connection
$\nabla$ on $E$, holomorphic on~$E|_{S^o}$, such that the geodesic field $G$
extends holomorphically from $E|_{S^o}$ to the whole of~$E$. 
Let $p_0\in\Sing(X)$ be a Fuchsian singularity of order~$\mu_X\ge 1$, 
and with vanishing
resonant index if $\mu_X-1-\rho\in\N^*$,  where $\rho=\Res_{p_0}(\nabla)$. 
Put $\mu_Y=\mu_X-1$. Then
there is a neighbourhood $U\subseteq S$ of~$p_0$ such that:
\smallskip
\itm{(i)} if $\Re\rho<\mu_Y$ then all geodesics but one issuing from any point $p\in U\cap 
S^o$ 
tend to~$p_0$ staying inside $U$ (the only exception escapes $U$); furthermore
for every geodesic $\sigma$ going to $p_0$ inside $U$ we have
\itemitem{\rm(a)} if $\mu_Y\Re\rho<|\rho|^2$ then 
$X^{-1}\bigl(\sigma'(t)\bigr)\to O$ as $\sigma(t)\to p_0$;
\itemitem{\rm(b)} if $\mu_Y\Re\rho>|\rho|^2$ then
$\bigl|X^{-1}\bigl(\sigma'(t)\bigr)\bigr|\to +\infty$ as $\sigma(t)\to p_0$;
\itemitem{\rm(c)} if $\mu_Y\Re\rho=|\rho|^2$ then 
$X^{-1}\bigl(\sigma'(t)\bigr)$ accumulates a circumference in~$E_{p_0}$;
\itm{(ii)} if $\Re\rho>\mu_Y$ then all geodesics but one issuing from any point~$p\in
U\cap S^o$ escapes $U$; furthermore, the exceptional geodesic $\sigma_0$ 
tends to~$p_0$ in finite time with 
$\bigl|X^{-1}\bigl(\sigma'_0(t)\bigr)\bigr|\to+\infty$ as $\sigma_0(t)\to p_0$;
\itm{(iii)} if $\Re\rho=\mu_Y$ but $\rho\ne\mu_Y$ then the geodesics not
escaping $U$
are either closed (with $X^{-1}(\sigma')$ either tending to~$O$ or
diverging to infinity) or accumulate the support of a closed geodesic in~$U$
(with $X^{-1}(\sigma')$ tending to $O$);
\itm{(iv)} if $\rho=\mu_Y$ then for all $p\in S^o\cap U$ there is a 
non-zero $v_p\in E_p$ such that
\itemitem{\rm(a)} if $v=\zeta v_p\in E_p$ with $\Re\zeta<0$ then the geodesic~$\sigma_v$ issuing from $p$ tangent to~$X(v)$ converges to~$p_0$ staying in~$U$
but with $\bigl|X^{-1}\bigl(\sigma'_v(t)\bigr)\bigr|\to+\infty$;
\itemitem{\rm(b)} if $v=\zeta v_p\in E_p$ with $\Re\zeta>0$ then the geodesic~$\sigma_v$ 
issuing from $p$ tangent to~$X(v)$ escapes from $U$;
\itemitem{\rm(c)} if $v=\zeta v_p\in E_p$ with $\Re\zeta=0$ then the geodesic~$\sigma_v$ 
issuing from $p$ tangent to~$X(v)$ is periodic and\break\indent\indent surrounds $p_0$.

\pf By Theorem~\sttre\ we know that, in a suitable local chart $(U,z,e)$
centered in~$p_0$, a curve $\sigma$ is a geodesic if and only if
$X^{-1}\bigl(\sigma'(t)\bigr)=\bigl(z(t),v(t)\bigr)$ satisfies
$$
\cases{z'=z^{\mu_X}v\;,\cr
v'=-\rho z^{\mu_Y}v^2\;,
\cr}
\neweq\eqstsei
$$
where $\rho=\Res_{p_0}(\nabla)$ and $\mu_Y=\mu_X-1$. 

Assume $\rho\ne\mu_Y$. Then the solution of \eqstsei\ is
$$
\cases{\displaystyle
z(t)=z_0(1+c t)^{1/(\rho-\mu_Y)}=z_0\exp\left({1\over\rho-\mu_Y}\log(1+ct)
\right)\;,\cr\displaystyle
v(t)=v_0(1+c t)^{-\rho/(\rho-\mu_Y)}=v_0\exp\left({-\rho\over\rho-\mu_Y}
\log(1+ct)\right)\;,
\cr}
$$ 
with $c=(\rho-\mu_Y)z_0^{\mu_Y}v_0$, where we have chosen the principal
determination of the logarithm (that is $\log 1=0$). In particular,
$$
|z(t)|=|z_0|\exp\left[\Re\left({1\over\rho-\mu_Y}\right)\log|1+ct|
-\Im\left({1\over\rho-\mu_Y}\right)\arg(1+ct)\right]
$$
and 
$$
|v(t)|=|v_0|\exp\left[\Re\left({-\rho\over\rho-\mu_Y}\right)\log|1+ct|
-\Im\left({-\rho\over\rho-\mu_Y}\right)\arg(1+ct)\right]\;.
$$
Notice that $\arg(1+ct)$ is always bounded.  

Suppose $\Re\rho<\mu_Y$, so that $\Re(\rho-\mu_Y)^{-1}<0$. Given $z_0$, if
$v_0$ is such that $(\rho-\mu_Y)z_0^{\mu_Y}v_0\notin\R^-$ then 
$\sigma(t)$ is defined for all $t>0$ and $\sigma(t)\to p_0$ as $t\to+\infty$. Moreover,
$\rho\ne 0$ because $p_0$ is Fuchsian, and
$$
\Re\left({-\rho\over\rho-\mu_Y}\right)={-|\rho|^2+\mu_Y\Re\rho\over|\rho-
\mu_Y|^2}\;;
$$
hence (a)--(c) follow. If instead $v_0$ is such that $c=(\rho-\mu_Y)z_0^{\mu_Y}v_0\in\R^-$,
then $|z(t)|$ explodes as $t\to -c^{-1}$, which means that $\sigma(t)$ escapes $U$;
so (i) is proved.

If instead $\Re\rho>\mu_Y$ the situation is reversed: if $(\rho-\mu_Y)z_0^{\mu_Y}v_0
\notin\R^-$ then the geodesic escapes, while if $c=(\rho-\mu_Y)z_0^{\mu_Y}v_0\in\R^-$
then $\sigma(t)\to p_0$ as $t\to-c^{-1}$. Moreover, $\Re\rho>\mu_Y$ implies
$|\rho|^2\ge(\Re\rho)^2>\mu_Y\Re\rho$, and so $|v(t)|\to+\infty$
as $t\to -c^{-1}$. This completes the proof of (ii).

In case (iii) we have $1/(\rho-\mu_Y)=i\gamma\in
i\R^*$, and thus
$$
\cases{\displaystyle
z(t)=z_0\exp\bigl[-\gamma\arg(1+ct)+i\gamma\log|1+ct|\bigr]\;,\cr\displaystyle
v(t)=v_0\exp\bigl[-\log|1+ct|+\mu_Y\gamma\arg(1+ct)\bigr]\exp\bigl[-i
\bigl(\arg(1+ct)+\mu_Y\gamma\log|1+ct|\bigr)\bigr]\;.
\cr}
$$ 
If $c=-i\gamma^{-1}z_0^{\mu_Y}v_0\in\R^-$ (respectively, $c\in\R^+$) then
$\arg(1+ct)$ is constant, and we get a closed geodesic, with $v(t)$ exploding as $t\to -c^{-1}$ (respectively, with $v(t)\to O$ as
$t\to+\infty$). 
If instead $c\notin\R^*$ we get a geodesic accumulating a circumference
of positive radius, which is the support of the geodesic issuing from
$z_1=z_0\exp(-\gamma\arg c)$ in the direction $v_1=-iz_0^{-\mu_Y}$. We also have $v(t)\to O$. Notice that $\Re\rho=\mu_Y$ implies that the induced residue
has real part~$-1$, in accord with the results of Section~4.

Finally, when $\rho=\mu_Y>0$ (when $\mu_Y=0$ we have $\rho\ne0$ by definition), the solution of \eqstsei\ is
$$
\cases{\displaystyle
z(t)=z_0\exp\left(z_0^{\mu_Y}v_0 t\right)\;,\cr\displaystyle
v(t)=v_0\exp\left(-\mu_Y z_0^{\mu_Y}v_0 t\right)\;.
\cr}
$$ 
In this case, $\sigma(t)\to p_0$ as $t\to+\infty$ if and only if 
$\Re(z_0^{\mu_Y}v_0)<0$, and then $|v(t)|\to+\infty$. If $\Re(z_0^{\mu_Y}v_0)>0$ then
the geodesic escapes, and if $\Re(z_0^{\mu_Y}v_0)=0$ then the geodesic is periodic,
completing the proof of (iv).\qedn

\newrem In the resonant case, a curve $\sigma$ is a geodesic if
and only if, in suitable local coordinates, $X^{-1}\bigl(\sigma'(t)\bigr)=\bigl(z(t),v(t)\bigr)$ satisfies
$$
\cases{z'=z^{n+\rho+1}v\;,\cr
v'=-\rho z^{n+\rho}(1+a z^n) v^2\;,
\cr}
$$
where $a\in\C^*$, $n\in\N^*$ and $\rho\in\Z$ with $\rho\ge-n$. Then
$$
{v'\over v}=-\rho\left({1\over z}+az^{n-1}\right)z'\qquad\Longrightarrow
\qquad v=c_0 z^{-\rho}\exp\left(-{\rho a\over n}z^n\right)
$$
with 
$$
c_0=v_0 z_0^\rho\exp\left({\rho a\over n} z_0^n\right)\ne 0\;,
$$
where $(z_0,v_0)=\bigl(z(0),v(0)\bigr)$. So $z(t)$ satisfies
$$
z'=c_0 z^{n+1}\exp\left(-{\rho a\over n}z^n\right)\;.
$$
Setting $w={\rho a\over n}z^n$ we find
$$
w'={n^2\over\rho a}c_0 w^2 e^{-w}\qquad\Longrightarrow\qquad
-{e^{w(t)}\over w(t)}+\hbox{\rm Ei}\bigl(w(t)\bigr)={n^2 c_0\over\rho a}t+c_1\;,
$$
where Ei$(w)$ is the holomorphic primitive of $w^{-1}e^w$ vanishing at 
$w_0=w(0)$, and $c_1=-w_0^{-1}e^{w_0}$. A numerical study of 
this equation suggests that for every $z_0$ one has 
$w(t)\to O$ (and hence $z(t)\to O$) as $t\to+\infty$ for an open (and possibly dense) set of~$v_0\in\C^*$. So we conjecture that 
Proposition~\stqua.(i) holds in this case too (in the resonant case
$\Re\rho=\rho<\mu_Y$ necessarily).

\newrem Notice that a Fuchsian singularity with $\Re\rho<\mu_Y$ cannot appear as a vertex in a simple cycle
of saddle connections which is accumulated by a geodesic. Indeed, such a behaviour
requires the existence of geodesics arbitrarily close to the singularity and escaping in both forward and backward time. So case (iv) of Theorem~\qqua\ cannot involve Fuchsian
singularities --- or apparent singularities of order~1, for the same reason.

\Figuraeps 1 (geod1)

{\sc Example 8.1.} A computation with Mathematica shows that,
if we take $\mu_X=1$ and $\rho=0.1$, then the geodesic
issuing from $z_0=1$ in the direction~$v_0=1+i$ intersect itself
twice before escaping to infinity; see Fig.~1. On the other hand,
if we take $\mu_X=1$ and $\rho=i$, the geodesic
issuing from $z_0=(1+i)/2$ in the direction~$v_0=1$ accumulates a
closed geodesic; see Fig.~2.

\Figurascaledeps 2 (geod3)(500)

Translating Proposition~\stqua\ to the case of homogeneous vector fields we get:

\newthm Corollary \stcin: Let $Q$ be a homogeneous vector
field on $\C^2$ of degree~$\nu+1\ge 2$. 
Let $[v_0]\in\P^1(\C)$ be a Fuchsian singularity  
of~$X_Q$ of order~$\mu_X\ge 1$, residue $\rho\in\C^*$ (and resonant index $a\in\C$ if 
$\mu_X-1-\rho\in\N^*$). Put $\mu_Y=\mu_X-1$. Then:
\smallskip
\itm{(i)} if the direction $[\gamma(t)]\in\P^1(\C)$ of an integral curve
$\gamma\colon[0,\eps)\to\C^2\setminus\{O\}$ of $Q$ tends to $[v_0]$
as $t\to\eps$ and $\gamma$ is not contained in the characteristic leaf $L_{v_0}$ then 
\itemitem{\rm(a)} if $\Re\rho<\mu_Y$ and $\mu_Y\Re\rho<|\rho|^2$ then $\gamma(t)$ tends to the origin;
\itemitem{\rm(b)} if $\rho=\mu_Y>0$, or $\Re\rho>\mu_Y$, or $\Re\rho<\mu_Y$ and $\mu_Y\Re\rho>|\rho|^2$, then $\|\gamma(t)\|$ tends to $+\infty$;
\itemitem{\rm(c)} if $\Re\rho<\mu_Y$ and $\mu_Y\Re\rho=|\rho|^2$ then $\gamma(t)$ accumulates a
circumference in $L_{v_0}$.
\smallskip
\noindent Furthermore there is
a neighbourhood $U\subset\P^1(\C)$ of $[v_0]$ such that an integral curve
$\gamma$ issuing from a point $z_0\in\C^2\setminus L_{v_0}$ with 
$[z_0]\in U\setminus\{[v_0]\}$ 
can have one of the following behaviors, where $\hat U=\{z\in\C^2
\setminus\{O\}\mid[z]\in U\}$:
\smallskip
\itm{(ii)} if $\Re\rho>\mu_Y$ then
\smallskip
\itemitem{\rm (a)}either $\gamma(t)$ escapes $\hat U$, and this happens
for a Zariski open dense set of initial conditions in~$\hat U$; or
\itemitem{\rm(b)} $[\gamma(t)]\to[v_0]$ but $\|\gamma(t)\|\to+\infty$;
\smallskip
in particular, no integral curve outside $L_{v_0}$ converge to the origin tangent to~$[v_0]$;
\itm{(iii)} if $\Re\rho=\mu_Y$ but $\rho\ne\mu_Y$ then 
\smallskip
\itemitem{\rm(a)} either $\gamma(t)$ escapes $\hat U$; or
\itemitem{\rm(b)} $\gamma(t)\to O$ without being tangent to any direction,
and $[\gamma(t)]$ is a closed curve or accumulates a closed
curve in~$\P^1(\C)$ surrounding $[v_0]$; or
\itemitem{\rm(c)} $\|\gamma(t)\|\to+\infty$ without being tangent to any direction,
and $[\gamma(t)]$ is a closed curve in~$\P^1(\C)$ surrounding $[v_0]$;
\smallskip
in particular, no integral curve outside $L_{v_0}$ converge to the origin tangent to~$[v_0]$;
\itm{(iv)} if $\rho=\mu_Y>0$ then 
\smallskip
\itemitem{\rm(a)} either $\gamma(t)$ escapes $\hat U$, and this 
happens for an open set $\hat U_1\subset\hat U$ of initial conditions; or
\itemitem{\rm(b)} $[\gamma(t)]\to[v_0]$ with $\|\gamma(t)\|\to+\infty$,
and this happens for an open set $\hat U_2\subset\hat U$ of initial conditions
such that $\hat U_1\cup\hat U_2$ is dense in $\hat U$; or
\itemitem{\rm(c)} $\gamma$ is a periodic integral curve with $[\gamma]$
surrounding $[v_0]$;
\smallskip
in particular, no integral curve outside $L_{v_0}$ converge to the origin tangent to~$[v_0]$, but we have periodic\break\indent integral curves of arbitrarily long period accumulating the origin; 
\itm{(v)} if $\Re\rho<\mu_Y$ and $a=0$ then $[\gamma(t)]\to[v_0]$
for an open 
dense set~$\hat U_0$ of initial conditions in~$\hat U$, and $\gamma$
escapes $\hat U$ for $z\in\hat U\setminus\hat U_0$; more precisely,
\smallskip
\itemitem{\rm (a)} if $\mu_Y\Re\rho<|\rho|^2$ then $\gamma(t)\to O$
tangent to~$[v_0]$ for all $z\in\hat U_0$;
\itemitem{\rm (b)} if $\mu_Y\Re\rho>|\rho|^2$ then $|\gamma(t)|\to
+\infty$ tangent to $[v_0]$ for all $z\in\hat U_0$;
\itemitem{\rm(c)} if $\mu_Y\Re\rho=|\rho|^2$ then $\gamma(t)$
accumulates a circumference in $L_{v_0}$.

\pf Notice that $\mu_Y-\rho\in\N^*$ implies $\mu_Y>\rho\in\Z$, and so
in cases (ii), (iii) and (iv) the resonant index vanishes by 
definition. Then
the only part that does not follow immediately from Proposition~\stqua\ is
part (i) when $a\ne 0$. But in that case $\mu_Y\Re\rho<|\rho|^2$
if and only if $\rho<0$, and (c) can never happen. The assertion then follows from 
Proposition~\squa.\qedn

\newrem This result must be compared with a theorem due to Hakim [H2], saying that if $[v_0]$ is a non-degenerate characteristic direction whose director has positive real part then there is an open set of points whose orbits converge to the origin tangentially to~$[v_0]$. A non-degenerate characteristic direction $[v_0]$ with non-zero
director~$\delta$ is a Fuchsian singularity of order~1. Then 
Corollary~\stcin\ says that if $\Re\delta<0$ (that is $\Re\rho>0=\mu_Y$)
then no orbit of the time-1 map of $Q$ outside of~$L_{v_0}$
(that is, outside of the parabolic curve whose existence is
ensured by \'Ecalle and Hakim's results [\'E1--4], [H1]) converges to the origin tangent to~$[v_0]$,
whereas if $\Re\delta>0$ (that is $\Re\rho<0=\mu_Y$) and $a=0$ then
the orbits under the time 1-map of $Q$ converge to the origin tangent
to~$[v_0]$ for an open (and dense in a conical neighbourhood of~$[v_0]$)
set of initial conditions.

\newrem If the conjecture mentioned in Remark~8.2 is true then Corollary~\stcin.(v) holds in the resonant case too.

\newrem Since (as already observed in the proof of Corollary~\stre) the 
periods of the periodic integral curves in Corollary~\stcin.(iv) tend
continuously to infinity as the curves approach the origin, this yields
for the time-1 map of $Q$ both periodic orbits accumulating at the
origin (when the period of the periodic integral curve is rational)
and orbits whose closure is a circle accumulating the origin 
(when the period of the periodic integral curve is irrational). 

Putting together the Poincar\'e-Bendixson theorems discussed in Section~4
together with the local results in this and the previous sections, one can now say a lot about the dynamics of homogeneous vector fields. In the
next section we shall discuss in detail quadratic vector fields; we end
this section with an example of application of our results giving a complete description of the dynamics of a substantial class of homogeneous vector fields:

\newthm Corollary \stsei: Let $Q$ be a non-dicritical 
homogeneous vector field on $\C^2$ of degree~$\nu+1\ge 2$. Assume
that all characteristic directions of $Q$ are non-degenerate with non-zero director (or, equivalently, that they are Fuchsian singularities of
order~$1$). Assume moreover that for no set of $g\ge 1$ characteristic
directions the real part of the sum of the residues is equal to~$g-1$.\hfil\break\indent
Let $\gamma\colon[0,\eps_0)\to\C^2$ be a maximal integral curve of $Q$. Then:
\smallskip
\itm{(a)} If $\gamma(0)$ belongs to a characteristic leaf~$L_{v_0}$, then the image of $\gamma$ is contained 
in~$L_{v_0}$. Moreover, either $\gamma(t)\to O$ (and this happens for
a Zariski open dense set of initial conditions in~$L_{v_0}$), or $\|\gamma(t)\|\to+\infty$.
\itm{(b)} If $\gamma(0)$ does not belong to a characteristic leaf then either
\itemitem{\rm (i)} $\gamma$ converges to the origin tangentially to a characteristic 
direction $[v_0]$ whose residue has negative real part (and hence 
whose director has positive real part); or
\itemitem{\rm (ii)} $\|\gamma(t)\|\to+\infty$ tangentially to a characteristic
direction $[v_0]$ whose residue has positive real part (and hence
whose director has negative real part); or
\itemitem{\rm(iii)} $[\gamma]\colon[0,\eps_0)\to\P^1(\C)$ intersects itself
infinitely many times. 
\smallskip
Furthermore, if {\rm (iii)} never occurs then {\rm (i)} holds for a Zariski open dense set
of initial conditions.

\pf Statement (a) follows immediately from Lemma~\dcl.

For (b), first of all notice that
by assumption all characteristic directions of $Q$ are Fuchsian singularities
of order~1 (see Remark~6.5). In particular, the induced residues are one less than
the residues of the connection induced by $Q$
for all characteristic directions. So the assumption on the residues 
implies that for no set of characteristic
directions the sum of the induced residues has real part equal to~$-1$; therefore
Theorem~\qqua\ implies that either (iii) holds or $[\gamma(t)]$ tends to a characteristic
direction $[v_0]$, whose residue cannot be purely imaginary.
If the real part of the residue of $[v_0]$ is positive, then we apply  
Corollary~\stcin.(ii), showing in particular that this can happen
only for a nowhere dense set of initial conditions.
If instead the real part of the residue
of $[v_0]$ is negative, we apply Corollary~\stcin.(i).(a) to finish the proof.\qedn

{\sc Example 8.2} By Theorem~\qqua, if for no set of characteristic directions the real
part of the sum or the induced residues belongs to $(-3/2,-1/2)$ then
the case (b.iii) cannot occur, and thus we get a complete description of
the dynamics of $Q$. For instance, assume that $Q$ is a quadratic field
(that is $\nu=1$) with three characteristic directions, necessarily of
order~1 (see Remark~5.4), with residues $\rho_1$, $\rho_2$ and~$\rho_3$
respectively. By the classical residue theorem for 
meromorphic connections (see, e.g., [IY, Theorem~I\negthinspace 
I\negthinspace I.17.33]) we know that
$$
\rho_1+\rho_2+\rho_3=1\;.
$$
Then an easy computation shows that case (b.iii) cannot occur if
$$
\Re\rho_1,\ \Re\rho_2\notin\left(-{1\over 2},{1\over 2}\right)\qquad
\hbox{and}\qquad \Re(\rho_1+\rho_2)\notin\left({1\over2},{3\over 2}\right)\;.
$$
As a consequence, at least one (and at most two) of the three residues must have negative real part, and thus case (b.i) can actually occur. In
particular, if only one residue has negative real part then almost all
integral curves converge to the origin tangentially to that characteristic
direction.

\bigskip
{\sc Example 8.3} The vector field
$$
Q=\left(-{1\over 3}(w^1)^2+{2\over3}w^1w^2\right){\de\over\de w^1}+
\left({2\over3}w^1w^2-{1\over3}(w^2)^2\right){\de\over\de w^2}
$$
has three Fuchsian characteristic directions ($[1:0]$, $[0:1]$ and $[1:1]$) 
of order $1$ and residue $1/3$. In particular, Corollary~\stcin\ says
that for most integral curves~$\gamma$ the geodesic $[\gamma]$ does not converge to a 
characteristic direction. Since no sum of induced residues has real
part equal to~$-1$, Theorem~\qqua\ implies that for most integral curves
the induced geodesic must intersect itself infinitely many times.
For instance, Fig.~3 shows $[\gamma]$ for the integral curve~$\gamma$
issuing from $(i,i-1)$. 

\Figuraeps 3 (geod2)

\smallsect 9. Quadratic vector fields

In this section we shall present, as an example of application of our methods, what we can infer on the dynamics of homogeneous quadratic vector fields.


Arguing as in [A3], it is not difficult to see that any not identically zero homogeneous quadratic vector field is linearly conjugated to one of the following:
\smallskip
\itemitem{$(\infty)$} $Q(z,w)=z^2{\de\over\de z}+zw{\de\over\de w}$;
\itemitem{$(1_{00})$} $Q(z,w)=-z^2{\de\over\de w}$;
\itemitem{$(1_{10})$} $Q(z,w)=-z^2{\de\over\de z}-(z^2+zw){\de\over\de w}$;
\itemitem{$(1_{11})$} $Q(z,w)=-zw{\de\over\de z}-(z^2+w^2){\de\over\de w}$;
\itemitem{$(2_{001})$} $Q(z,w)=zw{\de\over\de w}$;
\itemitem{$(2_{011})$} $Q(z,w)=zw{\de\over\de z}+(zw+w^2){\de\over\de w}$;
\itemitem{$(2_{10\rho})$} $Q(z,w)=-\rho z^2{\de\over\de z}+(1-\rho) zw{\de\over\de w}$, with $\rho\ne0$;
\itemitem{$(2_{11\rho})$} $Q(z,w)=(\rho z^2+zw){\de\over\de z}+\bigl((1+\rho) zw+w^2\bigr){\de\over\de w}$, with $\rho\ne0$;
\itemitem{$(3_{100})$} $Q(z,w)=(z^2-zw){\de\over\de z}$;
\itemitem{$(3_{\rho10})$} $Q(z,w)=\rho(-z^2+zw){\de\over\de z}+(1-\rho)(zw- w^2){\de\over\de w}$, with
$\rho\ne0$,~$1$;
\itemitem{$(3_{\rho\tau1})$} $Q(z,w)=\bigl(-\rho z^2+(1-\tau)zw\bigr)
{\de\over\de z}+\bigl((1-\rho)
zw-\tau w^2\bigr){\de\over\de w}$, with $\rho$,~$\tau\ne0$ and $\rho+\tau\ne1$.
\smallskip
\noindent In this list, the vector field $(\infty)$ is dicritical; the
vector fields $(1_{\bullet\bullet})$ have exactly one characteristic direction
$[0:1]$; the
vector fields $(2_{\bullet\bullet\bullet})$ have exactly two characteristic directions
$[1:0]$ and $[0:1]$; and the
vector fields $(3_{\bullet\bullet\bullet})$ have exactly three characteristic directions
$[1:0]$, $[1:1]$ and $[0:1]$. Let us see what we can say in the various cases. In our
computations, we shall often use Remark~1.7.
\medskip
\noindent {\bf Case $\hbox{\bf(}\bfm\infty\hbox{\bf )}$.} The field 
$$
Q(z,w)=z^2{\de\over\de z}+zw{\de\over\de w}
$$
is dicritical; in particular all complex lines $L_v$ with $v\in\C^2\setminus{O}$ are invariant, and the dynamics on each $L_v$ is described by
Lemma~\dcl. So the line $L_{(0,1)}$ is pointwise fixed, whereas on any
other line the integral curves of $Q$ goes to~$O$ in both forward 
and backward time, with the exception of one integral curve going to~$O$ in forward time and diverging to infinity in backward time, and of one integral curve diverging to infinity in forward time and going to~$O$ in backward
time.

\medbreak
\noindent {\bf Case $\hbox{\bf(1}_{\bf 00}\hbox{\bf )}$.} The field 
$$
Q(z,w)=-z^2{\de\over\de w}
$$
has only one (degenerate) characteristic direction $[0:1]$, necessarily of order~3,
index~$-1$, residue~$\rho=1$ and induced residue~$-2$. Then Theorem~\qqua\ and
Theorem~\ctreb\ imply that the direction~$[\gamma]$ of an integral
curve~$\gamma$ of~$Q$ outside the characteristic leaf~$L_{(0,1)}$ goes to
$[0:1]$ in both forward and backward time, whereas the characteristic leaf
is pointwise fixed. 

In the canonical coordinates $(\zeta_\infty,v_\infty)$ centered at~$[0:1]$
the geodesic field~$G$ is given by
$$
G=\zeta_\infty^3 v_\infty{\de\over\de\zeta_\infty}-\zeta^2_\infty 
v_\infty^2{\de\over\de v_\infty}\;;
$$
therefore $G$ is already in normal form with a Fuchsian singularity of order~3 (and vanishing resonant index). Furthermore $\mu_Y=2>1=\rho$;
hence Corollary~\stcin.(i.b) applies, and it follows that the norm $\|\gamma\|$
of any integral curve outside the characteristic leaf goes to~$+\infty$ 
in forward and backward time.  

In this particular case it is easy to explicitely write down the 
integral curves. Indeed, the integral curve issuing from~$(z_0,w_0)$ is given by
$$
\gamma(t)=(z_0,w_0-z_0^2 t)\;,
$$
whose behavior is exactly as predicted.

\medbreak
\noindent {\bf Case $\hbox{\bf(1}_{\bf 10}\hbox{\bf )}$.} The field 
$$
Q(z,w)=-z^2{\de\over\de z}-(z^2+zw){\de\over\de w}
$$
has only one (degenerate) characteristic direction $[0:1]$, necessarily of order~3,
index~$-1$, residue~$\rho=1$ and induced residue~$-2$. Then Theorem~\qqua\ and
Theorem~\ctreb\ again imply that the direction~$[\gamma]$ of an integral
curve~$\gamma$ of~$Q$ outside the characteristic leaf~$L_{(0,1)}$ goes to
$[0:1]$ in both forward and backward time, whereas the characteristic leaf
is again pointwise fixed. 

In the canonical coordinates $(\zeta_\infty,v_\infty)$ centered at~$[0:1]$
the geodesic field~$G$ is given by
$$
G=\zeta_\infty^3 v_\infty{\de\over\de\zeta_\infty}-\zeta_\infty(1+\zeta_\infty)
v_\infty^2{\de\over\de v_\infty}\;;
$$
therefore $G$ is in (formal) normal form with an irregular singularity of order~3, irregularity~2 and resonant index~1. We do not have (yet) 
general results about the dynamics nearby irregular singularities; however,
in this case we can work directly. Let $\sigma(t)=\bigl(\zeta_\infty(t),v_\infty(t)
\bigr)$ be an integral curve of~$G$. We know that $\sigma(t)$ is
contained in the horizontal foliation, and a quick computation shows that in  this case the horizontal foliation is given by
$$
\exp\left(-{1\over\zeta_\infty}\right)\zeta_\infty v_\infty\equiv c_0\;.
$$
It follows that $\sigma$ must satisfy the equation
$$
\zeta_\infty'=c_0\zeta_\infty^2\exp\left({1\over\zeta_\infty}\right)\;.
$$
Separating the variables we get $\exp\bigl(-1/\zeta_\infty(t)\bigr)=
c_0t+c_1$, that is
$$
\zeta_\infty(t)=-{1\over\log(c_0t+c_1)}\;, \qquad
v_\infty(t)=-{c_0\log(c_0 t+c_1)\over c_0t+c_1}\;,
$$
and $c_0$ and $c_1$ are determined by the initial conditions
$$
\log c_1=-{1\over\zeta_\infty(0)}\;,\qquad c_0=c_1 \zeta_\infty(0)
v_\infty(0)\;.
$$
It follows that if $\zeta_\infty(0)v_\infty(0)\notin\R$ then $v_\infty(t)$
tends to~0 for $t\to\pm\infty$; if instead $\zeta_\infty(0)v_\infty(0)\in
\R^*$ then $|v_\infty(t)|$ diverges to~$+\infty$ in finite time on one side 
and converges to~$0$ on the other side. Recalling that $\chi_1^{-1}
(\zeta_\infty,v_\infty)=(\zeta_\infty v_\infty,v_\infty)$, it follows that
if $(z_0,w_0)\in\C^2\setminus\{O\}$ is such that $z_0\notin\R$
then the integral curve~$\gamma$ of~$Q$ issuing from~$(z_0,w_0)$ goes to
the origin tangent to~$[0:1]$ in forward and backward time; if
instead $z_0\in\R^*$ then $\gamma(t)$ goes to the origin in forward time
and diverges in backward time, or conversely. 

In particular, the origin has an open basin of attraction even though 
the index of the characteristic direction has negative real part; this
cannot happen in the Fuchsian case (see Remark~8.3).

\medbreak
\noindent {\bf Case $\hbox{\bf(1}_{\bf 11}\hbox{\bf )}$.} The field 
$$
Q(z,w)=-zw{\de\over\de z}-(z^2+w^2){\de\over\de w}
$$
has only one (non-degenerate) characteristic direction $[0:1]$, necessarily of order~3,
index~$-1$, residue~$\rho=1$ and induced residue~$-2$. Then Theorem~\qqua\ and
Theorem~\ctreb\ still imply that the direction~$[\gamma]$ of an integral
curve~$\gamma$ of~$Q$ outside the characteristic leaf~$L_{(0,1)}$ goes to
$[0:1]$ in both forward and backward time; however this time the characteristic leaf is not pointwise fixed, and the dynamics there is
described by Lemma~\dcl.

In the canonical coordinates $(\zeta_\infty,v_\infty)$ centered at~$[0:1]$
the geodesic field~$G$ is given by
$$
G=\zeta_\infty^3 v_\infty{\de\over\de\zeta_\infty}-(1+\zeta_\infty^2)
v_\infty^2{\de\over\de v_\infty}\;;
$$
therefore $G$ is in (formal) normal form with an irregular singularity of order~3, irregularity~3 and resonant index~1. Let $\sigma(t)=\bigl(\zeta_\infty(t),v_\infty(t)
\bigr)$ be an integral curve of~$G$. We know that $\sigma(t)$ is
contained in the horizontal foliation, and a quick computation shows that in  this case the horizontal foliation is given by
$$
\exp\left(-{1\over2\zeta_\infty^2}\right)\zeta_\infty v_\infty\equiv c_0\;.
$$
It follows that $\sigma$ must satisfy the equation
$$
\zeta_\infty'=c_0\zeta_\infty^2\exp\left({1\over2\zeta_\infty^2}\right)\;.
$$
If $F$ denotes the primitive of the function $\exp(-w^2/2)$ with $F(0)=0$, separating
the variables we get
$$
F\left(-{1\over\zeta_\infty(t)}\right)=c_0 t+c_1\;.
$$
Since $F'$ is never vanishing, we can find a well-defined branch of $F^{-1}$
along the line $t\mapsto c_0t+c_1$, and thus
$$
\zeta_\infty(t)=-{1\over F^{-1}(c_0t+c_1)}\;,\qquad
v_\infty(t)=-c_0 F^{-1}(c_0 t+c_1)\exp\left({1\over 2}F^{-1}(c_0 t+c_1)^2
\right)\;,
$$
where $c_0$ and $c_1$ are determined by the initial conditions
$$
F^{-1}(c_1)=-{1\over\zeta_\infty(0)}\;,\qquad c_0=\exp\left(
-{1\over 2\zeta_\infty(0)^2}\right) \zeta_\infty(0)
v_\infty(0)\;.
$$
Writing $F^{-1}(c_0t+c_1)=R(t)+i I(t)$, with $R(t)$, $I(t)\in\R$, we have
$$
|v(t)|^2=|c_0|^2[R(t)^2+I(t)^2]\exp[R(t)^2-I(t)^2]\;.
$$
We know that $\zeta_\infty(t)\to 0$ as $t\to\pm\infty$ (if $c_1/c_0\in\R$
there is one value of~$t$ in which $\zeta_\infty(t)$ is not defined,
because the geodesic has left the canonical chart, but beyond that point
the geodesic re-enters the canonical chart); it follows that
$|F^{-1}(c_0t+c_1)|\to+\infty$. Furthermore, if $\limsup\limits_{t\to\pm\infty}
|R(t)|/|I(t)|<1$ then $|v(t)|\to 0$. Numerical experiments seem to suggest
that this can happen for an open set of initial conditions; if this is
correct, then we would have an open set of integral curves of $Q$ 
converging to the origin tangent to~$[0:1]$.

\medbreak
\noindent {\bf Case $\hbox{\bf(2}_{\bf 001}\hbox{\bf )}$.} The field 
$$
Q(z,w)=zw{\de\over\de w}
$$
has two (both degenerate) characteristic directions $[1:0]$ and~$[0:1]$. 
In the canonical coordinates $(\zeta_0,v_0)$ centered at~$[1:0]$
the geodesic field~$G$ is given by
$$
G=\zeta_0 v_0{\de\over\de\zeta_0}\;;
$$
therefore $[1:0]$ is an apparent singularity of order~1, residue~$\rho=0$,
induced residue~$-1$, and $G$ is in normal form. In the canonical coordinates $(\zeta_\infty,v_\infty)$ centered at~$[0:1]$
the geodesic field~$G$ is given by
$$
G=-\zeta_\infty^2 v_\infty{\de\over\de\zeta_\infty}+\zeta_\infty
v_\infty^2{\de\over\de v_\infty}\;;
$$
thus $[0:1]$ is a Fuchsian singularity of order~2, residue~$\rho=1$,
induced residue~$-1$, and $G$ is in normal form (up to a change of
sign in $v_\infty$). Both characteristic leaves are pointwise fixed.

This time Theorem~\qqua\ and Theorem~\ctreb\ say that the direction~$[\gamma]$ of an integral
curve~$\gamma$ of~$Q$ outside the characteristic leaves~$L_{(1,0)}$ 
and~$L_{(0,1)}$ either goes to $[1:0]$ or
$[0:1]$ in both forward and backward time, or is a closed geodesic.

In the coordinate chart centered at $[1:0]$ we can apply Corollary~\sdue; it follows that the geodesics are either a saddle connection between $[1:0]$ and $[0:1]$ or periodic; furthermore, $v_0(t)$ is constant. Recalling that
$\chi_1^{-1}(\zeta_0,v_0)=(v_0,\zeta_0v_0)$, it follows that the first
coordinate of an integral curve~$\gamma$ is always constant, and
the second coordinate is either (constant or) periodic or goes to~0 on
one side and diverges to infinity on the other side. So we have periodic integral curves; and the non-periodic integral curves go from a point
in the characteristic leaf~$L_{(1,0)}$ off to infinity toward the other
characteristic leaf~$L_{(0,1)}$; in particular, no integral curve
converges to the origin, and we have periodic integral curves (and hence periodic points for the time-1 map) accumulating at the origin. These results are confirmed by the explicit expression of the integral curve~$\gamma$ issuing from $(z_0,w_0)$:
$$
\gamma(t)=\bigl(z_0,w_0 e^{z_0 t}\bigr)\;.
$$

\medbreak
\noindent {\bf Case $\hbox{\bf(2}_{\bf 011}\hbox{\bf )}$.} The field 
$$
Q(z,w)=zw{\de\over\de z}+(zw+w^2){\de\over\de w}
$$
has two characteristic directions: $[1:0]$ is degenerate 
whereas~$[0:1]$ is non-degenerate. 
In the canonical coordinates $(\zeta_0,v_0)$ centered at~$[1:0]$
the geodesic field~$G$ is given by
$$
G=\zeta_0 v_0{\de\over\de\zeta_0}+\zeta_0v_0^2{\de\over\de v_0}\;;
$$
therefore $[1:0]$ is an apparent singularity of order~1, residue~$\rho=0$,
induced residue~$-1$, but $G$ is not in normal form. In the canonical coordinates $(\zeta_\infty,v_\infty)$ centered at~$[0:1]$
the geodesic field~$G$ is given by
$$
G=-\zeta_\infty^2 v_\infty{\de\over\de\zeta_\infty}+(1+\zeta_\infty)
v_\infty^2{\de\over\de v_\infty}\;;
$$
thus $[0:1]$ is an irregular singularity of order~2, irregularity~2, residue~$\rho=1$,
induced residue~$-1$, resonant index $1$, and $G$ is in normal form (up to a change of sign in $v_\infty$). The characteristic leaf~$L_{(1,0)}$ is pointwise fixed, while the dynamics on the characteristic leaf~$L_{(0,1)}$
is described by Lemma~\dcl. 

Theorem~\qqua\ and Theorem~\ctreb\ say again that the direction~$[\gamma]$ of an integral
curve~$\gamma$ of~$Q$ off the characteristic leaves~$L_{(1,0)}$ 
and~$L_{(0,1)}$ either goes to $[1:0]$ or
$[0:1]$ in both forward and backward time, or is a closed geodesic. In particular, self-intersecting geodesics are necessarily closed.

The description of the behavior of the integral curves of~$G$ nearby~$[1:0]$
is provided by Corollary~\sdue: we have periodic integral curves, 
integral curves converging to a non-zero element of the fiber over the
singularity, and integral curves escaping toward~$[0:1]$. In the chart
centered at~$[0:1]$, an integral curve $\sigma(t)=\bigl(\zeta_\infty(t),
v_\infty(t)\bigr)$ of $G$ satisfies 
$$
\exp\left(-{1\over\zeta_\infty(t)}\right)\zeta_\infty(t)v_\infty(t)\equiv
c_0\quad\hbox{and}\quad \zeta'_\infty=c_0\zeta_\infty e^{1/\zeta_\infty}\;;
$$
we leave to the reader an analysis similar to the one we did in case $(1_{11})$. 

\medbreak
\noindent {\bf Case $\hbox{\bf(2}_{\bf 10\rho}\hbox{\bf )}$.} The field 
$$
Q(z,w)=-\rho z^2{\de\over\de z}+(1-\rho) zw{\de\over\de w}
$$
with $\rho\ne 0$
has two characteristic directions: $[1:0]$ is non-degenerate 
whereas~$[0:1]$ is degenerate. 
In the canonical coordinates $(\zeta_0,v_0)$ centered at~$[1:0]$
the geodesic field~$G$ is given by
$$
G=\zeta_0 v_0{\de\over\de\zeta_0}-\rho v_0^2{\de\over\de v_0}\;;
$$
therefore $[1:0]$ is a Fuchsian singularity of order~1, residue~$\rho\ne0$,
induced residue~$\rho-1$, and $G$ is in normal form. In the canonical coordinates $(\zeta_\infty,v_\infty)$ centered at~$[0:1]$
the geodesic field~$G$ is given by
$$
G=-\zeta_\infty^2 v_\infty{\de\over\de\zeta_\infty}+(1-\rho)\zeta_\infty
v_\infty^2{\de\over\de v_\infty}\;;
$$
thus $[0:1]$ is a Fuchsian singularity of order~2, residue~$1-\rho$,
induced residue~$-1-\rho$, and $G$ is in normal form (up to a change of sign in $v_\infty$). The characteristic leaf~$L_{(0,1)}$ is pointwise fixed, while the dynamics on the characteristic leaf~$L_{(1,0)}$
is described by Lemma~\dcl. 

Since both singularities are already in normal form with vanishing resonant 
index, instead of using Theorem~\qqua\ we can rely directly on 
Proposition~\stqua\ and Corollary~\stcin\ to describe the behavior of the integral curves of $Q$. We have:
\smallskip
\item{--} if $\Re\rho<0$ then almost all integral curves of~$Q$ converge
to the origin tangent to $[1:0]$ both in forward and backward time;
each complex line~$L_v$ that is not a characteristic leaf contains exactly one real line of initial values of exceptional integral curves, which are converging to the origin tangent to~$[1:0]$ on one side and 
diverging to infinity toward~$L_{(0,1)}$ on the other side;
\item{--} if $\Re\rho>0$ then the roles of $[1:0]$ and $[0:1]$ are exchanged;
\item{--} if $\Re\rho=0$ then almost all integral curves of~$Q$ converge
to the origin both in forward and backward time without being tangent
to any direction; each complex line~$L_v$ that is not a characteristic leaf contains exactly one real line of initial values of exceptional integral curves, which cover closed geodesics and are converging to the origin on one side and diverging to infinity on the other side.

\medbreak
\noindent {\bf Case $\hbox{\bf(2}_{\bf 11\rho}\hbox{\bf )}$.} The field 
$$
Q(z,w)=(-\rho z^2+zw){\de\over\de z}+\bigl((1-\rho) zw+w^2\bigr){\de\over\de w}
$$
with $\rho\ne 0$
has two characteristic directions, both non-degenerate: $[1:0]$ and~$[0:1]$. 
In the canonical coordinates $(\zeta_0,v_0)$ centered at~$[1:0]$
the geodesic field~$G$ is given by
$$
G=\zeta_0 v_0{\de\over\de\zeta_0}-(\rho-\zeta_0) v_0^2{\de\over\de v_0}\;;
$$
therefore $[1:0]$ is a Fuchsian singularity of order~1, residue~$\rho\ne0$,
induced residue~$\rho-1$, and $G$ is not in normal form (unless $\rho=-1$, when $G$ is in normal form with resonant index~1). In the canonical coordinates $(\zeta_\infty,v_\infty)$ centered at~$[0:1]$
the geodesic field~$G$ is given by
$$
G=-\zeta_\infty^2 v_\infty{\de\over\de\zeta_\infty}+\bigl(1+(1-\rho)\zeta_\infty\bigr)
v_\infty^2{\de\over\de v_\infty}\;;
$$
thus $[0:1]$ is an irregular singularity of order~2, irregularity~2, residue (and resonant index)~$1-\rho$,
induced residue~$-1-\rho$, and $G$ is in normal form (up to a change of sign in $v_\infty$). The dynamics on the characteristic leaves~$L_{(1,0)}$ 
and~$L_{(0,1)}$ is described by Lemma~\dcl. 

Theorem~\qqua\ says that we can have closed geodesics only if $\Re\rho=0$,
and that if $\Re\rho\notin(-1/2,1/2)$ then no geodesic is self-intersecting;
so if $\Re\rho\notin(-1/2,1/2)$ then necessarily all geodesics are
saddle connections. 

If $\Re\rho>0$, Corollary~\stcin.(ii) applies, and we see that for almost 
all integral curves~$\gamma$ of~$Q$ the direction~$[\gamma]$ is
escaping from~$[1:0]$ (the exceptional curves are escaping to infinity
toward $[1:0]$). In particular, if $\Re\rho>1/2$ it follows that
for almost all integral curves~$\gamma$ the direction~$[\gamma]$ is
going to~$[0:1]$ both in forward and backward time. 

If $\Re\rho<0$ and $\rho\ne-1$ then Corollary~\stcin.(v.a) applies, and almost all integral curves whose direction starts close enough to~$[1:0]$
converge to the origin tangent to~$[1:0]$ (and if the conjecture
mentioned in Remark~8.2 is true than this holds for $\rho=-1$ too).

If $\Re\rho=0$ then Corollary~\stcin.(iii) applies, and we have integral
curves going to the origin or escaping to infinity without being tangent to 
any direction.

To complete the picture of this case, one needs to understand what happens nearby the irregular singularity. We sketch the approach suggested in Remark~1.7. The horizontal foliation is given by
$$
e^{-1/\zeta_\infty}\zeta_\infty^{1-\rho}v_\infty\equiv c_0\;;
$$
and a geodesic $\zeta_\infty(t)$ must satisfy
$$
\zeta'_\infty=-c_0\zeta_\infty^{\rho+1}\exp(1/\zeta_\infty)\;.
$$
Separating the variables one gets
$$
\Gamma\left(\rho,{1\over\zeta_\infty(t)}\right)=-(c_0t+c_1)\;,
$$
where $\Gamma(\rho,w)$ is the incomplete Gamma function. So one is left with 
studying the behavior of the inverse of the incomplete Gamma function.

\medbreak
\noindent {\bf Case $\hbox{\bf(3}_{\bf 100}\hbox{\bf )}$.} The field 
$$
Q(z,w)=(z^2-zw){\de\over\de z}
$$
has three characteristic directions: $[1:0]$ is non-degenerate, whereas
$[1:1]$ and~$[0:1]$ are degenerate. 
In the canonical coordinates $(\zeta_0,v_0)$ centered at~$[1:0]$
the geodesic field~$G$ is given by
$$
G=\zeta_0(\zeta_0-1)v_0{\de\over\de\zeta_0}-(\zeta_0-1) v_0^2{\de\over\de v_0}\;;
$$
therefore $[1:0]$ is a Fuchsian singularity of order~1, residue~$1$,
induced residue~$0$, and $G$ is not in normal form; on the other hand,
$[1:1]$ is an apparent singularity of order~1, residue~0, induced residue~$-1$, and $G$ is not in normal form. In the canonical coordinates $(\zeta_\infty,v_\infty)$ centered at~$[0:1]$
the geodesic field~$G$ is given by
$$
G=\zeta_\infty(\zeta_\infty-1) v_\infty{\de\over\de\zeta_\infty}\;;
$$
thus $[0:1]$ is an apparent singularity of order~1, residue~0,
induced residue~$-1$, and $G$ is not in normal form. The characteristic leaves $L_{(1,1)}$ and $L_{(0,1)}$ are pointwise fixed, while the dynamics on the characteristic leaf~$L_{(1,0)}$ is described by Lemma~\dcl. 

Putting together Theorem~\qdue, Corollary~\sdue\ and Corollary~\stcin.(ii)
we see that almost all integral curves goes from a
non-zero element of~$L_{(1,1)}$ to a non-zero element of~$L_{(0,1)}$;
the exceptions are periodic integral curves surrounding~$L_{(1,1)}$,
and integral curves diverging to infinity toward~$[1:0]$. This description
is confirmed by the explicit expression of the integral curve of~$Q$
issuing from $(z_0,w_0)\notin L_{(1,0)}\cup L_{(1,1)}\cup L_{(0,1)}$ given by
$$
\gamma(t)=\left({z_0w_0\over z_0-(z_0-w_0)e^{w_0t}},w_0\right)\;.
$$

\medbreak
\noindent {\bf Case $\hbox{\bf(3}_{\bf \rho10}\hbox{\bf )}$.} The field 
$$
Q(z,w)=\rho (-z^2+zw){\de\over\de z}+ (1-\rho)(zw-w^2){\de\over\de w}
$$
with $\rho\ne 0$, $1$
has three characteristic directions: $[1:0]$ and $[0:1]$ are non-degenerate, whereas $[1:1]$ is degenerate. 
In the canonical coordinates $(\zeta_0,v_0)$ centered at~$[1:0]$
the geodesic field~$G$ is given by
$$
G=\zeta_0(1-\zeta_0)v_0{\de\over\de\zeta_0}-\rho(1-\zeta_0) v_0^2{\de\over\de v_0}\;;
$$
therefore $[1:0]$ is a Fuchsian singularity of order~1, residue~$\rho$,
induced residue~$\rho-1$, and $G$ is not in normal form; on the other hand,
$[1:1]$ is an apparent singularity of order~1, residue~0, induced residue~$-1$, and $G$ is not in normal form. In the canonical coordinates $(\zeta_\infty,v_\infty)$ centered at~$[0:1]$
the geodesic field~$G$ is given by
$$
G=\zeta_\infty(1-\zeta_\infty) v_\infty{\de\over\de\zeta_\infty}-
(1-\rho)(1-\zeta_\infty)v_\infty^2{\de\over\de v_\infty}\;;
$$
thus $[0:1]$ is a Fuchsian singularity of order~1, residue~$1-\rho$,
induced residue~$-\rho$, and $G$ is not in normal form. The characteristic leaf $L_{(1,1)}$ is pointwise fixed, while the dynamics on the characteristic leaves~$L_{(1,0)}$ and $L_{(0,1)}$ is described by Lemma~\dcl. 

If $\Re\rho\in(0,1)$ then Corollary~\sdue\ and Corollary~\stcin.(ii) say
that almost all integral curves connect non-zero elements of~$L_{(1,1)}$;
the exceptions diverge to infinity toward $[1:0]$ or~$[0:1]$, or are
periodic integral curves around $L_{(1,1)}$. The local dynamics around
the singularities allows to exclude geodesics accumulating closed
geodesics or simple cycles of saddle connections; there might exist geodesics self-intersecting infinitely many times, however.

If $\Re\rho<0$ (respectively, $\Re\rho>1$) then $[1:0]$ (respectively, $[0:1]$) becomes attracting. Again, the local dynamics around
the singularities allows to exclude geodesics accumulating closed
geodesics or simple cycles of saddle connections. If $\Re\rho\le-1/2$
(respectively, $\Re\rho\ge3/2$) then we cannot have geodesics self-intersecting infinitely many times. Therefore almost all integral curves
either converge to $O$ tangentially to $[1:0]$ (respectively, to
$[0:1]$) or converge to a non-zero element of~$L_{(1,1)}$; the 
exceptions either diverge to infinity toward $[0:1]$ (respectively,
$[1:0]$) or are periodic integral curves around $L_{(1,1)}$. If
instead $\Re\rho\in(-1/2,0)$ (respectively, $\Re\rho\in(1,3/2)$) then
we have the same description, but there might exist geodesics
self-intersecting infinitely many times. But every simple loop in such
a geodesic must surround $[1:0]$ (respectively, $[0:1]$), must have
the same external angle, and cannot get too close to $[1:0]$ or 
$[1:1]$, and so the existence of such a geodesic seems unlikely. 

\medbreak
\noindent {\bf Case $\hbox{\bf(3}_{\bf \rho\tau1}\hbox{\bf )}$.} The field 
$$
Q(z,w)=\bigl(-\rho z^2+(1-\tau)zw\bigr)
{\de\over\de z}+\bigl((1-\rho)
zw-\tau w^2\bigr){\de\over\de w}$$
with $\rho$, $\tau\ne 0$ and $\rho+\tau\ne1$
has three characteristic directions: $[1:0]$ and $[0:1]$ and $[1:1]$, all non-degenerate. 
In the canonical coordinates $(\zeta_0,v_0)$ centered at~$[1:0]$
the geodesic field~$G$ is given by
$$
G=\zeta_0(1-\zeta_0)v_0{\de\over\de\zeta_0}-\bigl(\rho-(1-\tau)\zeta_0\bigr) v_0^2{\de\over\de v_0}\;;
$$
therefore $[1:0]$ is a Fuchsian singularity of order~1, residue~$\rho$,
induced residue~$\rho-1$, and $G$ is not in normal form. Also
$[1:1]$ is a Fuchsian singularity of order~1, residue~$1-\rho-\tau$, induced residue~$-\rho-\tau$, and $G$ is not in normal form. In the canonical coordinates $(\zeta_\infty,v_\infty)$ centered at~$[0:1]$
the geodesic field~$G$ is given by
$$
G=\zeta_\infty(1-\zeta_\infty) v_\infty{\de\over\de\zeta_\infty}-
\bigl(\tau-(1-\rho)\zeta_\infty\bigr)v_\infty^2{\de\over\de v_\infty}\;;
$$
thus $[0:1]$ is a Fuchsian singularity of order~1, residue~$\tau$,
induced residue~$\tau-1$, and $G$ is not in normal form. The  dynamics on the characteristic leaves~$L_{(1,0)}$, $L_{(1,1)}$ and $L_{(0,1)}$ is described by Lemma~\dcl. 

Thanks to the symmetry of the situation, and recalling that the sum of the residues is~1, it suffices to consider the following cases:
\smallskip
\item{(a)} $\Re\rho$, $\Re\tau>0$ and $\Re\rho+\Re\tau<1$ (three residues with positive real part). In this case Theorem~\qqua\ and 
Corollary~\stcin.(ii) implies that almost all geodesics intersect
themselves infinitely many times (as in Example~8.3); the exceptions are saddle connections corresponding to integral curves escaping to infinity
in both forward and backward time.
\item{(b)} $\Re\rho<0$, $\Re\tau>0$ and $\Re\rho+\Re\tau<1$ (two residues with positive real part, one residue with negative real part). Then $[1:0]$ is attracting; this means that we have at least an open set of initial conditions whose integral curves converge to the origin tangent to $[1:0]$.
If no induced residue belongs to $(-3/2,-1/2)$, that is $\Re\rho\le-1/2$,
$\Re\tau\ge1/2$ and $\Re\rho+\Re\tau\le1/2$, then Theorem~\qqua\ implies that almost all integral curves converge to the origin tangent to $[1:0]$;
the exceptions diverge to infinity toward $[1:1]$ or $[1:0]$, and thus in this case Corollary~\stsei\ yields a complete description of the dynamics. If instead
there is at least one induced residue belonging to $(-3/2,-1/2)$ there
might also be geodesics intersecting themselves infinitely many times.
\item{(c)} $\Re\rho$, $\Re\tau<0$ (two residues with negative real part).
A similar description applies to this case. The only differences are: 
we have two attracting characteristic directions, $[1:0]$ and $[0:1]$; the condition excluding induced residues belonging to $(-3/2,-1/2)$ is
$\Re\rho$, $\Re\tau\le -1/2$; and the exceptional integral curves
diverge to infinity toward $[1:1]$.
\item{(d)} $\Re\rho=0$, $\Re\tau\in(0,1)$ (one purely imaginary residue, two residues with positive real part). In this case we always have exactly one induced residue whose real part belongs to $(-3/2,-1)\cup(-1,-1/2)$
except when $\Re\tau=1/2$. So Theorem~\qqua\ and Corollary~\stcin\ say
that we have integral curves converging to the origin without being tangent
to any direction, and some exceptional integral curve diverging to infinity toward $[1:1]$ or $[0:1]$ or without being tangent to any direction;
and we might have integral curves corresponding to geodesics intersecting
themselves infinitely many times (this case is excluded if $\Re\tau=1/2$). 
\item{(e)} $\Re\rho=0$, $\Re\tau>1$ (one purely imaginary residue, one residue with positive real part, one residue with negative real part).
In this case we have integral curves converging to the origin tangent to $[1:1]$ or without being tangent
to any direction, some exceptional integral curve diverging to infinity toward $[0:1]$ or without being tangent to any direction;
and, if $1<\Re\tau<3/2$, we might have integral curves corresponding to geodesics intersecting
themselves infinitely many times.
\item{(f)} $\Re\rho=0$, $\Re\tau=1$ (two purely imaginary residues). Finally, in this case almost all integral curves converge to the origin without being
tangent to any direction; the exceptional integral curves diverge to
infinity either toward $[0:1]$ or without being tangent to any direction.

\setref{ABT2}

\beginsection References

\art A1 M. Abate: Diagonalization of non-diagonalizable discrete holomorphic
dynamical systems! Amer. J. Math.! 122 2000 757-781

\art A2 M. Abate: The residual index and the dynamics of holomorphic maps
tangent to the identity! Duke Math. J.! 107 2001 173-207

\art A3 M. Abate: Holomorphic classification of $2$-dimensional quadratic maps tangent to
the identity! S\=urikaisekiken\-ky\=usho K\=oky\=uroku! 1447 2005 1-14

\pre A4 M. Abate: Local discrete holomorphic dynamical systems! Preprint,
arXiv:0903.3289. To appear in {\bf Holomorphic dynamics,} G. Gentili, J. Guenot, G. Patrizio eds., Lect. Notes in Math., Springer Verlag, Berlin! 2010

\art AT1 M. Abate, F. Tovena: Parabolic curves in $\C^3$! Abstr. Appl. 
Anal.! 2003 2003 275-294

\book AT2 M. Abate, F. Tovena: Curve e superfici! Springer Italia, Milano,
2007

\art ABT1 M. Abate, F. Bracci, F. Tovena: Index theorems for holomorphic
self-maps! Ann. of Math.! 159 2004 819-864

\art ABT2 M. Abate, F. Bracci, F. Tovena: Index theorems for holomorphic
maps and foliations! Indiana Univ. Math. J.! 57 2008 2999-3048

%
\art C C. Camacho: On the local structure of conformal mappings and holomorphic
vector fields! Ast\'e\-risque! 59--60 1978 83-94

\book dC M.P. do Carmo: Differential geometry of curves and surfaces!
Prentice-Hall, Englewood Cliffs, 1976

\book \'E1 J. \'Ecalle: Les fonctions r\'esurgentes. Tome
I: Les alg\`ebres de fonctions r\'esurgentes! Publ. Math. Orsay {\bf 81-05,}
Universit\'e de Paris-Sud, Orsay, 1981 

\book \'E2 J. \'Ecalle: Les fonctions r\'esurgentes. Tome
I\negthinspace I: Les fonctions r\'esurgentes appliqu\'ees \`a l'it\'eration!
Publ. Math. Orsay {\bf 81-06,} Universit\'e de Paris-Sud, Orsay, 1981 

\book \'E3 J. \'Ecalle: Les fonctions r\'esurgentes. Tome I\negthinspace
I\negthinspace I: L'\'equation du pont et la classification analytique des
objects locaux! Publ. Math. Orsay {\bf 85-05,} Universit\'e de Paris-Sud, 
Orsay, 1985 

\coll \'E4 J. \'Ecalle: Iteration and analytic classification of local
diffeomorphisms of $\C^\nu$! Iteration theory and its functional equations!
Lect. Notes in Math., 1163, Springer Verlag, Berlin, 1985, 41--48

\art H1 M. Hakim: Analytic transformations of $(\C^p,0)$ tangent to the
identity! Duke Math. J.! 92 1998 403-428

\pre H2 M. Hakim: Transformations tangent to the identity. Stable pieces
of manifolds! Preprint! Universit\'e de Paris-Sud, 1997

\book HK B. Hasselblatt, A. Katok: Introduction to the modern theory of
dynamical systems! Cambridge Univ. Press, Cambridge, 1995

\book IY Y. Ilyashenko, S. Yakovenko: Lectures on analytic differential
equations! Graduate Studies in Mathematics 86, American Mathematical
Society, Providence, RI, 2008

\book K S. Kobayashi: Differential geometry of complex vector bundles! 
Princeton University Press, Princeton, NJ, 1987

\book M B. Malgrange:
\'Equations diff\'erentielles \`a coefficients polynomiaux! Birkh\"auser, Boston, 1991
 
\book Mi J. Milnor: Dynamics in one complex variable! 
Princeton University Press, Princeton, 2006

\book R M. Rivi: Local behaviour of discrete dynamical systems! Ph.D. 
Thesis, Universit\`a di Firenze, 1999

\book S S. Segal: Nine introductions in complex analysis! North Holland,
Amsterdam, 1981

\art Sh A.A. Shcherbakov: Topological classification of germs of conformal
mappings with identity linear part! Moscow Univ. Math. Bull.! 37 1982 60-65

\bye